\documentclass[leqno,10pt]
{article}
\usepackage{amsmath, amsfonts, amssymb}
\usepackage{psfig}

 \def\firstpage{1}
 \def\sign{\mathrm {sign}}
\def\Fc{{\mathcal F}}
\def\Pc{{\mathcal P}}
\def\Cc{\mathcal{C}}
\def\a{\alpha}
\def\ora{\overrightarrow}
\def\ola{\overleftarrow}
\def\g{\gamma}
\def\d{\delta}
\newcommand{\up}{{}^}

\newcommand\hook{\mathbin{\hbox{\vrule height .5pt width 3.5pt depth 0pt
\vrule height 6pt width .5pt depth 0pt}}}

\newcommand{\qed}{\hfill {$\blacktriangleright$}}
\newcommand{\fW}{{\mathfrak W}}
\newcommand{\fA}{{\mathfrak A}}
\newcommand{\fP}{{\mathfrak P}}
\newcommand{\fK}{\mathfrak {K}}
\newcommand{\fF}{{\mathfrak F}}
\newcommand{\fD}{{\mathfrak D}}
\newcommand{\fL}{{\mathfrak L}}

\newcommand{\fE}{{\mathfrak E}}

\newcommand{\cont}{{\mathrm {cont}}}
\newcommand{\restr}{\mathrm {restr}}
 \newcommand{\bx}{{\mathbf x}}
\newcommand{\by}{{\mathbf y}}
\newcommand{\bz}{{\mathbf z}}
\newcommand{\bgg}{{\mathbf g}}
\newcommand{\bff}{{\mathbf f}}

\newcommand{\bo}{\boldsymbol{\omega}}
\newcommand{\bmu}{\boldsymbol{\mu}}

\newcommand{\bh}{{\mathbf h}}
\newcommand{\bF}{{\mathbf F}}
\newcommand{\bG}{{\mathbf G}}

\newcommand{\bH}{{\mathbf H}}
\newcommand{\bK}{{\mathbf K}}
\newcommand{\bk}{{\mathbf k}}

\newcommand{\bV}{ \mathbf {V}}
 \newcommand{\CZ}{{\mathrm {CZ}}}
 \newcommand{\SFT}{{\mathrm {SFT}}}
 \newcommand{\RSFT}{{\mathrm {RSFT}}}
\newcommand{\bPhi}{\boldsymbol{\Phi}}
\newcommand{\bPsi}{\boldsymbol{\Psi}}

\newcommand{\wt}{\widetilde}

\newcommand{\Mc}{\mathcal {M}}
\newcommand{\cM}{\mathcal {M}}
\newcommand{\cD}{\mathcal {D}}
\newcommand{\cS}{{\mathcal S}}
\newcommand{\cK}{{\mathcal K}}
\newcommand{\cP}{{\mathcal P}}
\newcommand{\cE}{{\mathcal E}}

\newcommand{\CR}{\mathcal{CR}}

\newcommand{\Ker}{{\mathrm {Ker}}}
\newcommand{\cV}{{\mathcal V}}
\newcommand{\cN}{{\mathcal N}}

 \newcommand{\dist}{{\mathrm {dist}}}
\newcommand{\Coker}{{\mathrm {Coker}}}
\newcommand{\Morse}{{\mathrm {Morse}}}
\newcommand{\Int}{{\mathrm {Int}}}
\newcommand{\Id}{{\mathrm {Id}}}
\newcommand{\ind}{\mathrm{ind\,}}

\newcommand{\im}{\mathrm {Im}}
\newcommand{\comp}{\mathrm {comp}}
\newcommand{\D}{{\rm D}}

\newcommand{\contr}{{\mathrm {contr}}}
\newcommand{\punct}{{\mathrm {punct}}}
\newcommand{\even}{{\mathrm {even}}}
\newcommand{\cl}{{\mathrm {closed}}}
 \newcommand{\proof}{{\sl Proof. }}

\newcommand{\Z}{{\mathbb Z}}
\newcommand{\R}{{\mathbb R}}
\newcommand{\Q}{\mathbb {Q}\,}
\newcommand{\C}{\mathbb {C}\,}
\newcommand{\const}{{\mathrm const}}

\newcommand{\p}{{\partial}}

 \newcommand{\DM}{{\overline{\cM}}}
   \newcommand{\pr}{\mathrm{pr}}
 \newcommand{\ct}{\mathrm{ct}}
 \newcommand{\ev}{\operatorname{ev}}
 \newcommand{\disk}{\mathrm{disk}}
\newcommand{\gt}{{\tau}}
\newcommand{\go}{{\omega}}
\newcommand{\gd}{{\delta}}
\newcommand{\gb}{{\beta}}
\newcommand{\w}{{\wedge}}

 \newtheorem{theorem}{Theorem}[subsection]
\newtheorem{corollary}[theorem]{Corollary}
\newtheorem{lemma}[theorem]{Lemma}
\newtheorem{proposition}[theorem]{Proposition}

\newtheorem{problem}[theorem]{Problem}

\newtheorem{remark}[theorem]{Remark}
\newtheorem{definition}[theorem]{Definition}

\newtheorem{example}[theorem]{Example}

\begin{document}
\title{Introduction to Symplectic Field Theory}
\author{Yakov Eliashberg\\Stanford
University \and Alexander Givental \\  UC  Berkeley
\and Helmut Hofer\\New York University}
\date{October 2000}
 \maketitle

\begin{abstract}

We sketch in this article  a new theory, which we call {\sl
Symplectic Field Theory } or SFT,  which provides an approach to
Gromov-Witten invariants of symplectic manifolds and their
Lagrangian submanifolds in the spirit of topological field theory,
and at the same time
 serves as a rich source
of new invariants of contact manifolds and their Legendrian
submanifolds. Moreover, we  hope that the applications of  SFT
go far beyond this framework.\footnote{The research is   partially
supported by the National Science Foundation.}
\end{abstract}

 \tableofcontents
 \bigskip

 \centerline{\bf Disclosure}
 \medskip

Despite its length, the current paper presents only a very sketchy
overview of  Symplectic Field Theory. It contains practically no proofs,
and in a few places where the proofs  are given
their role is    just to
illustrate the involved ideas,  rather than to give  complete rigorous
arguments.

 The ideas, the algebraic formalism, and some of the applications
 of this new theory were presented  and popularized
 by the authors at several conferences
 and seminars (e.g. \cite{Eliash-ICM}). As a result,
   currently there exists a  significant
 mathematical community which is in some form familiar with
   the subject.
  Moreover, there are
  many mathematicians, including several former
   and current students of the authors,
  who are actively working  on foundational aspects
  of the theory and its applications, and    even
  published papers on this subject. Their results
    show that already the simplest versions of the theory
  have some remarkable corollaries (cf. \cite{Ustilovsky}).
   We hope that   the present  paper will help
   attracting even more people
   to   SFT.

   Of course, our ideas give just a small new twist
    to many other active directions of research in
    Mathematics and Physics
    (Symplectic topology, Gromov-Witten invariants and quantum cohomology,
    Floer homology theory, String theory,
     just to mention few), pioneered by  V.I. Arnold,
     C.~Conley-E.~Zehnder,
      M.~Gromov, S.K. Donaldson, E.~Witten, A.~Floer,
        M.~Kontsevich and others
     (see \cite{Arnold-conjecture, Conley-Zehnder,
     Gromov-holomorphic, Donaldson-invariants, Floer, Ruan, Witten1, Witten2, Kontsevich-CP2,
     Konts-Manin}).
     Many people independently contributed results and ideas,
      which may be considered as  parts of  SFT.
      Let us just mention here the work of  Yu.Chekanov
      \cite{Chekanov},
       K.Fukaya--K.Ono--Y.-G.Oh--H.Ohta  \cite{FKO3},
       A.~Gathmann  \cite{Gathmann}, E.Ionel--T.Parker
      \cite{Ionel-Parker, Ionel}  , Y.Ruan--A.-M.Li  \cite{Ruan-Li}.
      It also  draws on  other results of   the  current authors
      and their coauthors (see \cite{EHS, EH-3ball, Givental, Givental-Maslov,
      Givental-toric, Hofer-Weinstein, Hofer-Wysocki-Zehnder, AbbasHofer,
      HoferWysockiZehnder:94d, HoferWysockiZehnder:99}). The contact-geometric
      ingredient of our work is greatly motivated by two
      outstanding conjectures in contact geometry: Weinstein's conjecture about periodic
      orbits of Reeb fields  \cite{Weinstein},
       and  Arnold's chord conjecture \cite{Arnold-steps}.

   Presently, we are working on a series of papers devoted
   to the foundations,
   applications,  and further development of  SFT.  Among the applications,
   some of which are mentioned in this paper, are new invariants of contact manifolds
   and Legendrian knots and links, new methods for computing Gromov-Witten invariants,
   new restrictions on the topology of Lagrangian submanifolds, new non-squeezing type
   theorems in contact geometry etc. We are expecting new links with the low-dimensional topology
   and, possibly, Physics. It seems,
   however,
   that what we see at the moment is just a tip of an iceberg.
   The main   body of Symplectic Field Theory  and its applications is yet to be discovered.

     \bigskip
     \noindent{\bf Guide for an impatient reader.}
     The paper consists of two parts. The first part, except Section \ref{sec:Floer}
     and the end of
       Section \ref{sec:orientation},
     contains some  necessary   background symplectic-geometric and analytic information.
     An impatient reader can try to begin reading with Section
     \ref{sec:Floer},and use
     the rest of the first part for  the references.
     The second part begins with its own introduction
     (Section \ref{sec:informal})  where we present a very rough sketch of SFT. At the end of Section
         \ref{sec:informal} we describe the plan of the remainder of the paper.

   \bigskip
   \noindent{\bf Acknowledgements.} The authors benefited a lot from discussions
    with many mathematicians, and from the ideas which they corresponded to
    us. We are especially grateful to C. Abbas,  P. Biran, F. Bourgeois, K. Cieliebak,
    T. Coates, T. Ekholm, K. Fukaya, E. Getzler,
    M. Hutchings, E.-N.
    Ionel, V.M. Kharlamov,
    K. Mohnke, L. Polterovich, D. Salamon, M. Schwarz, K. Wysocki and E. Zehnder.
      A  part of the paper
    was written,  when  the first author visited
    RIMS at Kyoto University.
     He wants  to thank  RIMS  and  K. Fukaya and K. Saito,
    the organizers of a special program
    in Geometry and String Theory,   for the  hospitality. He also thanks
    T. Tsuboi and J.-L. Brylinski for organizing   cycle of lectures on  SFT at
     the    University of
    Tokyo and the Pennsylvania State University. The second author is thankful to
     A. Kirillov, R. Donagi, and
     the Department of Mathematics of the University
    of Pennsylvania  for the possibility to present   SFT in a series of lectures.
     The third author acknowledges the hospitality of   FIM at the ETH Z\"urich, where   some
     of the work was
carried
out. Our special thanks to N.M. Mishachev for drawing the
    pictures and to J. Sabloff for
    the computer verification of the formula
    (\ref{eq:CP2-final}).

\section{Symplectic  and analytic setup}

\subsection{Contact preliminaries}\label{sec:prelim}
A $1$-form $\a$ on a $(2n-1)$-dimensional manifold $V$
is called {\it contact} if the
restriction of $d\a$ to  the    $(2n-2)$-dimensional tangent
distribution $\xi =
\{\a = 0\}$ is non-degenerate (and hence symplectic).
A  codimension $1$ tangent distribution $\xi$ on $V$ is called
 a {\it contact structure} if
it can be locally (and in the co-orientable case globally) defined
by  the Pfaffian equation $\a
= 0$ for some choice of a contact form $\a$.  The pair  $(V,\xi)$   is called
a {\it contact manifold}.  According to  Frobenius' theorem the contact
condition is a condition of maximal
 non-integrability of the tangent hyperplane field $\xi$. In particular, all integral
submanifolds of $\xi$ have
dimension $\le n-1$. On the other hand, $(n-1)$-dimensional integral submanifolds, called
{\it Legendrian}, always exist in abundance.    We will be dealing in this paper
only with co-orientable, and moreover co-oriented contact structures.
Any non-coorientable contact structure can be canonically double-covered
 by a coorientable one.
If a contact form $\a$ is fixed then one can associate with it the {\it Reeb vector field}
$R_\a$, which is transversal
to the contact structure $\xi=\{\alpha=0\}$.
The field $R_\a$ is uniquely determined by the equations
$R_\a\hook d\a=0;\,\a(R_\a)=1\,.$   The flow of $R_\a$ preserves the contact form $\a$.

 The  $2n$-dimensional manifold $M=(T(V)/\xi)^*  \setminus V$ ,  called the
{\it symplectization} of $(V,\xi)$, carries a natural
symplectic structure $\omega$ induced by an embedding $M\to  T^*(V)$
 which assigns to each linear form $T(V) / \xi \to \R$ the
corresponding form $T(V) \to T(V)/\xi \to \R$. A choice of a contact
form $\a$ (if $\xi$ is co-orientable) defines a splitting $M=V\times(\R\setminus 0)$. As     $\xi$ is
 assumed  to be co-oriented  we can
 pick
 the positive half  $V\times\R_+$ of  $M$, and call it symplectization
as well.  The symplectic structure $\omega$ can be written in terms of this splitting as
$d(\tau\a), \tau>0$. It will be more convenient for us, however, to use additive notation
and write $\omega$ as $d(e^t\a),\,t\in\R$, on $M=V\times\R$.
Notice that the vector field $T=\frac{\p}{\p t}$ is conformally symplectic:
 we have ${\cal L}_T\omega=\omega$, as well as ${\cal L}_T
 (e^t\a)=e^t\a$,
 where ${\cal L}_T$ denotes the Lie derivative along the vector
 field $T$.
All the notions of contact geometry can be formulated as the
corresponding symplectic notions, invariant or equivariant with
respect to this conformal action. For instance, any contact
diffeomorphism of $V$ lifts to an equivariant symplectomorphism of
$M$;  contact vector fields on $V$ (i.e. vector fields preserving
the contact structure) are   projections of $\R$-invariant
symplectic (and automatically Hamiltonian) vector fields on $M$;
Legendrian submanifolds in $M$ correspond to cylindrical (i.e.
invariant with respect to the $\R$-action) Lagrangian submanifolds
of $M$.

 Notice that the Hamiltonian vector field on $V\times\R$, defined by the Hamiltonian function
$H=e^t$ is invariant under translations $t\mapsto t+c$, and projects to the Reeb vector field
$R_\alpha$ under the projection $V\times\R\to \R$.

The symplectization of a contact manifold is an example of a
symplectic manifold with {\it cylindrical} (or rather conical)
ends. We mean by that    a possibly non-compact symplectic
manifold $(W,\omega)$ with    ends of the form
$E^+=V^+\times[0,\infty)$ and  $E^-=V^-\times(-\infty,0]$, such
that  $V^\pm$ are compact manifolds, and
$\omega|_{V^\pm}=d(e^t\alpha^\pm)$, where $\a^\pm$ are  contact
forms on $V^\pm$. In other words,  the ends  $E^\pm$ of
$(W,\omega)$ are symplectomorphic, respectively,
  to the positive or negative halves of the
symplectizations of    contact manifolds $(V^\pm,\xi^\pm=\{\a^\pm=0\})$.
We will consider   the splitting of the ends and the  the contact forms
$\a^\pm$ to be parts of the structure of a symplectic manifold with cylindrical ends.
We will also call
$(W,\omega)$ a {\it directed symplectic cobordism} between
 the contact
manifolds $(V^+,\xi^+)$ and $(V^-,\xi^-)$, and denote it  by $\ora{V^-V^+}$.

Sometimes we will have to consider the compact part
$W^0=W\setminus(\Int E^+\cup\Int E^-)$ of a directed symplectic
cobordism $\ora{V^-V^+}$.  If it is not clear from the context we
will refer to $W^0$   as a {\it compact}, and to $W$ as a {\it
completed} symplectic cobordism.

Let us point out that ``symplectic cobordism" {\it is not an
equivalence relation, but rather a partial order. }
Existence of a directed
symplectic cobordism $\ora{V^-V^+}$  does not imply
the existence of a directed symplectic cobordism   $\ora{V^+V^-}$,
even if one does not fix contact forms for the contact structures $\xi^\pm$. On the other
hand,
directed
symplectic cobordisms $\ora{V_0V_1}$ and $\ora{V_1V_2}$ can be glued, in an obvious
way, into a directed symplectic cobordism  $\ora{V_0V_2}=\ora{V_0V_1}
\circledcirc\ora{V_1V_2}$.

Contact  structures have no local invariants. Moreover, any contact form is locally
isomorphic to the
form $\alpha_0= dz-\sum\limits_1^{n-1}y_idx_i$ (Darboux' normal form).
The contact structure $\xi_0$ on $\R^{2n-1}$ given by the form $\a_0$ is called {\it
standard}.
 The standard contact structure on $S^{2n-1}$ is formed by complex tangent hyperplanes to the
unit sphere in $\C^n$. The standard contact structure on $S^{2n-1}$
is isomorphic in the complement of a point to the standard contact structure  on
$\R^{2n-1}$.
According to a theorem of J. Gray (see \cite{Gray}) contact structures
 on closed
manifolds   have  the following stability property:
{\it Given a family  $\xi_t$, $t\in[0,1]$, of contact  structures on  a closed manifold
$M$, there exists
an isotopy $f_t:M\to M$,
such that $df_t(\xi_0)=\xi_t; t\in [0,1]$.}
Notice that for contact {\it forms} the analogous statement is  wrong. For
instance,
the topology of the $1$-dimensional foliation determined by the  Reeb vector field $R_\a$
is very sensitive to deformations of the contact form $\a$.

The conformal class of the symplectic form $d\a|_\xi$ depends only on
the cooriented contact structure $\xi$ and not on the
choice of the contact form $\a$. In particular, one can associate with $\xi$
 an almost complex structure $J:\xi\to\xi$, compatible with $d\a$
which means that $d\alpha(X,JY); X,Y\in\xi,$ is an Hermitian metric on $\xi$.
The space of   almost complex structures $J$ with this property is contractible, and hence
 the choice of $J$ is homotopically canonical.
Thus a  co-oriented contact structure  $\xi$ defines on $M$ a {\it
 stable almost complex
structure}
 $\widetilde J=\widetilde J_\xi$,
i.e. a splitting of
the tangent bundle $T(V)$ into the Whitney sum of a complex bundle of (complex) dimension
$(n-1)$ and a trivial $1$-dimensional real bundle.
The existence of a
 stable almost complex
structure  is  necessary for  the existence of a contact structure on $V$. If $V$
is open (see \cite {Gro-PDR}) or ${\mathrm {dim}}\,V=3$ (see \cite{Martinet,Lutz})
 this property is also sufficient for  the existence of a contact structure in the
prescribed homotopy
class. It  is still unknown
whether this condition
is sufficient for  the existence of a contact structure on a closed manifold of dimension
$>3$. However,   a positive answer  to this question is extremely unlikely.
The homotopy class of $\widetilde J_\xi$, which we denote by
$[\xi]$ and call the {\it formal }homotopy class  of $\xi$,
serves as an
invariant of $\xi$.
For an open $V$ it is a complete invariant (see \cite {Gro-PDR})   up to   homotopy
of contact structures,
but not up to a contact diffeomorphism. For closed manifolds this is known to be false
in many, but not all dimensions. The theory discussed in this
paper serves as a rich source of contact invariants, both of
closed and open contact manifolds.

\subsection{Dynamics of Reeb vector fields}\label{sec:dynamics}

Let $(V,\xi)$ be a $(2n-1)$-dimensional manifold with  a
co-orientable contact structure with a fixed contact form $\alpha$.
 For a generic choice of $\a$ there are only countably many periodic trajectories
  of the
vector field $R_\a$.  Moreover, these trajectories can be assumed {\it non-degenerate}
in the sense  that the linearized Poincar\'e return map $A_\g$
 along any  closed
 trajectory $\g$, including multiples,
has no eigenvalues equal to $1$.  Let  us denote by $\Pc=\Pc_\a$
  the set of  all periodic trajectories of $R_\a$, including multiples.
  \footnote{As it is explained below in Section \ref{sec:orientation}
  the orientation issues require us to exclude certain multiple periodic orbits out
  of consideration. Namely, let us recall that real eigenvalues of symplectic matrices
  different from $\pm 1$ come in pairs $\lambda,\lambda^{-1}$.
 Let $\g\in \Pc$ be a simple periodic orbit and $A_\g$ its linearized Poincar\'e return map.
If
  the total multiplicity of   eigenvalues of $A_\g$ from the interval
  $(-1,0)$ is odd, then we   exclude from $\Pc$ all even
  multiples of $\g$.}

The reason for a such choice
is discussed in Section \ref{sec:orientation} below. We will also fix a point
$m_\gamma$ on each {\it simple} orbit from $\Pc$.
Non-degenerate trajectories can be divided into {\it odd} and {\it even} depending
on the sign of the Lefshetz number ${\mathrm {det}}  (I- A_\g)$.
Namely, we call $\g$ odd if $  {\mathrm {det}}  (I-  A_\g)<0$, and even otherwise.
The parity of  a periodic orbit $\g$  agrees   with the parity of a certain
integer grading  which is defined   if certain additional choices are made, as
it is described below.

If $H_1(V)=0$
 then for each $\g\in\Pc$ we can choose  and
 fix a surface $F_{\g}$ spanning
the trajectory $\g$ in $V$.
  We will allow the case $H_1(V)\neq 0$,
   but will require in most of the paper that the torsion part is
   trivial
   \footnote{ The case when $H_1(V)$ has torsion elements is
   discussed in Section \ref{sec:torsion} below.}.
In this case   we choose a basis of $H_1(V)$,
  represent    it by oriented curves ${C}_1,\dots,{C}_K$, and    choose
  a symplectic trivialization of the bundle $\xi|_{{C}_i}$ for each chosen curve.
  We recall that the bundle $\xi$ is endowed with the symplectic
  form $d\alpha$ whose   conformal class depends only on $\xi$.
       For  any periodic orbit $\gamma\in\Pc$
   let us choose a surface $F_\gamma$ with $[\partial F_\gamma]=
   [\gamma]-\sum n_i[{C}_i]$.
  The coefficients $n_i$ are uniquely defined because of our  assumption that $H_1(V)$
  is
  torsion-free.

The above choices enable us to define the {\it Conley-Zehnder
index} $\CZ(\g)$ of  $\g$ as follows. Choose a homotopically
unique trivialization of the symplectic vector bundle $(\xi,d\a)$
over each trajectory $\g\in\Pc$ which extends to $\xi|_{F_{\g}}$
(and coincides with a chosen trivialization of $\xi|_{{C}_i}$
if ${C}_i$ is not homologically trivial). The linearized flow
of $R_\a$ along $\g$ defines then a path in the group
$Sp(2n-2,\R)$ of symplectic matrices, which begins at the unit
matrix and ends at a matrix with all eigenvalues different from
$1$. The Maslov index of this path (see \cite{Arnold-Maslov,
Robbin-Salamon}) is, by  the definition, the Conley-Zehnder index
$\CZ(\g)$ of the trajectory $\g$. See also
 \cite{HoferWysockiZehnder:94d}, Section 3, for an axiomatic
        description of the Conley-Zehnder index using our normalization
        conventions.

Notice that by changing the spanning surfaces for the trajectories
from $\Pc$ one can change Conley-Zehnder indices by  the value of
the cohomology class $2c_1(\xi)$, where $c_1(\xi)$ is the first
Chern class of the contact bundle $\xi$. In particular, ${\mathrm
{mod}}\,2$ indices can be defined independently of any spanning
surfaces, and even in the case when $H_1(V)\neq 0$. In fact,
$$(-1)^{\CZ(\g)}=(-1)^{n-1}\sign \left({\mathrm
{det}}(I-A_\gamma)\right).$$

  \subsection{Splitting of a symplectic manifold along \\
  a contact submanifold}
\label{sec:splitting} Let $V$ be a hypersurface of contact type,
or in a different terminology, a symplectically convex
hypersurface  in a symplectic manifold $(W,\omega)$. This means
that $\omega$ is exact, $\omega=d\beta$, near $V$, and the
restriction $\alpha=\beta|_V$ is a contact form on $V$.
Equivalently, one can say that the conformally symplectic vector
field $X$, $\omega$-dual to $\beta$, is transversal to $V$. Let us
assume that $V$ divides $W,\, W=W_+\cup W_-$, where
 the notation of the parts are chosen in such a way that   $X$
serves as an inward transversal for $W_+$, and an outward
transversal for $W_-$. The manifolds $W_\pm$ can be viewed as
compact directed symplectic cobordisms such that $W_-$ has only positive
contact boundary $(V,\alpha)$, while the same contact manifold
serves as a negative boundary of $W_+$.

Let $$(W_-^\infty,\omega_-^\infty)=(W_-,\omega)\cup\left(V\times[0,\infty),d(e^t\alpha)\right)$$ and
$$(W_+^\infty,\omega_+^\infty)=\left(V\times(-\infty,0],d(e^t\alpha)\right)\cup(W_+,\omega)$$
 be the  completions, and
$$(W_-^\tau,\omega_-^\tau)=(W_-,\omega)\cup\left(V\times[0,\tau],d(e^t\alpha)\right)$$ and
$$(W_+^\tau,\omega_+^\tau)=\left(V\times[ -\tau,0],d(e^t\alpha)\right)\cup(W_+^0,\omega)$$
  {\it partial completions} of
$W_\pm$. Let us observe that the symplectic manifolds
$$(W_-,e^{-\tau}\omega_-^\tau),\quad (V\times [-\tau,\tau],
d(e^t\alpha))\quad\hbox{and}\quad (W_+,e^\tau\omega_+^\tau)$$ fit
together into a symplectic manifold $(W^\tau, \omega^\tau)$, so
that $W^0=W$. Hence when $\tau\to\infty$ the deformation
$(W^\tau,\omega^\tau)$ can be viewed as a decomposition of the
symplectic manifold $W=W^0$ into the union   of two completed
symplectic cobordisms  $W^\infty_+$ and $ W^\infty_-$ . We will
write $W=W_-\circledcirc W_+$ and also
$W^\infty=W_-^\infty\circledcirc W_+^\infty$.

\begin{figure}
\centerline{\psfig{figure=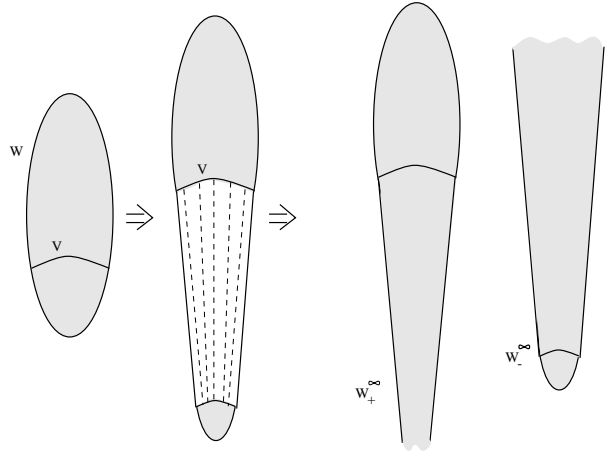,height=60mm}}
 \caption{\small Splitting of a closed symplectic manifold $W$
 into two completed symplectic cobordisms $W_-^\infty$ and $W^\infty_+$}
\label{fig:splitting}
\end{figure}
 Let us give here two important examples of the above
 splitting construction.

 \begin{example}\label{ex:cotangent}
 {\rm Suppose $L\subset W$ is a Lagrangian submanifold.
  Its neighborhood  is symplectomorphic to a
   neighborhood of the $0$-section  in the cotangent bundle
   $T^*(L)$. The boundary $V$ of an  appropriately chosen neighborhood
   has  contact type, and thus we can apply    along $V$
   the  above splitting construction. As the result we split
   $W$ into  $W_+^\infty$ symplectomorphic
   to $W\setminus L$, and $W_-^\infty$ symplectomorphic to $T^*(L)$.
    }
    \end{example}
    \begin{example}\label{ex:section}
        {\rm Let $M$ be a hyperplane section of a K\"ahler manifold
   $W$, or
   more generally a symplectic hyperplane section of a
   symplectic manifold
   $W$, in the sense of Donaldson (see \cite{Donaldson}). Then $M$ has a neighborhood
   with a contact boundary $V$. The affine part
   $W\setminus M$ is a Stein manifold in the K\"ahlerian case,
   and in any case
   has a structure of   a symplectic Weinstein
   manifold $\widetilde W$ (notice that the symplectic structure
   of $\widetilde W$ does not coincide with  the induced
   symplectic structure
   on $W\setminus M$ but contains $W\setminus M$ as an
   open symplectic
   submanifold).
   The Weinstein manifold $W\setminus M$ contains
   an isotropic  deformation retract  $\Delta$. The splitting of $W$ along $V$
   produces $W^\infty_-$ symplectomorphic to
    $\widetilde W,$ and $W^\infty_+$
    symplectomorphic to $W\setminus \Delta$.
    If $\Delta$ is a smooth Lagrangian submanifold, then
    we could get the same decomposition by
    splitting along the boundary of a tubular neighborhood
    of $L$, as in Example \ref{ex:cotangent}.  }
              \end{example}

\subsection{Compatible almost complex structures}\label{sec:almost}

According to M. Gromov (see \cite{Gromov-holomorphic}) an almost
complex structure $J$ is called {\it tamed} by a  symplectic form
$\omega$ if $\omega$ is positive on complex lines. If, in
addition,  one adds the calibrating condition that $\omega$ is
$J$-invariant, then $J$ is   said {\it compatible} with $\omega$.
For symplectic manifolds with cylindrical ends one needs further
compatibility conditions at infinity, as it is described below.

At each positive, or negative end $\left(V\times\R_\pm,d(e^t\alpha
)\right)$ we require   $J$ to be invariant with respect to
translations $t\mapsto t \pm c$, $c>0$ at least for sufficiently
large $t$. We also require   the contact structure
$\xi^\pm|_{V\times t}$ to be invariant under $J$, and define
$J\frac{\partial}{\partial t}=R_\alpha$, where
 $R_\alpha$ is
 the Reeb vector field (see \ref{sec:dynamics} above)
  of the contact form $\alpha$. In the case when
$W=V\times\R$ is the symplectization of a manifold $V$, i.e. $W$ is    a
cylindrical manifold, we additionally  require   $J$ to be globally
invariant under all translations along the second factor.

To define a compatible almost complex structure $J$ in the above
Examples \ref{ex:cotangent} and \ref{ex:section}
 one needs to specify a contact form $\alpha$ on the contact manifold
 $V$. In the case of the boundary of a tubular neighborhood    of a Lagrangian
 submanifold $L$ a natural choice of a contact form is provided by
 a  Riemannian metric on $L$. The Reeb vector field for such a
 form $\alpha$ generates on $V$ the geodesic flow of the metric.

 When $V$ is the boundary of a neighborhood of a hyperplane section  $M$
 then there exists another good choice of a contact form. It is
 a $S^1$-invariant connection form $\a$ on the principal $S^1$-bundle
 $V\to M$, whose curvature equals the symplectic form $\omega|_M$.
 The contact manifold $(V,\xi=\{\a=0\})$ is called the {\it
 pre-quantization} of the symplectic  manifold $(M,\omega)$.
 Orbits of the  Reeb field $R_\a$ are all closed and coincide
 with the fibers of the fibration, or their multiples. Notice that though
 the Reeb flow in
 this case looks extremely nice and simple, all its  periodic orbits are
 highly degenerate, see Section \ref{sec:Bott} below.
 Notice that  the symplectization $W$  of $V$
 can be viewed  as the total space of a complex  line bundle
 $L$ associated with
the $S^1$-fibration $V\to M$ with the zero-section removed.
It is possible
 and convenient to choose
$J$  compatible with the structure of this bundle,
and in such a way that the projection   $W\to M$
becomes holomorphic with respect to a certain almost
complex structure $J_M$ on $M$
compatible with $\omega$.

   Let us describe now what the symplectic splitting
   construction from  Section \ref{sec:splitting}
   looks like from the point of view of  a compatible almost
   complex structure.

   First, we assume that the original almost complex structure
   $J$ on $W$ is chosen in such a way
   that     the contact structure $\xi=\{\alpha=0\}$ on $V$ consists
   of complex tangencies to $V$, and that $JX=R_\alpha$, where
   $X$ is a conformally symplectic vector field, $\omega$-dual
   to $\alpha$, and $R_\alpha$ is the Reeb vector field of $\alpha$.
   Next we define an almost complex structure
   $J^\tau$ on $W^\tau=W_-\cup V\times[-\tau,\tau]\cup W_+$ by setting
   $J^\tau|_{W_\pm}=J$ and requiring  $J^\tau$ to be independent
   of $t\in[-\tau,\tau]$ on  $V\times[-\tau,\tau]$.
   When $\tau\to\infty$  the almost complex structure
   $J_-^\tau$ on
   $W_-^\tau=W_- \cup V\times[-\tau,\tau]$   converges
    to an almost complex
   structure $ J^\infty_-$ on $W^\infty_- $ compatible   with
   $\omega^\infty_-$, and  $J_+^\tau$ on
   $W_+^\tau=V\times[-\tau,\tau]\cup W_+$   converges
    to an almost complex
   structure $ J^\infty_+$ on $W^\infty_+$ compatible   with
   $\omega^\infty_+$.

   \subsection{Holomorphic curves in symplectic
   cobordisms}\label{sec:holomorphic}

  Let $(V,\alpha)$ be a contact manifold with a fixed contact form
and $(W=V\times\R,\omega=d(e^t\alpha))$ its symplectization. Let
us denote by $\pi_\R$ and $\pi_V$ the projections $W\to\R$ and $W\to
V$, respectively. For a map $f:X\to W$ we write $ f_\R$ and $f_V$
instead of $\pi_\R\circ f$ and $\pi_V\circ f$.

Notice that given a trajectory $\gamma$ of the Reeb field
$R_\a$,
the cylinder $\R\times\gamma\subset W$ is a $J$-holomorphic curve.
 Let us also observe that
\begin{proposition}\label{prop:positivity}
  For a $J$-holomorphic curve  $C\subset W$
 the restriction $d\a|_C$ is non-negative, and if $d\a|_C\equiv 0$
 then $C$ is a (part of a) cylinder $\R\times\gamma$ over a
 trajectory $\gamma$ of the Reeb field $R_\a$.
  \end{proposition}

 Given a $J$-holomorphic map $f$ of a punctured disk $D^2\setminus
 0\to W$ we say that  the map $f$ is {\it asymptotically
 cylindrical } over a periodic orbit $\g$ of the Reeb field
$R_\a$ at $+\infty$ (resp. at $-\infty$) if
$\mathop{\lim}\limits_{r\to 0} f_\R(re^{i\theta})=+\infty$ (resp.
$=-\infty$), and $\mathop{\lim}\limits_{r\to
0}f_V(re^{i\theta})=\bar f(\theta)$, where the map $\bar
f:[0,2\pi]\to V$ parameterizes the trajectory $\gamma$.

 \bigskip

 The almost complex manifold $(W,J)$ is bad from the point of view of
the theory of
holomorphic curves: it has a {\it pseudo-concave} end
$V\times(-\infty,0)$, or using Gromov's terminology its geometry
at this end is not bounded. However, it was shown in
\cite{Hofer-Weinstein} that   Gromov compactness theorem can be
modified to accommodate this situation, see Theorems
\ref{thm:comp1} and \ref{thm:comp2}
below. We will mention in this
section only the following fact related to compactness, which motivates the usage
of holomorphic curves asymptotically cylindrical over orbits from
$\Pc_\a$.
\begin{proposition}
 Suppose that all periodic orbits of the Reeb field $R_\alpha$ are
non-degenerate.
 Let  $C$ be  a non-compact Riemann surface  without boundary and $f:C\to W$
   a proper $J$-holomorphic curve. Suppose that there
exists a constant $K>0$ such that
$\int\limits_{C}f^*d\a  < K$.
Then
 $C$ is conformally equivalent to a compact Riemann surface $S_g$
 of genus $g$ with $s^++s^-$ punctures
 $$ x^+_1,\dots,
 x^+_{s^+},x^-_1\dots,x^-_{s^-}\in S_g,$$
  such that near the punctures ${\bx^+}=(x^+_1,\dots,x^+_{s^+})$
  the map
 $f$ is asymptotically cylindrical over periodic orbits
  ${\Gamma^+}=\{\g^+_1,\dots,\g^+_{s^+}\}$ at $+\infty$, and near
 the punctures ${\bx^-}=\{x^-_1,\dots,x^-_{s^-}\}$ the map
 $f$ is asymptotically cylindrical over periodic orbits
  ${\Gamma^-}=\{\g^-_1,\dots,\g^-_{s^-}\}$ at $-\infty$.
\end{proposition}

 Thus holomorphic maps of punctured Riemann surfaces,
   asymptotically cylindrical over periodic
orbits of the Reeb vector field $R_\alpha$, form a natural class of
holomorphic curves to consider in symplectizations
as well as  more general symplectic manifolds with cylindrical ends.  We
will define now moduli spaces of such curves.

Let $W =\ora{V^-V^+}$ be a (completed) directed cobordism,
$\alpha^\pm$    corresponding contact forms on $V^-$ and $V^+$,
 $\Pc^\pm$  the sets of all
periodic orbits (including multiple ones) of the Reeb vector
fields
$R_{\alpha^\pm}$. We assume that $\alpha^\pm$ satisfies
the genericity assumptions from Section \ref{sec:dynamics}.
Choose a compatible almost complex  structure  $J$   on $W$.
 Let $\Gamma^\pm$ be   ordered sets of trajectories from
 $\Pc^\pm$ of cardinality $s^\pm$.
 We also assume that every {\it simple} periodic orbit $\gamma$ from
 $\Pc^\pm$ comes with a fixed marker
 $m_\gamma\in\g$.

 Let $S=S_g$ be a compact Riemann surface of genus $g$ with a conformal
 structure $j$, with $s^+$ punctures $\bx^+=\{x^+_1,\dots,x^+_{s^+}\}$, called
 positive,
 $s^-$ punctures $\bx^-=\{x^-_1,\dots,x^-_{s^-}\}$, called negative, and
 $r$ marked points $\by=\{y_1,\dots,y_r\}$.
  We will also   fix an {\it asymptotic marker} at each
   puncture. We mean by that  a ray originating at each
   puncture.  Alternatively, if one
   takes the cylinder $S^1\times[0,\infty)$ as a conformal model of
   the punctured disk $D^2\setminus 0$ then an asymptotic marker
   can be viewed as a point on the circle at infinity.
   If a holomorphic map $f:D^2\setminus 0\to V^\pm\times\R_\pm$ is asymptotically
   cylindrical over a periodic orbit $\g$,  we say that
    a marker $\mu=\{\theta=\theta_0\}$
    is mapped by $f$ to the
   marker $m_\gamma\in\overline\gamma$, where $\overline\g$ is the simple
     orbit which  underlines $\g$, if $\mathop{\lim}\limits_{r\to
   0}f_{V^\pm}(re^{i\theta_0})=m_\g$.
Let us recall (see Section \ref{sec:dynamics} above) that we
provided each periodic orbit from $\Pc^\pm$ with a ``capping"
surface. This surface bounds $\g\in\Pc^\pm$ in $V_\pm$ if $\g$ is
homologically trivial, or realizes a homology    between $\g$ and
the corresponding linear combination of  basic curves
${C}_i^\pm$. We will  continue to rule out torsion elements in
the first homology
 (see the discussion of torsion in Section \ref{sec:torsion} below)   and choose
  curves ${C}_k\subset W$ which represent a basis of
the image $H_1(V^-\cup V^+)\to H(W)$ and for each curve
${C}^\pm_i$ fix a surface $G^\pm_i$ which realizes a homology
in $W$
between ${C}^\pm_i$ and the corresponding linear combination of
curves ${C}_k$.
All the choices enable us to associate with a  relative
 homology class $A'\in H_2(W,\Gamma^-\cup\Gamma^+)$,
 $\Gamma^\pm\subset\Pc^\pm$,
 an absolute integral  class
 $A\in H_2(W)$, which we will view as an element of $H_2(W;\C)$.

Let us denote by $\Mc^A_{g,r}(\Gamma^-,\Gamma^+;W,J)$ the moduli
space of $(j,J)$-holomorphic curves
$S_g\setminus(\bx^-\cup\bx^+)\to W$ with $r$ marked points, which
are asymptotically cylindrical over the periodic orbit $\g^+_i$
from $\Gamma^+$ at    the positive end at the puncture $x^+_i$,
and asymptotically cylindrical over the periodic orbit $\g^-_i$
from $\Gamma^-$ at    the negative end  at the puncture $x^-_i$,
and which send asymptotic markers to the
markers on the
corresponding periodic orbits. The curves from
$\Mc^A_{g,r}(\Gamma^-,\Gamma^+;W,J)$ are additionally required to
satisfy a stability condition, discussed in the next section. We
write $\Mc^A_{g}(\Gamma^-,\Gamma^+;W,J)$ instead of
$\Mc^A_{g,0}(\Gamma^-,\Gamma^+;W,J)$, and
$\Mc^A_{g,r}(\Gamma^-,\Gamma^+)$ instead of
$\Mc^A_{g,r}(\Gamma^-,\Gamma^+;W,J)$ if   it is clear from the
context which target almost complex  manifold $(W,J)$ is
considered.

 Notice, that we are not fixing $j$, and the configurations of
 punctures, marked points  or asymptotic markers.
 Two maps are called equivalent if they differ by a conformal map
 $S_g\to S_g$ which preserves all punctures, marked points and
 asymptotic markers.
     When the manifold $W=V\times\R$ is cylindrical,
     and hence the almost complex structure      $J$ is invariant under
     translations along the second factor, then
     all the moduli spaces $\Mc^A_{g,r}(\Gamma^-,\Gamma^+;W,J)$ inherit the $\R$-action.
     We will denote the quotient moduli space  by
       $\Mc^{A}_{g,r}(\Gamma^-,\Gamma^+;W,J)/\R$, and by
$\Mc^A_{g,r,s^-,s^+}(W,J)$ the union $\bigcup
\Mc^A_{g,r}(\Gamma^-,\Gamma^+;W,J)$ taken over
all sets of periodic orbits $\Gamma^\pm\subset\Pc^\pm$   with the
prescribed numbers $s^\pm$ of elements.
We will also need to consider the moduli space of disconnected
  curves of  Euler characteristic $2-2g$, denoted by
  $\wt{\Mc}_{g,r}^A(\Gamma^-,\Gamma^+)$.

\subsection{Compactification of the moduli spaces
$\Mc_{g,r}^A(\Gamma^-,\Gamma^+)$}\label{sec:compact}

To describe the compactification we need an appropriate notion of
a {\it stable holomorphic curve}.

Given a completed symplectic cobordism $W=\ora{V^-V^+}$ we first
define a {\it   stable curve of height $1$}, or a {\it $1$-story stable curve }
 as a ``usual" stable curve in a sense of M. Kontsevich (see \cite{Kontsevich-CP2}),
i.e. a collection of of holomorphic curves $h_i:S_i\to
W$  from  moduli spaces
$\Mc_{g_i,r}^A(\Gamma^-_i,\Gamma^+_i)$ for various genera $g_i$
 which realize  homology classes
$A_i$, and collections of periodic orbits $\Gamma_i^\pm$, such
that certain pairs of marked points  (called special) on these
curves are required to be mapped to one point in $W$. The
stability condition means the absence of infinitesimal symmetries
of the moduli space. Let us point out, however, that in the case
when
 $W$ is a cylindrical cobordism, and in particular the almost complex structure
 $J$ is translationally invariant, we would need to consider along
 with the above moduli space its quotient under the $\R$-action.
 The stability for this new moduli space still means an absence of
 infinitesimal deformations, but it translates into an additional restriction
 on holomorphic curves.
  Namely, in the first case the stability
condition means that each constant curve has, after   removal
of the marked points, a negative Euler characteristic. In the second
case it additionally requires  that when {\it all} connected
components of the curve are straight cylinders $\g\times\R,
\g\in\Pc$ then  at least one of these cylinders should have a
marked point.

  One can define an arithmetic genus $g$ of the
resulting curve, the total sets $\bx^\pm$ and $\by$ (equal to the
union of sets $\bx_i^\pm$ and $\by_i$ for the individual curves of
the collection), and the total  absolute  homology class  $A\in
H_2(W )$ (see the discussion in Section \ref{sec:holomorphic}
above), realized by the union of all curves of the collection.

Moduli of stable curves of height $1$, denoted by
${}_1\Mc_{g,r}^A(\Gamma^-,\Gamma^+)$, form  a part of the
compactification of the moduli space
$\Mc_{g,r}^A(\Gamma^-,\Gamma^+)$.
  However, unlike the case of  closed symplectic manifolds, the
  stable curves of height $1$ are not   sufficient  to describe  the
 compactification of the moduli space $\Mc_{g,r}^A(\Gamma^-,\Gamma^+)$.

A finite sequence  $(W_1,\dots, W_k)$ of symplectic manifolds with
cylindrical ends is called a {\it chain} if  the positive end of
$W_i$ matches with the negative end of $W_{i+1}$, $i=1,\dots,k-1$.
This means that all data, assigned to an end, i.e. a contact
form, marking of periodic orbits, and an almost complex structure,
are the same for the matching ends.

Let us first suppose that
 none of the cobordisms  which form a chain $(W_1,\dots, W_k)$
is cylindrical.    Then
a {\it stable curve of height
$k$}, or a $k$-story stable curve  in the
chain $(W_1,\dots, W_k)$ is a $k$-tuple
$f=(f_1,\dots,f_k)$, where $f_{
i } \in
{}_1\wt{\Mc}_{g_i,r_i}^{A_i}(\Gamma_{i}^-,\Gamma_{i}^+;W_i,J_i)$,
 such that the  boundary data  of the curve $f_i$ at the positive end
match  the boundary data of $f_{i+1}$ on the negative one.
 One
also needs to impose the following additional equivalence relation regarding
the asymptotic markers on multiple orbits.
                    Suppose that $\g $ is
   a $k$-multiple periodic orbit, so that the holomorphic curve $f_i$ is
   asymptotically cylindrical over $\g$ at the positive end at a puncture $x^+ $, and
   $f_{i+1}$ is asymptotically cylindrical over $\g$ at  the negative end  at a puncture $x^-$.
    There are $k$ possible  positions $\mu^+_1,\dots,\mu^+_k$
    and $\mu^-_1,\dots,\mu^-_k$
  of   asymptotic  markers at  each of the punctures $x^\pm$.
   We assume here that the markers are numbered cyclically
  with respect to the orientation defined by the Reeb vector field
  at each of the punctures, and that the markers $\mu^+_1$ and
  $\mu^-_1$ are chosen for the curves $f_i$ and
  $f_{i+1}$.         Then we identify
  $f=\{\dots,f_i,     f_{i+1},\dots\}$
   with $(k-1)$ other stable curves of height $k$
  obtained by simultaneous cyclic  shift of
   the asymptotic markers at the punctures $x^+$ and $x^-$.

   The curves $f_i$, which form  a $k$-story stable curve $f=(f_1,\dots,f_k)$ are called
{\it floors},  or {\it levels} of $f$.
\medskip

   If some of the cobordisms which form the chain   $(W_1,\dots, W_k)$, say $W_{i_1},\dots,W_{i_l}$,
are cylindrical  then we will assume that the corresponding floors
of a $k$-story curve in $W=(W_1,\dots,W_k)$ are defined {\it only
up to translation.} In other words, if $W_i$ is cylindrical for some $i=1,\dots, k$
(i.e. $W_i=V_i\times\R$ and $J_i$ is translationally invariant) then
$f_i$ should be viewed as an element of
${}_1\wt{\Mc}_{g_i,r_i}^{A_i}(\Gamma_{i}^-,\Gamma_{i}^+;W_i,J_i)/\R$,
rather than ${}_1\wt{\Mc}_{g_i,r_i}^{A_i}(\Gamma_{i}^-,\Gamma_{i}^+;W_i,J_i)$.
 It will be convenient for us, however, to introduce the following convention.
 When speaking about stable holomorphic  curves  in chains which contain
 cylindrical cobordisms, we will treat the corresponding floors as curves representing
 their equivalence classes from
  ${}_1\wt{\Mc}_{g_i,r_i}^{A_i}(\Gamma_{i}^-,\Gamma_{i}^+;W_i,J_i)/\R$. Any statement about such
  curves should be understood in the sense, that {\it there exist} representatives for
  which the statement is true.

\bigskip
Let us define now the meaning of convergence of a
sequence of holomorphic curves to a stable curve of height $l$.
For $l=1$ this is   Gromov's standard  definition (see \cite{Gromov-holomorphic}). Namely, with
each stable curve
 $h=\{S_i,h_i\}\in {}_1\Mc_{g,r}^A(\Gamma^-,\Gamma^+)$   of height $1$ we associate a nodal
   surface $\widehat S$
obtained by identifying special pairs of marked points on $\coprod S_i$.  The maps $h_i$ fits together
to a continuous map $\widehat S\to W$ for which we will keep the notation
$h$.
 Let us  consider
also a smooth surface $S$ obtained by smoothing the nodes
of $\widehat S$. There exist a partitioning of  $S$
 by circles
 into  open parts   diffeomorphic to  surfaces $S_i$    with removed special points, and a map
  $g:S\to\widehat S$  which is a diffeomorphism from the complement $\tilde S$
  of the dividing circles
  in $S$ to the complement of  the double points
    in $\hat  S$, and  which collapses the partitioning circles to
  double points.
  A sequence of holomorphic  $\varphi_l:(S,j_l)\to (W,J)$ is said to converge
   to a stable curve
 $h$
    if
      the sequence $\varphi_l|_{\tilde S}$ converges
      to $h_i\circ g|_{\tilde S}$, and
      $j_l$ converges to $g^*(j)$  uniformly on compact sets,
      where $j$ is the conformal structure on the stable curve.
      Of course, we also require  convergence of marked points and asymptotic markers.
    A sequence of stable curves
    $h^j=\{S^j_i,h^j_i\}_{i=1,\dots,k}\in {}_1\Mc_{g,r}^A(\Gamma^-,\Gamma^+), j=1,\dots,$  is said to converge
    to a stable curve $h$, if $h$ can be presented as a collection
    of stable curves $h_i$, $i=1,\dots, k$, such that  $h^j_i$ converges to $  h_i$ in the above
    sense for each $i=1,\dots,k$.
    \medskip

 The convergence of
  a sequence of smooth curves to a stable curve of height  $l>1$  is understood
  in a similar sense.
  Let us  assume here     for simplicity that $l=2$ and that the   floors
    $f_{1} :     S_1\to W_1$ and $f_2 :
  S_2\to W_2$  of  a stable curve $f$ in a chain $(W_1,W_2)$
  are smooth, i.e.
  have no  special marked points.
  As in the height $1$ case let us consider
      \begin{description}
      \item{-} the smooth surface $S$
      partitioned  according to the combinatorics of our stable curve
      by circles into two open  (possibly disconnected) parts $U_1$ and $U_2$  diffeomorphic to
      the punctured surfaces $S_1$
       and $S_2$,
       \item{-} the surface
        $\hat S$ with double points obtained by collapsing
        these circles to points, and
        \item{-} the projection $g:S\to\hat S$.
      \end{description}
       Let $(W,J)=(W_1,J_1)\circledcirc (W_2,J_2)$
   be the composition of (completed) directed symplectic    cobordisms $W_1$ and $W_2$
   with compatible almost complex structures $J_1$ and $J_2$. This means that
   \begin{description}
   \item{-}    there exists a contact    hypersurface $V\subset W$ which splits
      $W$ into two cobordisms $W^0_1$ and $W^0_2$;
      \item{-} $W_1=W^0_1\cup V\times[0,\infty),\;W_2=V\times(-\infty,0]\cup W^0_2$;
      \item{-} $J|_{W^0_j}=J_1|_{W^0_j},\;\; j=1,2$;
      \item{-} $J_1$ and $J_2$ are translationally invariant at the ends
      $V\times[0,\infty)$ and $(-\infty,0]\times V$.
      \end{description}

      We denote  by
            $W^k $, $k=1,\dots,$ the quotient space of the disjoint union
            $$W^0_1\coprod V\times[-k,k]\coprod W^0_2$$ obtained by identifying $V=\partial W^0_1$
            with  and $V\times( -k)$ and $V=\partial W^0_2$
            with  and $V\times k$,  and extend the   almost
      complex  structures $J_1|_{W^0_1}$ and     $J_2|_{W^0_2}$
      to the unique almost complex structure  $J^k$ on $W^k$ which is translationally
       invariant on $V\times[-k,k]$. We also consider $W^k_1$ obtained
       by gluing $W^0_1$ and $V\times[0,k]$ along $V=\partial
       S^0_1=V\times 0$,
       and $W^k_2$ glued in a similar way from $V\times[-k,0]$ and $W^0_2$.
        We have $W_j=\bigcup\limits_{k=0}^\infty
       W^k_j$, $j=1,2$. On the other hand, $W^k_1$ and $W^k_2$
        can be viewed as submanifolds
       of $W^k$.

      \begin{definition}\label{def:convergence}
      {\rm
                Suppose that we are  given a sequence $j^k$
                 of conformal structures on the surface $S$
       and a sequence  of  $1$-story $(j^k,J^k)$-holomorphic curves
       $f^k:S\to W^k$.
       We say that this sequence converges to a stable curve $f=(f_1,
       f_2)$
       of height   $2$
       in $(W_1,W_2)$ if     there     exist two sequences of domains
       $U_1^1\subset\dots \subset U_1^i\subset \dots\subset U_1$
       and $U_2^1\subset\dots \subset U_2^i\subset \dots\subset U_2$,
        such that
        \begin{description}
        \item{} $\bigcup\limits_{k=1}^\infty U^k_i=U_i,\; i=1,2$ ;
        \item{} $f^k(U^k_i)\subset W^k_i$ for $i=1,2,\;\;k=1,\dots$;
        \item{}  for $i=1,2$  the  holomorphic curves
        $f^k|_{U^k_i}$ converges to $f_i\circ g: U_i\to W_i$,
                and the conformal structures   $j^k|_{U^k_1}$ converge to $g^*{j_i} $
                when $k\to\infty$
         uniformly on
        compact sets.
As in the case of stable curves of level $1$
we also require
convergence of marked points and asymptotic markers.
                \end{description} }
  \end{definition}
  Let us emphasize that  when some of the cobordisms  are cylindrical
  then according to  the convention which we introduced above
    one is allowed to compose  the corresponding curves  with translations to satisfy
  the above definition.

 \bigskip
 Notice that if the  cobordism $W_2$ is cylindrical, i.e.
 $W_2=V\times\R$ and $J_2$ is translationally invariant,
 then $W_1\circledcirc W_2$ can be identified with $W_1$, and thus one can talk about
 convergence of a sequence of curves
  $f^k\in {}_1{\Mc}_{g,r}^{A}(\Gamma^-,\Gamma^+;W_,J_1)$
  (where the almost complex structure
  $J_1$ is fixed!)
 to   a $2$-story curve $(f_1,f_2)$, where
 $f_1 \in{}_1\wt{\Mc}_{g_1,r_1}^{A_1}(\Gamma^-,\Gamma;W_1,J_1)$,
  $f_2\in{}_1\wt{\Mc}_{g_2,r_2}^{A_2}(\Gamma,\wt\Gamma^+;V\times\R,J_2)/\R$,
  $g=g_1+g_2,\,r=r_1+r_2, A=A_1+A_2$, and $J_2$ is translationally invariant.
 It is important to stress the point that the curve   $f_2$ is defined only up to translation.

\begin{figure}
\centerline{\psfig{figure=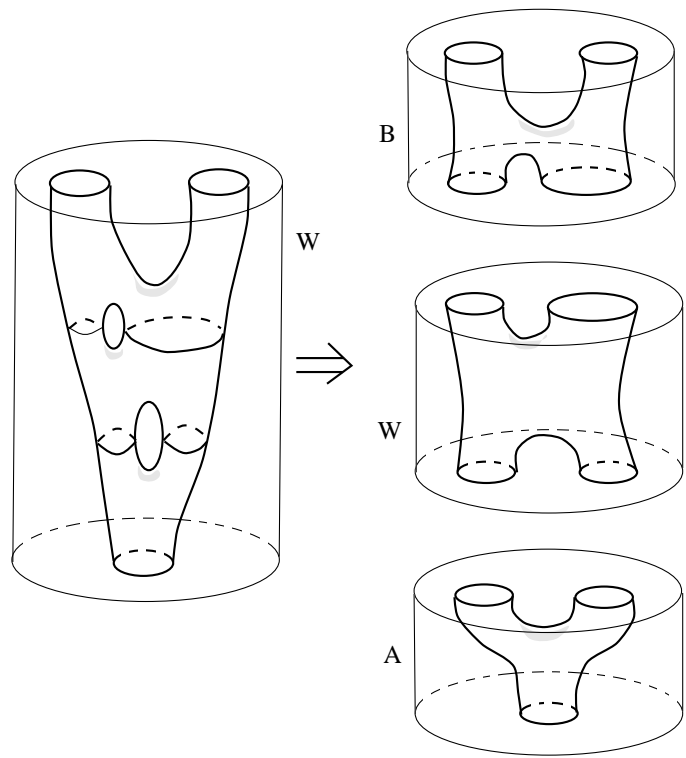,height=100mm}}
\caption{\small A possible splitting of a sequence of holomorphic curves in
 a completed symplectic cobordism}
\label{fig:split1}
\end{figure}

       \begin{theorem}\label{thm:comp1}
       Let $f_k\in {}_1\Mc_g^A(\Gamma^-,\Gamma^+)$, $k=1,\dots, $ be a
       sequence of stable holomorphic curves in
       a (complete)  directed symplectic cobordism $W$.
        Then there exists a chain of  directed symplectic
        cobordisms $$A_1,\dots, A_a,W,B_1,\dots,B_b,$$ where all
        cobordisms $A_i$ and $B_i$ are cylindrical,
        and a stable curve $f_\infty$ of height $a+b+1$ in this chain
        such that a subsequence of $ \{f_i\}$  converges to
        $f_\infty$. {\rm See Fig. \ref{fig:split1}}.
        \end{theorem}

\begin{figure}
\centerline{\psfig{figure=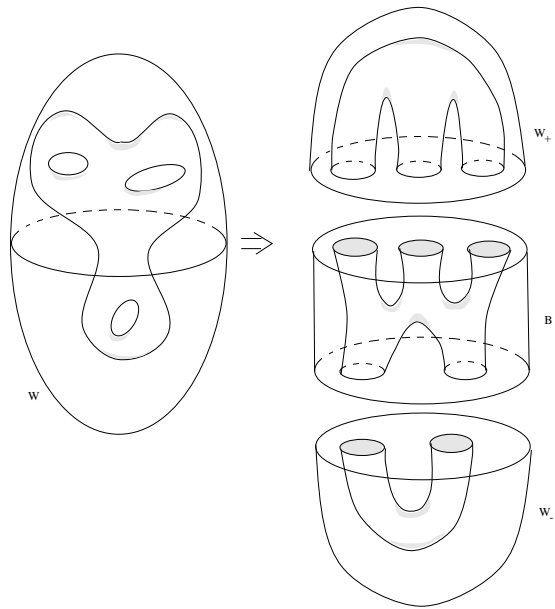,height=80mm}} \caption{\small A
possible splitting of a sequence of holomorphic curve   when
$J_k\to J_\infty$ } \label{fig:split2}
\end{figure}
        \begin{theorem}\label{thm:comp2}
        Let $W$ be a completed directed symplectic cobordism,
        $V\subset W$ a contact hypersurface, and $J_k$    a
        sequence of compatible almost complex structures on $W$
        which realizes the splitting of $W$ along $V$ into two
        directed symplectic cobordisms $W_-^\infty$ and
        $W_+^\infty$
        {\rm (see Section \ref{sec:almost} above)}.
        Let $f_k$ be a sequence of stable $J_k$-holomorphic curves
          from $ {}_1\Mc_g^A(\Gamma^-,\Gamma^+;W,J_k)$. Then there exists a
        chain of directed symplectic cobordisms
        $$A_1,\dots,A_a,W_-,B_1,\dots,B_b,W_+,C_1,\dots,C_c$$  where
        all cobordisms $A_i, B_j, C_l$ are cylindrical, such that a
        subsequence of
        $\{f_i\}$   converges to a stable curve of height $a+b+c+2$
        in
         the chain $$A_1,\dots,A_a,W_-,B_1,\dots,B_b,W_+, C_1,\dots,
         C_c.$$ {\rm See Fig. \ref{fig:split2}. The reader
         may consult \cite{HoferWysockiZehnder:99}
         for the analysis of splitting $\C P^2$
         along the boundary of a tubular
          neighborhood
         of $\C P^1\subset\C P^2$.}
        \end{theorem}

The definition of convergence can be extended in an obvious way to a sequences of stable curves
of height $l>1$. Namely, we say that a sequence of $l$-story curves
$f^k=(f^k_1,\dots,f^k_l)$, $k=1,\dots,\infty$, in a chain $(W_1,\dots,W_l)$ converge to a stable $L$-story,
$L=m_1+\dots+m_l$, curve  $f=(f_{11},\dots,,f_{1m_1},\dots,f_{l_1}\dots f_{lm_l}) $ in a chain
 $$(W_{11},\dots,W_{1m_1},\dots,W_{l_1}\dots W_{lm_l})$$
if for each $i=1,\dots,l$  the cobordism $W_i$ splits into the composition
$$W_i=W_{i1}\circledcirc\dots\circledcirc W_{im_i}$$
and the sequence $f^k_i$, $k=1,\dots,\infty$, of stable curves of height $1$ converges
to the $m_i$-story curve $f_i=(f_{i1},\dots,f_{im_i})$ in the chain $(W_{i1},\dots, W_{im_i})$ in the sense of Definition
\ref{def:convergence}.
\bigskip

It is important to combine Theorems \ref{thm:comp1} and
\ref{thm:comp2} with the following
observation which is a
corollary of  Stokes' theorem combined with Proposition \ref{prop:positivity}.

\begin{proposition}\label{prop:one-end}
 A holomorphic curve in an exact
  directed symplectic cobordism
 (and in particular in a cylindrical one)
 must have at least one positive puncture.
 \end{proposition}

 In particular,  we have
 \begin{corollary}\label{cor:pos-puncture}
Let $f^n\in {\Mc}_0(W,J)$  be a sequence of rational holomorphic curves with  one positive, and possibly
several negative punctures.
 Suppose that the sequence  converges to a stable
 curve $$F=\{g_1,\dots,g_a,f,h_1,\dots,h_b\}$$ of height $a+b+1$ in a chain
 $$A_1,\dots,A_a,W,B_1,\dots,B_b.$$ Then the  $W$-component   $f\in {\Mc}_0(W,J)$ of the stable curve
 $F$ has precisely
 one positive puncture as well.
 \end{corollary}

 \subsection{Dimension of the moduli spaces
$\Mc^A_{g,r}(\Gamma^-,\Gamma^+)$}
One has the following index formula for the corresponding
$\bar\partial$-problem which compute the dimension
of the moduli space $\Mc^A_{g,r}(\Gamma^-,\Gamma^+;W,J)$ for a {\it
generic}
 choice of $J$.
\begin{proposition}

  \begin{equation}\label{eq:dim}
  \begin{split}
 & {\mathrm {dim}}\Mc^A_{g,r}(\Gamma^-,\Gamma^+;W,J)
 =\sum\limits_1^{s^+}\CZ(\gamma^+_i)-
 \sum\limits_1^{s^-}\CZ(\gamma^-_k)\\
 &  +  (n-3)(2-2g-s^+-s^-)+2c_1(A)+2r,\\
 \end{split}
 \end{equation}
  where
$s^\pm$ are the cardinalities of the sets $\Gamma^\pm$, and $c_1\in
H^2(W)$ is the first Chern class of the almost complex manifold
$(W,J)$
\end{proposition}

  Making the moduli spaces non-singular by picking
   generic $J$ is needed for the purpose of curve counting but
            does not always work properly. It is therefore crucial
             that the moduli spaces of stable $J$-holomorphic
curves are non-singular {\em virtually}. This means that for {\em
any} $J$ the moduli spaces, being generally speaking singular, can
be equipped with some canonical additional structure that make
them {\em function in the theory the same way as if they were
orbifolds with boundary and had the dimension prescribed by the
Fredholm index}. In particular, the moduli spaces come equipped
with rational fundamental cycles relative to the boundary (called
{\em virtual fundamental cycles}) which admit pairing with
suitable de Rham cochains and allow us to use the Stokes
integration formula.

Technically the virtual smoothness is achieved by a
finite-dimensional reduction of the following picture: a singular
moduli space is the zero locus of a section defined by the
Cauchy-Riemann operator in a suitable orbi-bundle over a moduli
orbifold of stable $C^{\infty}$-maps. More general {\em virtual
transverality} properties for families of
 $J$'s also hold true
(cf. \cite{Fukaya-Ono, Fukaya-Ono2, Liu-Tian, Li-Tian, Ruan-virtual, Siebert, McDuff} at al.).
We
are reluctant to provide in this quite informal exposition
 precise
formulations   because of numerous not entirely innocent
subtleties this would entail.  Fortunately, what we intend to say in the rest of this
paper does not depend much on the details we are omitting.

\bigskip

As it was explained in Section \ref{sec:compact} above the moduli
space $\Mc^A_{g,r}(\Gamma^-,\Gamma^+;W,J)$  can be compactified by
adding strata which consist of stable holomorphic curves of
different height.  This compactification looks quite similar to
the Gromov-Kontsevich compactification  of   moduli spaces of
holomorphic curves in a closed symplectic manifold  with a
compatible almost complex structure.
There is, however, a major difference. In the case of a closed
manifold   all the strata   which one needs to add to compactify the moduli space
of smooth holomorphic curves have (modulo virtual cycle complications) codimension
$\geq 2$. On the other hand  in our case  the  codimension one strata are
present {\it generically}. Thus in this case the boundary of the compactified moduli space,
rather than the moduli space  itself, carries  the fundamental
cycle.

 In particular, this boundary
  is tiled by codimension one
 strata represented by stable curves $(f_-,f_+)$ of {\it height two}.
Each such a stratum can be described by the constraint matching
the positive ends of $f_-$ with the negative ends of $f_+$ in the
Cartesian product of the moduli spaces $\Mc_{\pm}$ corresponding to   the
curves
$f_{\pm}$  separately.
Proposition \ref{prop:boundary} below describes these top-dimensional boundary
 strata more precisely
in two important for our purposes situations.
Let us point out that  \ref{prop:boundary} literally holds only under certain
 transversality conditions.
Otherwise it should be understood only virtually.

 \begin{proposition}\label{prop:boundary}
     \begin{description}
     \item{1.} Let $(W=V\times\R,J)$ be a cylindrical cobordism.
     Then any
    top-dimensional stratum  $\cS$ on the boundary  of the compactified moduli space
     $\overline{\Mc^A_{g,r}(\Gamma^-, \Gamma^+;W,J)/\R}$ consists
     of
         stable curves  $(f_-,f_+)$   of height   two,
          $f_\pm\in\Mc_\pm/\R$, where
          $$\Mc_-=\wt\Mc^{A_-}_{g_-,r_-}(\Gamma^-,\Gamma;W,J)/\R ,
          \quad \Mc_+=\wt\Mc^{A_+}_{g_+,r_+}(\Gamma,\Gamma^+;W,J)/\R ,$$
           $g=g_-+g_+, \;r=r_-+r_+,\;A=A_-+A_+$,   $\Gamma=\{\g_1,\dots,\g_l\}\subset\Pc$.

                 All  but one
          connected components  of  each of the curves $f_-$ and $f_+$
           are trivial cylinders
          (i.e.  have the form   $\g\times\R, \g\in\Pc$) without
          marked points.
          \item{2.} Let $(W=\ora{V^-V^+},J)$  be any    cobordism,
          and $(W_\pm,J_\pm)= (V^\pm\times\R,J_\pm)$
          be the cylindrical cobordisms associated to its
          boundary.
          Then any  top-dimensional strata $\cS$ on the boundary  of the
          compactified moduli space
     $\overline{\Mc^A_{g,r}( \Gamma^-, \Gamma^+;W,J) }$ consists
     of
         stable curves  $(f_-,f_+)$   of height   two,
          $f_\pm\in\Mc_\pm$, where either
          $$\Mc_-=\wt\Mc^{A_-}_{g_-,r_-}(\Gamma^-,\Gamma;W_-,J_-)/\R,\;
 \Mc_+=\wt\Mc^{A_+}_{ g_+,r_+}(\Gamma,\Gamma^+;W,J)\;\hbox{
 and}\;
  \Gamma\subset\Pc^-,$$ or
$$\Mc_+=\wt\Mc^{A_+}_{g_-,r_-}(\Gamma ,\Gamma^+;W_+,J_+)/\R,\;
 \Mc_-=\wt\Mc^{A_-}_{g_-,r_-}(\Gamma^-,\Gamma;W,J ) \;\hbox{
 and}\;
  \Gamma\subset\Pc^+.$$
   In both cases  we have  $$g=g_-+g_+,
 \;r=r_-+r_+,\;\;\hbox{and}\;\;A=A_-+A_+.$$
  The part of the stable curve
 $(f_-,f_+)$ which is  contained in in  $W_\pm$ must have
 precisely one non-cylindrical connected component, while there
 are no restrictions on the number and the character of connected components or the
 other part of the stable curve.
 \item{~~} In both cases the stratum
         $\cS=\cS(\Gamma,g_-|g_+,r_-|r_+,A_-|A_+)$ is
         diffeomorphic to a $\kappa$-multiple cover of the product
         $\Mc_-\times\Mc_+$, where the multiplicity $\kappa$ is
         determined by the multiplicities of periodic orbits from
         $\Gamma$.
 \end{description}
              \end{proposition}

Proposition \ref{prop:boundary} is not quite sufficient for our
purposes, as we also needs to know the structure of the boundary
of   moduli spaces of $1$-parametric families of holomorphic
curves. However, we are not  formulating the corresponding
statement in this paper, because it  is   intertwined  in a much more serious way
with
the virtual cycle techniques and terminology. An algebraic description
of this boundary is given in
   Theorem
\ref{thm:SFT-chain}
below.

\bigskip

Let us consider some special cases of the formula (\ref{eq:dim}).
Suppose, for instance, that $W$ is the cotangent bundle of a manifold
$L$. Then $W$ is a symplectic manifold which has  only a positive
cylindrical end.
 If $L$ is orientable then there is a canonical way to
define Conley-Zehnder indices. Namely, one takes any
trivialization along orbits, which is  tangent to vertical
Lagrangian fibers. The resulting index is independent of a
particular  trivialization. For this trivialization, and a choice
of a contact form corresponding to  a metric on $L$ we have
 \begin{proposition}
Periodic orbits of the Reeb flow are lifts of closed geodesics,
and if $L$ is orientable their Conley-Zehnder indices   are equal
to Morse indices of the corresponding geodesics and we have
 $${\mathrm {dim}}\Mc_g^A( \Gamma^+)= \sum\limits_i\Morse(\g^+_i)+ (n-3)(2-2g-s^+) .$$
\end{proposition}\label{prop:Morse-Zehnder}
Notice that for a metric on  $L$ of  non-positive curvature we have  $s^+>1$, because in this case
there are no
contractible geodesics. Moreover, if the metric has negative curvature then all geodesics
have  Morse
indices equal to $0$.  Hence, we get
\begin{corollary}\label{cor:isolated} In  the cotangent bundle of a  negatively curved
manifold  of dimension $>2$ there  could be only isolated
holomorphic curves. If, in addition, $n\neq 3$ then these curves are spheres
with two positive punctures. Each of these curves is asymptotically
cylindrical at punctures over   lifts of the same geodesic with
opposite orientations.
\end{corollary}
Let us point out that the orientability is not required in
Corollary \ref{cor:isolated}. The corresponding result for a
non-orientable manifold  follows from \ref{prop:Morse-Zehnder}
applied to its orientable double  cover.

\subsubsection*{Absence of hyperbolic Lagrangian submanifolds  in uniruled manifolds}

As the first  application of the above compactness theorems let us
prove here the following theorem of C. Viterbo. Let us recall that
a  complex projective manifold $W$ is called uniruled, if there is a
rational holomorphic curve
  through each point of $W$. For instance,  according to
  Y. Myaoka--S. Mori
 \cite{Mori-Miayoka}  and J. Kollar  \cite{Kollar}  Fano manifolds are uniruled.

\begin{theorem} {\rm (C. Viterbo, \cite{Viterbo})}
Let $W$ be a uniruled manifold of complex dimension $>2$,  $\omega$ its K\"ahler
sympletic form, and $L\subset W$ an embedded Lagrangian
submanifold. Then $L$ does not admit a Riemannian metric of
negative sectional  curvature.
\end{theorem}
\proof J. Kollar  \cite{Kollar}  and in a more general case Y.
Ruan   \cite{Ruan-virtual} proved  that there exists a homology
class $A\in H_2(W)$, such that for any almost complex structure
compatible with $\omega$ and any point $z\in W$ there exists
$f\in\Mc_{0,1}^A(W,J)$ with $f(y)=z$, where $y$ is the marked
point. Let us identify a neighborhood $U$ of $L$ in $W$ with  a
neighborhood of the zero-section in $ T^*(L)$. Suppose $L$ admits
a Riemannian metric of negative
  curvature. We can assume that  $U$ is the round neighborhood of radius $1$ in
  $T^*(L)$. Let us consider a sequence $J^m$ of almost complex structures
  on $W$, which realizes the splitting along the contact  type hypersurface
  $(V=\partial U,\alpha=pdq|_V)$.
  (see Section \ref{sec:almost}). Then according to Example \ref{ex:cotangent}
  $W$ splits into $W_-=T^*(L)$ and $W_+=W\setminus L$. The almost
  complex structure $J_-$ on $T^*(L)$ is compatible at infinity  with the
  contact $1$-form   $\alpha=pdq|_{V}$. According to Corollary
  \ref{cor:isolated} for any choice $\Gamma=\{\g_1,\dots,\g_k\}$ and
  any $g\geq 0$ the
  moduli spaces $\Mc_g(\Gamma;W_-,J_-)$ are empty, or
  $0$-dimensional. One the other hand, Theorem \ref{thm:comp2}
  together with Ruan's theorem guarantee the existence of a
  rational holomorphic curve  with punctures through every point
  of $L$. This contradiction proves that
  $L$ cannot admit a metric
  of negative curvature.  \qed

\subsection{ Coherent orientation of the moduli spaces
  of holomorphic curves}\label{sec:orientation}

    To get started with the algebraic formalism,
    one first needs to orient moduli spaces
    $\Mc(\Gamma^-,\Gamma^+)$ of holomorphic curves with punctures.
    This problem is much easier in the case
    of moduli spaces of closed holomorphic curves, because in that
    case moduli spaces are even-dimensional and carry a canonical almost
     complex structure (see Section \ref{sec:coherent-closed} below).
    In our case we have to adapt the philosophy of {\it coherent orientations}
     of the moduli spaces
    borrowed from Floer homology theory (see \cite{Floer-Hofer}).
    We sketch this approach in this section.

\subsubsection{Determinants}

In order to separate the problems of orientation and
transversality we are going to orient the determinant line bundles
of the linearized $\overline\partial$-operators, rather than  the moduli spaces themselves.

 For a linear Fredholm operator $F:A\rightarrow
B$ between Banach spaces we can define its determinant line
$\det(F)$ by $$ \det(F)= (\Lambda^{max} \Ker (F))\otimes
(\Lambda^{max} \Coker (F))^{\ast}. $$ We note that for the trivial
vector space $\{0\}$ we have $\Lambda^{max}\{0\}= {\R}$.  An
orientation for $F$ is by definition an orientation for the line
$\det(F)$. In particular, given an isomorphism $F$ we can define a canonical
orientation given through the vector $1\otimes 1^{\ast}\in {\R
}\otimes {\R}^{\ast}$.

Given a continuous family $F=\{F_y\}_{y\in Y}$ of Fredholm
operators, parameterized by a topological space $Y$,  the
determinants of operators $F_y$ form   a line bundle
$\det(F)\rightarrow Y$. The fact that this is a line bundle in a
natural way might be surprising since the dimensions of kernel and
cokernel vary in general. This is however a standard fact, see for
example \cite{Floer-Hofer}.

\subsubsection{Cauchy-Riemann Type Operators on Closed Riemann Surfaces}
\label{sec:coherent-closed} Let $(S,j)$ be a closed, not
necessarily connected Riemann surface and $E\rightarrow S$ a
complex vector bundle. Denote by $X_E\rightarrow S$ the complex
$n$-dimensional vector bundle whose fiber over $z\in S$ consists
of all complex ant-linear maps $$ \phi:T_zS\rightarrow E_z,\ z\in
S,\ \hbox{i. e.}\ J\circ\phi+\phi\circ j=0, $$ where $J$ is the
complex structure on $E$. Fixing a connection $\nabla$ and a
smooth $a\in\hbox{Hom}_{\mathbb R}(E,X_E)$ we can define a
Cauchy-Riemann type operator $$L:C^{\infty}(E)\rightarrow
C^{\infty}(X_E)$$ by the formula $$
(Lh)(X)=\nabla_Xh+J\nabla_{jX}h +(ah)(X), $$ where $X$ is an
arbitrary vector field on $S$. Since the space of connections is
an affine space we  immediately see that the set ${\mathcal O}_E$
of all Cauchy-Riemann type operators on $E$ is convex. For a
proper functional analytic set-up, where we may chose H\"older or
Sobolev spaces, the operator $L$ is Fredholm. By elliptic
regularity theory the kernel and cokernel would be spanned always
by the same smooth functions, regardless which choice we have
made. The index of $L$ is given by the Riemann-Roch formula $$
\hbox{ind}(L)= (1-g)  \hbox{dim}_{\mathbb R}(E) + 2 c(E), $$ where
$c(E)$ the first Chern number $c_1(E)(S)$ of $E$. Here we assume
$S$ to be a connected closed surface of genus $g=g(S)$.

 Let   $\phi:(S,j)\rightarrow (T,i)$
 be a biholomorphic map and $\Phi:E\rightarrow F$ a ${\mathbb
C}$-vector bundle isomorphism covering $\phi$. Then $\Phi$ induces
an isomorphism $$ { \Phi}_{\ast}:{\mathcal O}_E\rightarrow
{\mathcal O}_F $$ in the obvious way. The operators   $(E,L)$ and
$(F,K)$ are called isomorphic  if there exists ${
\Phi}:E\rightarrow F$, so that ${\bf \Phi}_{\ast}(L)=K$. We will
denote by $[E,L]$ the equivalence class  of an operator  $(E,L)$
which consists of operators $(F,K)$, equivalent to $(E,L)$ under
isomorphisms, {\it isotopic to the identity}, and by $[[E,L]]$ the
equivalence class under the action of the full group of
isomorphisms. The moduli space of equivalence classes $[[E,L]]$
will be denoted by $\CR$, and the ``Teichmuller space" which
consists of classes $[E,L]$ will be denoted by $\wt{\CR}$. An
isomorphism $\Phi$ induces an isomorphism between the kernel
(cokernel) of $L$ and ${ \Phi}_{\ast}L$ for every $L\in {\mathcal
O}_E$, and hence one can canonically associate the determinant
line to an isomorphism class, and thus define the {\it determinant
line bundle} $\cV$ over the moduli space $\CR$.
 Given an orientation $o$ for $L$ we obtain an induced
orientation ${ \Phi}_{\ast}(o)$.  Let us note the following
\begin{lemma}
 The bundle $\cV$ is orientable.
\end{lemma}
\proof The lift $\wt\cV$ of the bundle $\cV$ to the Teichmuller
space $\wt\CR$ is obviously orientable, because each connected
component of the space $\wt\CR$ is contractible. However, one
should check that an arbitrary isomorpism $\Phi:(E,L)\to (F,K)$
preserves the orientation. This follows from the following
observation. Any connected component of $\wt\CR$ contains an
isomorphism class of a complex linear operator $(E,L_0)$,  and any
two complex linear operators representing points in a given
component of $\wt\CR$ are homotopic in the class of complex linear
operators. The determinant of $(E,L_0)$ can be oriented
canonically by observing that its kernel and cokernel are complex
spaces.  Any isomorphism  maps a   complex linear operator to a
complex linear operator and preserves its complex orientation.
Hence, it preserves an orientation of the determinant line  of any
operator $(E,L)$.
\qed
\medskip

We will call an orientation of $\cV$ {\it complex} if it coincides
with    the complex orientation of   determinants   of complex
linear operators.

\bigskip

The  components of the space $\CR$ are
parameterized by the topological type of the underlying surface
$S$ and the isomorphism class of the bundle $E$.
It turns out that the complex orientation of $\cV$ satisfies three
{\it coherency} Axioms A1--A3 which we formulate below.
They relate  orientations of $\cV$ over different components of
$\CR$.
Conversely, we will see that these axioms determine the
orientation uniquely up to a certain normalization.

Given $(E,L)$ and $(F,K)$ over surfaces $\Sigma_0$ and $\Sigma_1$
 we define a {\it disjoint union}
$$ [E,L]\dot{\cup}[F,K] := [G,M]  $$ of $(E,L)$ and $(F,K)$ as a pair
$(G,M)$, where $G$ is a bundle over     the disjoint union
$\Sigma=\Sigma_0\coprod\Sigma_1$, so that $(G,M)|_{\Sigma_0}$ is
isomorphic to $(E,L)$ and $(G ,M)_{\Sigma_1}$  is isomorphic to
$(F,K)$. Clearly, the
isomorphism class of a disjoint union is uniquely determined by the classes
of $(E,L)$ and $(F,K)$.
Thus, we have a well-defined construction called {\it disjoint
union}:
 The determinant $\det\Sigma$ is canonically isomorphic to $\det
 L\otimes\det K$, and hence the orientations    $o_K$ and $o_L$
 define
 an orientation $o_K\otimes o_L$ of    $ \Sigma $.
 Our first axiom  reads
 \medskip

\noindent{\textsf{Axiom C1}}. For any disjoint union $
[G,M]=[E,L]\dot{\cup}[F,K]$ the orientation $o_M$ equals
 $o_K\otimes o_L$.
 \bigskip

Given $(E,L)$ and $(F,K)$, where $E$ and $F$ are  bundles over $S$
of possibly different rank, we can define an operator $(E\oplus
F,L\oplus K)$. There is a canonical map $$
\hbox{det}(L)\otimes\hbox{det}(K)\rightarrow \hbox{det}(L\oplus
K),
$$ and thus  given orientations $o_L$ and $o_K$ we obtain $o_L\oplus o_K$.

\bigskip
\noindent{\textsf{Axiom C2}}.$$o_{L\oplus K}=o_L\oplus o_K.$$
 \bigskip

To formulate the third axiom, we need a construction, called {\it
cutting and pasting}.

Let $(E,L)$ be given and assume that
$\gamma_1,\gamma_2:S^1\rightarrow S$ be real analytic embeddings
with mutually disjoint images. Assume that $\Phi:E|_{\gamma_1}
\rightarrow E|_{\gamma_2}$ is a complex vector bundle isomorphism
covering $\sigma=\gamma_2\circ\gamma_{1}^{-1}$. The maps
$\gamma_1$ and $\gamma_2$   extends as holomorphic
embeddings $\bar{\gamma}_j:[-\varepsilon,\varepsilon]\times S^1
\rightarrow S$ for a suitable small $\varepsilon>0$, so that the
images are still disjoint. Locally, near $\gamma_j$ we can
distinguish the left and the right side of $\gamma_j$. These sides
correspond to the left or the right part of the annulus
$[-\varepsilon,\varepsilon]\times S^1$. Cutting $S$ along the
curves $\gamma_j$ we obtain a compact Riemann surface $\bar{S}$
with boundary. Its boundary components are $\gamma_{j}^{\pm}$,
$j=1,2$, where $\gamma^{\pm}_{j}$ is canonically isomorphic to
$\gamma_j$. The vector bundle $E$ induces a vector bundle
$\bar{E}\rightarrow \bar{S}$. We define a space of smooth sections
$\Gamma_{\Delta}(\bar{E})$ as follows. It consists of all
smooth sections $\bar{h}$ with the property that $$
\bar{h}|_{\gamma_{j}^{-} }= \bar{h}|_{\gamma_{j}^{+}}\ \hbox{for}\
j=1,2. $$ Then $L$ induces an operator
$\bar{L}:\Gamma_{\Delta}(\bar{E}) \rightarrow
\Gamma(X_{\bar{E}})$. The operators  $L$ and $\bar{L}$
have naturally isomorphic kernel and cokernel. So an orientation $o$ of $\hbox{det}(L)$
induces one of $\hbox{det}(\bar{L})$. The boundary condition
$\Delta$ can be written in the form
\begin{eqnarray*}
\left[\begin{array}{cc} 1&0\\ 0&1
\end{array}\right]
\cdot\left[\begin{array}{c}
\Phi(\gamma_{1}(t)\bar{h}\circ\gamma_{1}^{-}(t)\\
\bar{h}\circ\gamma_{2}^{-}(t)
\end{array}\right]
=\left[\begin{array}{c}
\Phi(\gamma_1(t))\bar{h}\circ\gamma^{+}_{1}(t)\\
\bar{h}\circ\gamma^{+}_{2}(t)
\end{array}\right]
\end{eqnarray*}
We introduce a parameter depending boundary condition by
\begin{eqnarray*}
\left[\begin{array}{cc} cos(\tau)&sin(\tau)\\ -sin(\tau)&cos(\tau)
\end{array}\right]
\cdot\left[\begin{array}{c}
\Phi(\gamma_{1}(t)\bar{h}\circ\gamma_{1}^{-}(t)\\
\bar{h}\circ\gamma_{2}^{-}(t)
\end{array}\right]
=\left[\begin{array}{c}
\Phi(\gamma_1(t))\bar{h}\circ\gamma^{+}_{1}(t)\\
\bar{h}\circ\gamma^{+}_{2}(t)
\end{array}\right]
\end{eqnarray*}
for $\tau\in [0,\frac{\pi}{2}]$. For all these boundary conditions
$L$ induces an operator, which is again Fredholm of the same
index. For every $\tau$ we obtain a Cauchy-Riemann type operator
from $\Gamma_{\Delta_{\tau}}(\bar{E}) $ to $\Gamma(X_{\bar{E}})$.
Note that for a section $h$ satisfying the boundary condition
$\Delta_{\tau}$ the section $ih$ satisfies the same boundary condition.
On the other hand for
$\tau=\frac{\pi}{2}$ we obtain a Fredholm operator whose
kernel and cokernel naturally isomorphic to the kernel on cokernel of
 a   Fredholm operator
on a new closed surface. Namely identify $\gamma^{+}_{1}$ with
$\gamma^{-}_{2}$ and $\gamma^{+}_{2}$ with $\gamma^{-}_{1}$. For
the bundle $\bar{E}$ we identify the part above $\gamma^{+}_{1}$
via $\Phi$ with the part over $\gamma_{2}^{-}$ and we identify the
part above $\gamma^{-}_{1}$ with $-\Phi$ to the part above
$\gamma^{+}_{2}$. The latter surface and bundle we denote by
$E_{\Phi}\rightarrow S_{\Phi}$ and the corresponding operator by
$L_{\Phi}$. Letting the parameter run we obtain starting with an
orientation $o$ for $L$ an orientation $o_{\Phi}$ for $L_{\Phi}$.
If $o$ is the complex orientation it is easily verfied that
$o_{\Phi}$ is the complex orientation as well. We  say that the operator
$L_{\Phi}$  is an operator obtained
from $L$ by cutting and pasting.  This operator $L_\Phi$ has the
same index as $L$,  and the component of $[L,\Phi]$ in $\wt\CR$
depends only the isotopy classes of the embeddings $\gamma_1$ and
$\gamma_2$.

\bigskip

\noindent{\textsf{Axiom C3}}. $$o_{L_\Phi}=o_\Phi\, .$$
 \bigskip

 Note that we have to require here that the parts of $L$ over the curves
  $\gamma_1$ and $\gamma_2$ are isomorphic via the gluing data.
It is straightforward to check that
\begin{theorem}\label{thm:complex-coherent}
The complex orientation of $\cV$ is coherent, i.e. it satisfies
Axioms C1--C3.
\end{theorem}

Let us point out a simple
 \begin{lemma}
  Let $(E,L)$ be  an isomorphism then the orientation by $1\otimes 1^{\ast}$
  of the $\det
L={\mathbb R}\otimes {\mathbb R}^{\ast}$ defines the complex
orientation of $\cV$ over the component of $[E,L]$.
\end{lemma}
 The following
theorem gives the converse of Theorem \ref{thm:complex-coherent}.

\begin{theorem}\label{thm:coherent-complex}
Suppose that a coherent orientation of $\cV$  coincides with
the complex orientation for  the trivial line bundle over $S^2$
and for
the line bundle over $S^2$ with Chern number $1$. Then the
orientation is complex.
\end{theorem}

\proof
Let us first observe that according to Theorem
\ref{thm:complex-coherent} the disjoint union, direct sum and
cutting and pasting procedures preserve the class of complex
orientations.  Consider the pair $(E_0,L_0)$, where  $E_0$
is the trivial bundle $ S^2\times {\mathbb C}\to S^2$ and $L_0$ is the
standard Cauchy-Riemann operator. Then the $\ind L_0=2$. Take
small loops around north pole and south pole on $S^2$ and
 identify the trivial bundles over these loops. Now
apply the cutting and pasting procedure and Axiom C3 to obtain the
disjoint union of the trivial bundle over the torus and the
trivial bundle over $S^2$. Hence we can use Axiom C1 to obtain an
induced orientation for the Cauchy-Riemann operator on the trivial
bundle over $T^2$. Taking appropriate loops we obtain orientations
for all trivial line bundles over Riemann surfaces of arbitrary
genus. Using direct sums
 and disjoint unions constructions, and applying Axioms C1 and C2  we see that  the orientation of
  all trivial bundles
of arbitrary dimensions over Riemann surfaces of arbitrary genus
are complex. Let $E_1$ be the bundle over $S^2$ with Chern number
$1$.  Then we can use C3 to glue  two copies of $(E_1,L_1)$ to
obtain the complex orientation of the disjoint union of a Cauchy
Riemann operator on the trivial bundle and one on the bundle with
Chern number $2$. Now it is clear that  the given coherent
orientation has to  be complex over all components of the moduli space
$\CR$. \qed

In the next section we extend the coherent orientation from
Cauchy-Riemann type operators over closed surfaces to
  a  special class of Cauchy-Riemann type
operators on Riemann surfaces with punctures.

\subsubsection{A special class of Cauchy-Riemann type
operators on punctured Riemann surfaces}\label{sec:coherent-punct}

Let us view ${\mathbb C}^n$ as a real vector space equipped with
the Euclidean  inner product which is the real part of the
standard Hermitian inner product. We define a class of
self-adjoint operators as follows. Their domain in
$L^2(S^1,{\mathbb C}^n)$
 is
$H^{1,2}(S^1,{\mathbb C}^n)$ of  Sobolev maps $h:S^1\rightarrow
{\mathbb C}^n$.  The operators have the form
\begin{equation}\label{eq:asymptotic}
(Ax)(t)=-i\frac{dx}{dt}-a(t)x,
\end{equation}
 where $a(t)$ is a smooth loop of
real linear self-adjoint maps. We assume that $A$ is
non-degenerate in the sense that $Ah=0$ only has the trivial
solution, which just means that the time-one map $\psi(1)$ of the
Hamiltonian flow
\begin{equation}\label{eq:asympt-flow}
\begin{split}
 \dot{\psi}(t)=&i a(t)\psi(t),\\
 \psi(0)=&\Id\\
 \end{split}
 \end{equation}
 has no eigenvalues equal  to $1$.
  In particular, $A:H^{1,2}\rightarrow L^2$ is an
isomorphism.

Given a smooth vector bundle $E\rightarrow S^1$ we can define
$H^{1,2}(E)$ and $L^2(E)$ and a class of operators $B$ by
requiring that $A=\Phi B\Phi^{-1}$ for  an Hermitian
trivialization $\Phi$ of the bundle $E$.  We shall call such
operators {\it asymptotic}, for reasons which will become clear
later.

 As it was defined in Section \ref{sec:holomorphic} above,
  an asymptotically marked
punctured Riemann surface is a triplet $(S, \bx , \mu )$, where
$S=(S,j)$ is a closed Riemann surface,
 $\bx =\{x _1,\dots,x _{s }\}$
is the set  of  punctures, some of them called positive, some
negative, and $\bmu=\{\mu_1,\dots,\mu_s\}$   is the set of
asymptotic markers, i.e. tangent rays,
 or equivalently oriented
tangent lines at the punctures.

One can introduce near each puncture $x_k\in\bx$ a holomorphic
parameterization, i.e. a holomorphic map $h_k:D\to S$ of the unit
disk $D$ such that $h_k(0)=x_k$ and the asymptotic marker $\mu_k$
is tangent to the ray $h_k(r),\;r\geq 0$. We assume that the
coordinate neighborhoods $\cD_k=h_k(D)$ of all the punctures are
disjoint. Then we define $\sigma_k:{\mathbb R}^+\times
S^1\rightarrow {\mathcal D}\setminus\{0\}$ by $$
\sigma_k(s,t)=h_k(e^{\pm{2\pi(s+it)}}), $$ where the sign $-$ is
chosen if the  puncture $x_k$ is positive, and the sign $+$ for
the negative puncture. We will refer to $\sigma_k$ as holomorphic
polar coordinates adapted to $(x_k,\mu_k)$. Given two adapted
polar coordinate systems $\sigma $ and $\sigma'$ near the same
puncture $x\in\bx$
  we observe that the transition map
(defined for $R$ large enough) $$ \sigma
^{-1}\circ\sigma':[R,\infty)\times S^1 \rightarrow
[0,\infty)\times S^1 $$ satisfies for every multi-index $\alpha$
$$ D^{\alpha}[\sigma^{-1}\circ\sigma'(s,t)-(c+s,t)] \rightarrow 0
$$ uniformly for $s\rightarrow\infty$, where $c$ is a suitable
constant. The main point is the fact that there is no phase shift
in the $t$-coordinate.

Given $(S, \bx , \mu )$ we associate to it a smooth surface
$\bar{S}$ with boundary compactifying the punctured Riemann
surface $S\setminus\bx$ by adjoining a circle for every puncture.
Each circle has a distinguished point $0\in S^1={\mathbb
R}/{\mathbb Z}$. Namely for each   positive puncture we compactify
${\mathbb R}^+\times S^1$ to $[0,\infty]\times S^1$, where
$[0,\infty]$ has the smooth structure making the map $$
[0,\infty]\rightarrow [0,1]:s\rightarrow s(1+s^2)^{-\frac{1}{2}},\
\infty\rightarrow 1 $$ a diffeomorphism. We call
$S^{+}_{k}=\{\infty\}\times S^1$ the circle at infinity associated
to $(x_k,\mu_k)$. For negative punctures we compactify at
$-\infty$ in a similar way.

\begin{definition}
A smooth complex vector bundle $E\rightarrow  (S, \bx , \mu )$ is
a smooth vector bundle over $\hat{S}$ together with Hermitian
trivializations $$ \Phi_k:E|_{S_{k}}\rightarrow S^1\times{\mathbb
C}^n. $$ An isomorphism between two bundles $E$ and $F$ over
surfaces $S$ and $T$ is a a complex vector bundle isomorphism
$\Psi:E\to F$ which covers   a biholomorphic map
$\phi:(S,j)\rightarrow (T,i)$, preserves punctures and the
asymptotic markers  (their numbering and signs)  and respects the
asymptotic trivializations.
\end{definition}

Define as in Section \ref{sec:coherent-closed} above   the bundle
$X_E\rightarrow \bar{S}$. Set $\dot{S}=S\setminus\Gamma$. We
introduce the Sobolev space $H^1(E)$ which  consists of all
sections $h$ of $E\rightarrow \dot{S}$ of class $H^{1,2}_{loc}$
with the following behavior  near   punctures.
  Suppose, that $x$ is a positive puncture and   $\sigma$
is an adapted system of  holomorphic polar coordinates. Pick a
smooth   trivialization $\psi$ of $E\rightarrow \bar{S}$ over
$[0,\infty]\times S^1$ (in local coordinates) compatible with the
given asymptotic trivialization. Then the map $(s,t)\rightarrow
\psi(s,t)h\circ\sigma(s,t)$ is assumed to belong to
$H^{1,2}({\mathbb R}^+\times S^1,{\mathbb C}^n)$. A similar
condition is required for   negative punctures. In a similar
way we define  the space $L^2(X_E)$. Observe that defacto we use
measures which are infinite on $\dot{S}$ and that the
neighborhoods of punctures look  like half-cylinders.

A Cauchy-Riemann type operator $L$ on $E$ has the form $$ (Lh)X =
\nabla_Xh+J\nabla_{jX}h +(ah)X, $$ where $X$ is a vector field on
$S$.   We require, however a particular behaviour of $L$  near the
punctures. Namely, regarding   $E$ as a trivial bundle
$[0,\infty]\times \C^n$ with respect to the chosen polar
coordinates and trivialization near say a positive  puncture we
require that $$ (Lh)(  s,t )(\frac{\partial }{\partial s})
=\frac{\partial h}{\partial s} - A(s)h, $$ where $A(s)\rightarrow
A_{\infty}$ for an asymptotic operator $A_{\infty}$,   as  it was
previously introduced.
\begin{theorem}
The operator $L$ is Fredholm.
\end{theorem}
 The index of $L$ can be computed in terms of Maslov indices of
 the asympotic operators (and, of course, the first Chern class of
 $E$ and the topology of $S$).

Similar to the case of closed surfaces we define the notion of
isomorphic pairs $(E,L)$ and $(F,K)$, where we emphasize the
importance of the compatibility of the asymptotic trivializations,
define the moduli space $\CR_{\punct} \supset \CR_{\cl}$ and the
Teichmuller spaces $\wt\CR_{\punct} \supset\wt\CR_{\cl}$, and
extend the determinant line bundle $\cV$ to $\CR_{\punct}$ and
$\wt\cV$ to $\wt\CR_{\punct}$. The bundle $\wt\cV$ is orientable
by the same reason as in the
  case of closed surfaces: each component of the space $\wt\CR_{\punct}$ is
contractible. However, unlike the closed case,  there is no
canonical (complex) orientation of $\wt\cV$. Still due to the
requirement
 that
isomorphisms  preserves the end structure of the operators, one
can deduce the fact that even isotopically non-trivial isomorphisms
preserve the orientation of $\wt\cV$, which shows that the bundle
$\cV$ over $\CR_{\punct}$ is orientable.

 Let us review now Axioms C1--C3 for the line bundle
$\cV$ over $\CR_\punct$. The formulation of Disjoint Union Axiom
C1
 should be appended by the  following requirement. Let $(E,L)$ and $(F,K)$ be
 operators
 over the punctured
 Riemann surfaces $(S,\bx=\{x_1,\dots,x_s\})$ and $(T,\by=\{y_1,\dots,y_t\})$, respectively.
 Then $(E,L)\dot{\cup}(F,K)$ is an operator over the surface
 $S\coprod T$ with the set of punctures ${\mathbf z}=\{x_1,\dots,x_s,y_1,\dots,y_t\}$.
 The disjoint union operation is associative, but not necessarily
 commutative (unlike the case of closed surfaces).
  Axioms C2 and C3  we formulate without any changes compared to the closed case.
 By a coherent orientation of the bundle $\cV$ over $\CR_{\punct}$
 we will mean any orientation of $\cV$ which satisfies Axioms
 C1--C3.

Take  the trivial (and globally trivialized) line bundle $E_0=
\C\times\C$ over the 1-punctured  Riemann sphere $\C=\C P^1\setminus \infty$.
For any
admissible asymptotic operator $A$ we choose a Cauchy-Riemann
operator $L^\pm_A$ on $E_0$ which has $A$ as its asymptotics at $\infty$. The
superscript $\pm$ refers to the choice of $\infty$ as the positive
or negative puncture.
 Note, that the   component  of $([E_0,L^\pm_A)$ in the  moduli spaces
  $\CR$   is uniquely determined
  by the  homotopy class $[A]$ of  the asymptotic operator $A$ in the space
  of {\it non-degenerate} asymptotic operators.

The following theorem describes all possible coherent orientations
of  the line bundle $\cV$ over $\CR$.

\begin{theorem}\label{thm:coherent}
Let us choose an  orientation $o^\pm_A$ of the operator
$(E_0,L^\pm_A)$ for a representative $A$ of each homotopy class
$[A]$ of non-degenerate  asymptotic operators. Then this choice
extends to the unique coherent orientation of the bundle $\cV$
over $\CR_{\punct}$, which coincide with the complex orientation
over $\CR_{\cl}$.
\end{theorem}
Thus  there are infinitely many coherent orientations of $\cV$
over $\CR_{\punct}$ unlike the  case of closed surfaces, when there are
precisely four.

We sketch below  the proof of Theorem \ref{thm:coherent}.
 First, similar to the case of closed surfaces, it is sufficient
 to consider only operators on the trivial, and even globally
 trivialized bundles.
 Next take the disjoint union of
 $(E_0,L^-_A)$ and $ (E_0,L^+_A)$, consider two circles $\g^\pm$ around the  punctures
 in the two copies of $\C$ and apply the cutting/pasting
 construction along these circles. As the result we get a disjoint
 union of an operator $\wt L_A$ on the trivial line  bundle over the closed Riemann
 sphere, and an operator $\overline L_A$ over  the cylinder $C=S^1\times\R$,
 which we view as the Riemann sphere
 with  two
 punctures $x_1=\infty$ and $x_2=0$ and consider $x_1$ as a positive puncture
 and $x_2$ as a negative one. The operator $\overline L_A$ has the
 same asymptotic operator $A$ at both punctures.
 Then Axioms C1 and C3 determine the orientation of $\overline
 L_A$, because for  the operator $\widetilde L_A$ we have chosen
 the complex orientation. Notice that if one glue $L_A^\pm$ in the opposite order,
 then we get an operator $\overline L_A'$ which has the reverse
  numbering of the punctures.
  The orientation of $\overline L_A'$  determined by the gluing  may be the same, or opposite as
  for the operator $\overline L_A$, depending on the parity of the Conley-Zehnder
  index of the asymptotic operator $A$.\footnote{
  The operator $\overline L_A$
 is homotopic to an isomorphism, and thus has a canonical
 orientation
 $1\otimes 1^{\ast}$. If we insist on that normalization, than
 our construction would determine the orientation of $L_A^-$ in
 terms of $L_A^+$.}

Consider now an arbitrary operator  $(E,L)$  acting on sections of a complex line bundle
$E$ over   a punctured
Riemann surface
 $(S,\bx,\bmu)$ with $\bx=\{x_1,\dots,x_s\}$, $E=S\times\C^n$, and the  asymptotic operators
 $A_1,\dots,A_s$ at the corresponding punctures. For each
 $i=1,\dots,s$ consider an operator $(E_0,L_i=L_{A_i}^\pm)$, where
 $E_0=\C\times\C^n$,  the sign $+$ is chosen if the puncture
 $x_i$ is negative, and the sign $+$ is chosen  otherwise.
 Using Axiom C1 we orient the operator $
 (E,L)\dot{\cup}(E_0,L_s)$, and then choosing circles around the
 puncture $x_s$ and $\infty$ apply the cutting/pasting procedure.
 As the result we get the disjoint union of an operator $L'$ over
 the Riemann surface with punctures $(x_1,\dots,x_{s-1})$ and the
 operator $\overline L_A$, or $\overline L_A'$ depending on
 whether the puncture $x_s$ was negative, or positive. Hence
 Axioms C1 and C3 determine the orientation of $L'$ in terms of the orientation
  of $L$. Repeating the
 procedure for the punctures $x_{s-1},\dots,x_1$ we express
  the orientation of $L$ in terms of the complex
 orientation of an operator over the closed surface.

 It remains to observe that if $E$ is a trivial complex bundle of rank $r>1$,
  then any asymptotic operator
 $A$ can be deformed through non-degenerate asymptotic operators to an operator
  $\widetilde A$ which is split into the direct
 sum of asymptotic operators on the trivial complex line. Hence we can use the direct sum
  axiom C2 to orient determinants of  operators  acting on bundles of arbitrary rank.
 \bigskip

\subsubsection{Remark about  the coherent orientation for
 asymptotic  operators with symmetries}
Let $A$ be an asymptotic operator
given by the formula (\ref{eq:asymptotic}),
where the loop $a(t),\,t\in S^1=\R/\Z$, of symmetric matrices has a symmetry
$a(t+1/2)=a(t),\;t\in\R/\Z$. Let $L$ be a Cauchy-Riemann type
operator on a bundle $E\to S$, which has $A$ as its
asymptotic
operator at a puncture $x\in S$ with an asymptotic marker $\mu$.
Let $L'$ be an operator which differs from $L$ by
rotating by the angle $\pi$ the marker $\mu$ to a marker $\mu'$,
with the corresponding change of the trivialization near the
puncture.
 Let $h:S\to S$ be  a diffeomorphism  which rotates
 the polar coordinate
 neighborhood $\cD$ of the punctures $x$ by $\pi$, and is fixed
 outside a slightly larger neighborhood.
 Then the operator $h_*L'$ has the same asymptotic data as $L$ and
 the  isomorphism classes $[E,L]$ and $[E,h_*L']$ belongs to the
 same component of the space $\wt\CR$. Given a coherent
 orientation of $\cV$, do the orientations $o_L$ and $o_{h_*L}$
 coincide? It turns out that
 \begin{lemma}\label{lm:bad-orbits}
 Let $\Psi$ be the time-one map of the linear Hamiltonian flow $\psi(t)$,
 defined by the equation {\rm (\ref{eq:asympt-flow})}. The orientations
 $o_L$ and $o_{h_*L}$ coincide if and only if the number of  real eigenvalues  of  $\Psi$
 (counted with multiplicities)  from the interval $(-1,0)$ is
 even.
 \end{lemma}
This lemma is the reason why we excluded certain periodic orbits
from $\Pc$ in Section \ref{sec:dynamics} above. See also Remarks
\ref{rem:bad-orbits1} and \ref{rem:bad-orbits2}.

\subsubsection{Coherent orientations of moduli spaces}
The moduli spaces of holomorphic curves which we need to orient are zero sets of nonlinear
Cauchy-Riemann type operators, whose linearizations are related to
operators of the kind we described (see below for more details).
In general, the moduli spaces are neither manifolds nor orbifolds,
due to the fact  that Fredholm sections cannot be made transversal
to the zero section by changing natural parameters like the almost
complex structure or the contact form. Such a transversality will
only be achievable by making abstract perturbations, leading to
virtual moduli spaces. Those virtual spaces will be the moduli
spaces which will provide us with the data for our constructions.
Nevertheless the Fredholm operators occurring in the description
of the virtual moduli spaces will only be compact perturbations of
the Cauchy-Riemann type operators, and hence the orientation scheme
for these virtual moduli spaces does not differ from the case of moduli spaces
of holomorphic curves.

A moduli space  $\Mc(\Gamma^+,\Gamma^-;W,J)$ of  holomorphic
curves in a directed symplectic cobordism $(W=\ora{V^-V^+},J)$ is
a fiber bundle over the corresponding moduli space of Riemann
surfaces. Its base is   a complex orbifold, and hence it is
 canonically oriented, while the
fiber over a point $S$, where $S$ is a Riemann surface with a
fixed conformal structure
 and positions
of punctures,  can be viewed as the space  solutions of the $\overline
\partial_J$-equation.
 If the  transversality is achieved than the  tangent bundle  of a
 moduli space  $\Mc(\Gamma^+,\Gamma^-;W,J)$  arise as  the kernel of the linearized
 surjective  operator $\overline\partial_J$.
The linearization  of $\overline
\partial_J$
  at a point $f\in \Mc(\Gamma^+,\Gamma^-;W,J)$
is a Fredholm operator in a suitable functional analytic setting.
This set-up involves Sobolev spaces with suitable asymptotic
weights derived from the non-degeneracy properties of the periodic
orbits. It is a crucial observation, again a corollary of the
behaviour near the punctures, that up to a compact perturbation,
the operator $L$ splits into two operators $L'$ and $L''$, where
$L'$ is   a complex linear operator acting on the complex line
bundle $T(S)$ of the Riemann surface $S$, and $L''$ is a
Cauchy-Riemann  type on the the bundle $E$, such that $T(S)\oplus
E=f^*(TM)$. This operator is usually only real linear, but most
importantly it is of the kind we just described in our linear
theory.
  The trivialization of $E$ near the punctures
is determined by the chosen in \ref{sec:holomorphic}
trivialization
 of the contact structure near
periodic orbit of   the Reeb vector fields on $V^\pm$, and the
asymptotic operators are determined by the linearized Reeb flow
near the periodic orbits. We have $\det L=\det L'\otimes\det L''$.  But  $\det L'$ has a canonical complex
orientation, and   hence the orientation for $\det L$ is  determined by   the orientation of
 $\det L''$. Therefore, a choice of
  a coherent orientation of $\cV$ over $\CR$ determines in the
transversal case the orientation of  all the moduli spaces
$\Mc(\Gamma^+,\Gamma^-;W,J)$.


   \subsection{First attempt of algebraization:
   Contact Floer homology}  \label{sec:Floer}
   \subsubsection{Recollection of finite-dimensional Floer theory} \label{sec:Morse}
   Let us  first recall   the basic steps in defining a Floer homology theory in the simplest case of a
Morse function $f$ on a finite-dimensional orientable closed manifold
$M$. We refer the reader
to  Floer's  original papers (see, for instance, \cite{Floer}), as well as an excellent exposition by D. Salamon  \cite{Salamon:Floer}
for the general theory.

 First,  one forms a  graded complex $C(f,g)$ generated by
critical points  $c_1,
\dots,c_N$ of      $f$, where the grading is given by the Morse index of critical points.
 Next, we  choose a  generic Riemannian metric $g$ on $M$  which satisfy the Morse-Smale condition of transversality
 of stable and unstable varieties
 of critical points. This enables us to define  a differential $d=d_{f,g}:C(f,g)\to C(f,g)$ by counting  gradient trajectories
 connecting critical points   of neighboring indices:
 $$d(c_i)=\sum L_i^jc_j,$$ where the sum is taken
 over all critical points $c_j$  with   $\ind c_j=\ind c_j-1$. The coefficient  $L_i^j$
 is the {\it  algebraic number}  of    trajectories    connecting $c_i$ and $c_j$.
  This means that the trajectories are counted
 with signs.   In the finite-dimensional case the signs could  be determined
 as follows. For each critical point we orient arbitrarily its stable manifold. Together with
 the orientation of $M$ this allows us to orient  all unstable manifolds,
  as well as the intersections
 of stable and unstable ones.   If $\ind c_j=\ind c_i-1$ then
 the stable manifold of $c_i$ and the unstable manifold
 of $c_j$ intersect along finitely many trajectories which   we want to count,
  and hence each of these trajectories
 gets an orientation. Comparing this orientation with the one given by the direction of the gradient $\nabla f$ we can associate
 with every trajectory a sign. \footnote{The generalization of this procedure to an infinite-dimensional
 case is not straightforward, because stable and unstable manifolds not only can become
 infinite-dimensional, but in most interesting cases cannot  be defined at all.  On the other hand,
  the moduli spaces of
 gradient trajectories connecting pairs of critical points  (which in the finite-dimensional case
 coincide with the intersection
 of stable and unstable manifolds of the critical points) are often defined, and one can use the coherent
 orientation scheme, similar to the one described in Section \ref{sec:orientation}
 above for the moduli spaces of holomorphic curves,
  to define their orientation.}

  To show that $d^2=0$, which then  would allow us to define the homology group
  $H_*(C(f,g),d)$, we proceed as follows.
Let us observe that  the coefficients $K_i^j$ in the expansion $d^2(c_i)=\sum K_i^jc_j$
  count  the algebraic number of broken gradient trajectories  $(\d_{il},\d_{lj})$  passing through an intermidiate
  critical point $c_l$, $l=1,\dots N$.
 But each broken trajectory    $(\d_{il},\d_{lj})$, which connects critical points whose
 indices differ by $2$, is a boundary point
 of the $1$-dimensional manifold of smooth trajectories connecting $c_i$ and $c_j$. The  algebraic number of boundary points of a compact
 $1$-dimensional  manifold is, of course, equal to $0$. Hence $K_i^j=0$, and thus $d^2=0$.

 Next we want to show that the homology group $H_*(C(f,g),d)$ is  an invariant of the manifold $M$
 (of course, in the case we consider it is just $H_*(M)$), i.e. it is independent of the choice of the function
 $f$ and the Riemannian metric $g$.  The proof of the invariance consists of three steps.
\medskip

\noindent \textsf{Step 1.}
Let us show that given a homotopy  of functions $F=\{f_t\}_{t\in[0,1]}$,   and a homotopy  of Riemannian metrics
 $G=\{g_t\}_{t\in[0,1]}$, one can define a homomorphism $\Phi=\Phi_{F,G}:C(f_1,g_1)\to C(f_0,g_0)$ which commutes with
 the boundary homomorphisms $d_0=d_{f_0,g_0}$ and $d_1=d_{f_1,g_1}$, i.e.
  \begin{equation}\label{eq:Morse-boundary}
  \Phi\circ d_1-d_0\circ\Phi=0\,.
  \end{equation}
         To construct $\Phi$ we consider the product
         $W=M\times\R $ and, assuming that the homotopies $\{f_t\}$ and $\{g_t\}$ are
         extended to all $t\in\R$  as independent  of $t$ on
         $(-\infty,-1]\cup[1,\infty)$, we
          define on $W$
        a function, still denoted by $F$,
        by the formula
         \begin{equation*}
           F(x,t)=
           \begin{cases}
           f_0(x)+ct,&  \;\;t\in (\infty,0) ;\\
           f_t(x)+ct,&  \;\;t\in[0,1]\;\;  ;\\
           f_1(x)+ct,&  \;\;t\in(0,\infty) ,\\
                        \end{cases}
           \end{equation*}
where the constant $c$ is chosen to ensure that $\frac{\partial
F}{\partial t}>0$. Similarly, we use the family of Riemannian
metrics $g_t$
to define a metric $G$ on $W$  which is equal to $g_t$ on $M\times
t$ for all $t\in\R$, and such that $\frac{\partial}{\partial t}$ is the unit vector field
orthogonal to the slices $M\times t,t\in\R$.
The gradient trajectories of $\nabla F$ converge to
critical points of $f_1$ at $+\infty$ , and  to the critical points
of $f_0$ at $-\infty$. For a generic choice of $G$ the moduli space of the
(unparameterized)
trajectories connecting two critical points, $c^1$ of $f_1$ and
$c^0 $ of $f_0$, is a compact $k$-manifold with boundary with corners,
where $k=\ind c^1-\ind c^0$. Hence,
similarly
 to the above definition of the differential $d$, we can define a homomorphism
$\Phi:C(f_1,g_1)\to C(f_0,g_0)$ by  taking an algebraic count of gradient trajectories between
the critical point of $f_1$ and $f_0$ of the same Morse index,
i.e.
$\Phi(c^1_j)=\sum \wt L^i_jc^0_j$. The identity
(\ref{eq:Morse-boundary}) comes from the  description of the  boundary
of the $1$-dimensional moduli spaces of trajectories of $\nabla
F$. Notice that the function $F$ has no critical points, and hence
a family of gradient trajectories cannot converge to a broken
trajectory in a usual sense. However, this can happen {\it at
infinity}. Let us recall that the function $F$ and the metric $G$
are cylindrical outside of  $M\times[-1,1]$. Hence
away from a compact set  a gradient
trajectory of $F$  projects to a gradient
trajectory of $f_0$ or $f_1$. When  the projection, say at $+\infty$,
 of a sequence  $\d_n:\R\to W$ of trajectories of $\nabla F$   converges  to
a broken trajectory  of    $\nabla f_1$ this can be interpreted as a
splitting at $+\infty$. This phenomenon is very similar to  the one  described
for the moduli spaces of holomorphic curves in Section
\ref{sec:compact}. Namely, there exist    gradient  trajectories $\d:\R\to W $  of $\nabla F$,
and $\d':\R\to M_1$ of $\nabla f_1$,  such that
\begin{description}
\item{$-\;$} $\d_n\to\d$ uniformly on $(-\infty,C]$ for all $C$;
\item{$-\;$} there exists a sequence $C_n\to+\infty$ such that
$\d'_n(t)=\d_n(t+C_n)$ converges to $(\d'(t),t)$ uniformly on all
subsets $[-C,\infty)$.
 \end{description}

In this sense broken trajectories of the form $(\d,\d')$ and
$(\d'',\d)$, where $\d''$ is a   trajectory of $\nabla f_0$ form
the boundary of the $1$-dimensional moduli spaces of trajectories
of $\nabla F$ connecting critical points  $c^1$ of $f_1$ and
$c^0 $ of $f_0$  with $\ind c^1-\ind c^0=1$. Therefore the
algebraic number of these trajectories equals $0$. On the other hand,
this number is equal   to $\Phi\circ d_1-d_0\circ\Phi $  which yields
the identity (\ref{eq:Morse-boundary}).
\medskip

\noindent\textsf{Step 2.} Our next goal is to check that  if
$(F_u,G_u), u\in[0,1],$  is a homotopy of homotopies which is
constant outside of a compact subset of $W$, then
  the homomorphisms $\Phi_0=\Phi_{F_0,G_0}$ and $\Phi_1=\Phi_{F_1,G_1}$ are related via the
chain homotopy formula
\begin{equation}\label{eq:Morse-chain}
\Phi_1-\Phi_0=K\circ d_1+d_0\circ K,
\end{equation}
for a homomorphism $K:C(f_1,g_1)\to C(f_0,g_0)$. The space of all homotopies $(F,G)$
connecting given pairs $(f_0,g_0)$ and $(f_1,g_1)$ is contractible,
 and hence (\ref{eq:Morse-chain}) implies
  that the homomorphism $\Phi_*:H_*(C(f_1,g_1),d_1)\to
H_*(C(f_0,g_0),d_0)$ is independent of the choice of a homotopy
$(F,G)$.

To prove (\ref{eq:Morse-chain}) one studies moduli spaces  of
gradient trajectories  of the whole $1$-parametric family of
functions $F_u$. For a generic choice of the homotopy one has
isolated critical values of the parameter $u$  when appear
{\it handle-slides}, i.e. gradient connections between critical
points with the index difference $-1$.  By  counting these
trajectories
one can then define a
homomorphism $K:C(f_1,g_1)\to\C(f_0,g_0)$ in exactly the same way as the homomorphism $\Phi$
was  defined in Step 1
   by counting   trajectories   with the index difference  $0$.

 The identity (\ref{eq:Morse-chain}) expresses  the fact that the broken
 trajectories  of the form $(\d,\d')$ and $(\d'',\d)$, where $\d$
 is a handle-slide trajectory and $\d'$ is a trajectory of $\nabla
 f_1$, form the boundary of the moduli space  of   index $0$
 trajectories in the family $(F_u,G_u)$. The difference in signs
 in formulas (\ref{eq:Morse-boundary}) and (\ref{eq:Morse-chain})
 is a reflection of the fact that the homomorphism $K$ raises the grading by
 $1$, while $\Phi$ leaves it unchanged.

\medskip
\noindent\textsf{Step 3.} Finally we need to show that
\begin{equation}
\label{eq:Morse-composition}
(\Phi_{F,G})_*=(\Phi_{F',G'})_*\circ(\Phi_{F'',G''})_*,
\end{equation}
if
 $(F,G)=\{f_t,g_t\}_{t\in[0,2]}$ is the composition of homotopies
 $(F'',G'')=\{f''_t,g''_t\}_{t\in[0,1]}$ and
 $(F',G')=\{f'_t,g'_t\}_{t\in[1,2]}$. To prove this we  view,
 as in Step 1, the  homotopy $(F,G)$ as a function and a metric
 on the cylinder $W=M\times\R$. Consider a deformation $(F_T,G_T)$ of $F,G$, by
 cutting $W$ open along $M\times 1$ and inserting a cylinder
 $M\times[0,T]$ of growing height $T$ with the function and the
 metric independent of the coordinate $t$. When $T\to+\infty$
 the gradient trajectories of $F_T$ with respect to $G_T$  split  in
 a  sense, similar to  the one explained in Step 2,
 \footnote{See
 also the discussion of  a similar phenomenon for the
 moduli spaces of holomorphic curves   in Section \ref{sec:compact}
 above.} into a ''broken trajectory"
 $(\d'',\d')$, where $\d'$  (resp.  $\d"$) is a   trajectory of $\nabla_{G'}F'$
( resp. $\nabla_{G''}F''$).  Consider the $1$-dimensional moduli  space $\Mc$ of
trajectories  of $\nabla_{G_T}F_T, T\in[0,\infty),$ connecting a
fixed
critical point  $c=c^2$ of $f_2$  with an arbitrary critical point
$c^0 $ of $f_0$  with $\ind c-\ind c^0=1$. Then the boundary of
$\Mc$ consists of
\begin{description}
\item{ a) } all the trajectories of
$\nabla_{G_0}F_0=\nabla_G F$ connecting $c^0$ and $c$; they are
given by the expression $\Phi(c)$;
\item{ b) } all the broken trajectories $(\d'',\d')$ described
above, such that $\d''$ begins at $c^0$ and ends  at a critical point
$c^1$ of $f^1$ which is, necessarily, of the same Morse index as
$c^0$ and $c^2$, $\d'$ begins at $c^1$ and ends at
$c$; these broken trajectories are described by the
expression
$\Phi_{F',G'}\big(\Phi_{F'',G''}(c)\big)$;
\item{ c) } broken trajectories defined according to Step 2 for the
$1$-dimensional family $F_T, T\in[0,\infty)$; they are described
by the expression $K(d_0(c))+d_2( K(c))$ for some homomorphism
$K:C(f_2,g_2)\to C(f_0,g_0)$.
\end{description}
 Thus the sum (taken with appropriate signs) of the three
 expressions defined in a)--c)  equals $0$, and thus we get
 \begin{equation*}
\Phi_{F',G'}\big(\Phi_{F'',G''}(c)\big)-\Phi(c )=K\circ d_0(c)+d_2\circ
K(c),
\end{equation*}
i.e. the homomorphisms $\Phi$ and
$\Phi_{F',G'}\circ\Phi_{F'',G''}$ are chain homotopic, which
yields formula (\ref{eq:Morse-composition}).

We can finish now the proof that the homology group  $H_*(C(f,g),d)$
is independent of the choice of $f$ and $g$ as follows.
 Given two pairs   $(f_0,g_0)$ and  $(f_1,g_1)$ we
 first take any homotopy $(F,G)$  connecting  $(f_0,g_0)$ with  $(f_1,g_1)$, and also take  the
 inverse homotopy $(\overline F,\overline G)$ connecting  $(f_1,g_1)$ with
 $(f_0,g_0)$.
  The composition $(\wt F,\wt G)$ of  the  homotopies $(F,G)$ and  $(\overline F,\overline G)$
  connects  the pair $(f_0,g_0)$ with itself. According to Step 3
  we have
  $(\Phi_{\wt F,\wt G})_*=(\Phi_{F,G})_*\circ(\Phi_{\overline
  F,\overline
  G})_*$. On the other hand, we have  shown in Step 3  that the
  homomorphism $(\Phi_{\wt F,\wt G})_*$ is independent of the
  choice of a homotopy, connecting $(f_0,g_0)$ with itself, and
  hence it equals the identity. Therefore, we conclude that $(\Phi_{F,G})_*$ is
 surjective, while $ (\Phi_{\overline
  F,\overline
  G})_*$ is injective. Taking the composition of homotopies  $(\overline F,\overline G)$
  and $(F,G)$ in the opposite order we prove that both
  homomorphisms are bijective.
  \bigskip

 A. Floer  discovered  that the
    finite-dimensional scheme  which we explained  in this section  works,
     modulo some
    analytic  complications, for several geometrically interesting functional on infinite-dimensional
    spaces. For instance,
  in the  symplectic Floer homology theory one deals with
  critical points of the action functional. Its critical points are periodic orbits of
  a Hamiltonian system, while for an appropriate choice of
  a metric  and an almost complex structure the  gradient trajectories    can be interpreted
  as holomorphic
  cylinders   which connect these trajectories.
  The role of broken trajectories is played here  by split holomorphic cylinders,
   and    finite-dimensional  compactness theorems are replaced by  the  highly
   non-trivial Gromov  compactness theorem for holomorphic curves.

   In the  rest of this section we explore  the  Floer-theoretic approach for
   the problem of defining invariants of contact manifolds. We will see that this approach
   works only  in a very special and restrictive   situation. However, the  general
   algebraic formalism of SFT, though quite different, has  a distinctive flavor    of
   a Floer homology theory.

\subsubsection{Floer homology for the Action functional}
   Let us make an attempt to define invariants of contact manifolds
   in the spirit of Floer homology theory. Let $(V,\xi)$     be a contact manifold
         with a fixed contact form $\alpha$ and an  almost complex structure
         $J:\xi\to\xi$,   compatible
    with the symplectic form $d\alpha|_\xi$.  Then $J$ and $d\alpha$
    define a Riemannian metric on the vector bundle $\xi$
    by the formula $g(X,Y)=d\alpha(X,JY)$ for any vectors $X,Y\in\xi$.
    We extend $g$ to the whole tangent bundle $T(V)$ by declaring the vector field $R_\alpha$ to be the unit normal
    field to $\xi$.  Consider the free loop space
    $$\Lambda(V)=\{u:S^1=\R/\Z\to V\},$$
and define the {\it action functional}
   \begin{equation}
   S:\Lambda(V)\to\R\quad\hbox{  by the formula}\quad
    S(\g) =\int\limits_\gamma \alpha .
    \end{equation}
   The least action principle tells us that the critical points of the functional $S$
   are, up to parameterization,   the
   periodic orbits of the Reeb field $R_\a$.

     The metric $g$ on $T(V)$ defines a        metric on $\Lambda(V)$
      and thus allows us to consider
    gradient trajectories of the action functional connecting critical points of $V$.
      The gradient direction $\nabla S(u)$, $u\in \Lambda(V)$, is given
by the vector field $J\pi( \frac{du}{dt})$, where $\pi:T(V)\to\xi$
is the projection along the Reeb direction, so that a gradient
trajectory $u(t,s),\,t\in\R/\Z,s\in\R,$ is given by the equation
\begin{equation}\label{eq:CR}
\frac{\partial u}{\partial s}(t,s)=J\pi\big(\frac{\partial u}{\partial
t}(t,s)\big)\,.
 \end{equation}
Equation (\ref{eq:CR}) has a flavor of a Cauchy-Riemann equation.
We want to modify it into a genuine one.
Namely, consider the Cauchy-Riemann equation
 \begin{equation*}
  \frac{\partial U}{\partial t}(t,s)=J\frac{\partial U}{\partial
  s}(t,s)
  \end{equation*}
  for $U(s,t)=   \big(u(s,t),\varphi(s,t)\big)\in V\times\R$.
  It  can be rewritten
  as a system
\begin{equation}\begin{split}\label{eq:CR2}
\frac{\partial u}{\partial s}(t,s)&=J\pi\big(\frac{\partial u}{\partial
t}(t,s)\big)+\frac{\partial \varphi}{\partial t}(t,s)R_\a(u(t,s)\\
 \frac{\partial \varphi}{\partial s}(s,t)&=-\big\langle\frac{\partial
u}{\partial t}(t,s),R_\a(u(t,s))\big\rangle\,.\\
\end{split}
 \end{equation}
Notice that $dS(\nabla S+\psi R_\a)\geq 0$ for any function $\psi(t,s)$.
 Hence,
   the first
  equation of the system (\ref{eq:CR2})  can be viewed as the flow equation of the
  gradient-like vector-field $\nabla S+\frac{\partial \varphi}{\partial
  t} R_\a$. Trajectories  of  this gradient like field   connecting critical points
  $\gamma^-,\gamma^+$ of the action functional   correspond to
   elements
  of the moduli space  $\Mc_0(\gamma^-,\gamma^+;W,J)$, and
  therefore the Floer homology philosophy (\cite{Floer}),
  which we described above in the finite-dimensional
  case, suggests the following construction.

  Let us associate a variable $q_\gamma$ with every periodic orbit $\gamma\in\Pc_\a$
  and assign to it the grading $$\deg{ q_\gamma}=\CZ(\gamma)+(n-3).$$
  The choice of the   constant $n-3$ is not important
  for purposes of this definition, but
 it will become important for generalizations considered in the second
 part of this paper.

  Let  $A$ be  the group algebra $ \C\,[H_2(V)]$. We will fix a basis
  $A_1,\dots, A_N$ of $H_2(V;\C\,)$ and identify each homology class $\sum d_iA_i$ with its {\it degree}
  $d=(d_1,\dots, d_N)$. Thus  we can view the algebra $A$ as
  the  algebra of  Laurent polynomials
   of $N$ variables $z_1,\dots,z_N$ with complex coefficients, and write its elements
   in the form $\sum a_dz^d,$ where $z^d=z_1^{d_1}\cdots
   z_N^{d_N}$. The variables $z_i$ are also considered graded,
   $\deg\,z_i=-2c_1(A_i),\,i=1,\dots,N$.
   Consider a complex $\fF$ generated by the (infinitely many) graded variables $q_\gamma$ with coefficients
  in the graded algebra $A$,
  and define  a differential $\partial: \fF\to\fF $ by the formula:
  \begin{equation} \label{eq:Floer1}
  \partial q_\gamma= \sum\limits_{\g',d}\frac{n_{\gamma,\gamma',d}}{\kappa_{\g'}}z^d  q_{\gamma'},
  \end{equation}
  where $\kappa_{\gamma'}$ denotes the multiplicity of the orbit
  $\gamma'$,
   the sum is taken over all trajectories
  $\g'\in\Pc_\a $ and $d=(d_1,\dots,d_N)$ with $$\CZ(\g')=\CZ(\g)+2\langle
  c_1,d\rangle-1,$$
   and the coefficient $n_{\gamma,\gamma',d}$ counts
  the algebraic number of components of the $0$-dimensional moduli space
  $\Mc^d_0(\gamma',\gamma;W,J)/\R
  $.\footnote{
  Let us recall that  according to our definition of the moduli
  space $\Mc^d_0(\gamma',\gamma;W,J)/\R
  $ the coefficient $n_{\gamma,\gamma',d}$ counts
  equivalence classes of holomorphic curves {\it with asymptotic
  markers}, and hence     each holomorphic cylinder connecting $\g$ and
  $\g'$ is counted $\kappa_\g\kappa_{\g'}$ times, unless the
  cylinder itself is
  multiply covered. The role of the denominators $\kappa_{\g'}$
  in formula (\ref{eq:Floer1}), as well in a similar
  formula (\ref{eq:Floer2}) below, is to correct this ``over-counting". }
  Notice that the Liouville flow of the vector field  $\frac{\partial}{\partial t}$
    defines  a $\R$-action on
  the moduli spaces    $\Mc^d_0(\gamma',\gamma;W,J)$, which makes the
   $1$-dimensional components of the moduli
  spaces
    canonically oriented.
   Comparing this orientation
  with the coherent orientation we produce  signs   which we use in the formula
  (\ref{eq:Floer1}).

To simplify the  assumptions in the propositions which we formulate below
we will  assume for the rest of this section that $c_1|_{\pi_2(V)}=0$.
 This assumption allows us to define for any
 contractible periodic orbit $\g$ the Conley-Zehnder index
 $\CZ_{\disk}(\g)$ computed with respect to {\it any disk}
 $\Delta$ spanned by $\g$ in $V$. We denote
 $\deg_{\disk}(\g)=\CZ_{\disk}(\g)+n-3$.

  \begin{proposition}\label{prop:Floer1}
  If for a  contact form $\alpha$ the Reeb
  field $R_\alpha$ has no contractible periodic orbits $\g\in\Pc_\a$  with
  $\deg_{\disk}(\g)=1$, then
  $\partial^2=0$.
  \end{proposition}

{\sl Sketch of the proof}. Similarly to the finite-dimensional
case considered in Section \ref{sec:Morse} above the identity
$\partial^2=0$ in Floer homology  is equivalent to the fact that
the codimension $1$ stratum of the compactified moduli spaces
$\Mc^d_0(\gamma',\gamma)$   consists of broken
trajectories, which in our case are represented by the height $2$
stable curves $(f_1,f_2)$, $f_1\in\Mc^{d'}_0(\gamma',\gamma'')/\R$, $f_2\in
\Mc^{d''}_0(\gamma'',\gamma)/\R$, where $d=d'+d''$. However, in  the general case a sequence of
holomorphic cylinders in $\Mc^d_0(\gamma',\gamma)$ can split into curves
different from cylinders, as it is stated in Proposition \ref{prop:boundary} and  Corollary
\ref{cor:pos-puncture}. But
 if this happens then the  first-floor curve
$f_1$ must have  a component which is conformally equivalent to $\C\,$ and
asymptotically cylindrical over a  contractible orbit at
$+\infty$. Moreover, if $(f_1,f_2)$ belongs to a top-dimensional
stratum of the boundary of the moduli space
$\Mc^d_0(\gamma',\gamma)$, then $\deg_{\disk}(\g)=1$, which
contradicts our assumption.

\begin{remark}\label{rem:bad-orbits1}{\rm
Let us recall that we excluded from $\Pc$ certain ``bad" periodic
orbits (see  the footnote in Section \ref{sec:dynamics}).
However on the boundary of the moduli space
$\Mc^d_0(\gamma',\gamma)$  there could be a stratum
 which consists of height $2$
stable curves $(f_1,f_2)$, $f_1\in\Mc^{d'}_0(\gamma',\gamma'')$, $f_2\in
\Mc^{d''}_0(\gamma'',\gamma')$, where the orbit $\gamma''$ is
one of  the bad orbits which we excluded from $\Pc$. The
orbit $\g''$ has even multiplicity $2k$, and hence on the boundary
of   $\Mc^d_0(\gamma',\gamma)$ there are $2k$ strata which
correspond
to
$2k$ different possible positions of the asymptotic marker at
the punctures mapped to $\g''$. The Poincar\'e return map of the Reeb flow
along the orbit $\g''$ has an odd number of eigenvalues in the interval $(-1,0)$,
and hence  according to Lemma \ref{lm:bad-orbits}  the coherent orientation
will automatically assign to these orbits opposite signs, which means that these strata
will not contribute to the sum  (\ref{eq:Floer1}). This explains why the exclusion of bad orbits
is {\em possible}. Remark \ref{rem:bad-orbits2} below explains why this
exclusion is {\em necessary}.}
\end{remark}

Now we follow Steps 1--3 in Section \ref{sec:Morse} above to show the independence
of   the homology group   $$\oplus H_k(\fF ,\partial )= \Ker\partial /{\mathrm {Im}}\partial, $$
    graded by the degree $k$,  of the choice of
  a  nice contact form $\a$ and  a compatible almost complex
  structure  $J$.

Suppose now that we have a  directed symplectic cobordism
$W=\ora{V^-V^+}$,  and $J$  is a compatible almost complex structure on $W$.
Suppose that the inclusions $V^\pm\hookrightarrow W$ induce
isomorphisms on $2$-dimensional homology.  Then we can define
  a homomorphism $\Phi=\Phi_W:\fF^+\to\fF^-$ by the formula
\begin{equation}\label{eq:Floer2}
\Phi( q_\gamma)= \sum\limits_{\g',d} \frac{1}{\kappa_{\gamma' }}n_{\gamma,\gamma',d}z^dq_{\gamma'},
\end{equation}
  where the sum is taken over all trajectories $\g'\in\Pc^- $  and $d$
  with
  $\CZ(\g')=\CZ(\g)+2\langle c_1,d\rangle$, and
  the coefficient $n_{\gamma,\gamma',d}$ counts
  the algebraic number of points of the compact $0$-dimensional moduli space $\Mc^d_0(\gamma',\gamma;W,J)
  $.
   If the   condition on the second homology is not satisfied then the above construction
  gives us only a correspondence, rather than
  a homomorphism. See Section \ref{sec:composition}  for the discussion of a
   more general case.

\begin{proposition}\label{prop:Floer2}
Suppose that the contact forms   $\a^\pm$ associated to the
ends satisfy the condition  $\deg_{\disk}(\g)\neq 0,1$ for any {\em contractible in $W$} periodic orbit
 $\g\in\Pc^\pm$. Then the homomorphism $\Phi_W$ commutes
with $\partial $.
\end{proposition}
\begin{proposition}\label{prop:Floer3}
 Let $J_t,t\in[0,1] $, be a family of  almost complex structures compatible with
 the directed symplectic cobordism
$W=\ora{V^-V^+}$. Suppose  that the forms   $\a^\pm$ associated to the
ends satisfies the condition  $\deg_{\disk}(\g)\neq -1,0,1$ for any  {\em contractible in $W$} periodic orbit
 $\g\in\Pc^\pm$. Then the  homomorphisms
 $\Phi_0=\Phi_{W,J_0}$ and  $\Phi_1=\Phi_{W,J_1}$  are chain
 homotopic, i.e there exists a homomorphism $\Delta:\fF^+\to\fF^-$
 such that $\Phi_1-\Phi_0=\partial \Delta +\Delta\partial$.
 \end{proposition}
\begin{proposition}\label{prop:Floer4}
Given two cobordisms $W_1$ and $W_2$, and a compatible almost complex structure
$J$ on the composition $W_1\circledcirc W_2$, the homomorphism
$\Phi_{W_1\circledcirc W_2}$ is chain-homotopic to
$\Phi_{W_1}\circ\Phi_{W_2}$.
\end{proposition}

Together with an obvious remark that for the cylindrical cobordism   $W_0$
the homomorphism $\Phi_{W_0}$ is the identity,
  Propositions \ref{prop:Floer1}--\ref{prop:Floer3}
  imply that if a contact structure $\xi$ on $V$ admits a {\it
  nice} contact form, i.e. a form without contractible periodic
  orbits of index $-1,0$ and $1$, then the {\it contact homology
  group}
  $$\oplus HC_k(V,\xi)=\oplus H_k(\fF ,\partial ) $$
   is well defined and
    independent of the choice of
  a  nice contact form and  a compatible almost complex
  structure
  (however, if $H_2(V)\neq 0$ and/or $H_1(V)\neq0$ it depends on a
  choice of spanning surfaces $F_\g$ and the  framing of the bundle
  $\xi$ over basic loops).
   Similarly to what was explained in the sketch of
     the proof  Proposition
      \ref{prop:Floer1} the  ``niceness" assumptions   guarantees
      that the  top codimension strata   on the boundary of the involved
      moduli spaces consist of height 2 cylindrical curves,
      and thus         the proofs of Propositions  \ref{prop:Floer2}--\ref{prop:Floer4}
      may precisely follow the standard scheme of the Floer
      theory (see \cite{Floer, Salamon:Floer}).
\begin{remark}\label{rem:bad-orbits2}{\rm
Similarly to what we explained in Remark \ref{rem:bad-orbits1} the
coefficient $n_{\g,\g'}$ in the definition (\ref{eq:Floer2}) of
$\Phi$ equals $0$ if at least one of the orbits $\g,\g'$ is ``bad''.
Hence, in the presence of ``bad" orbits the
homomorphism $\Phi$  could never be  equal to the  identity, even
for the cylindrical cobordism. This explains why
the exclusion of ``bad" periodic orbits is {\em necessary}.}
\end{remark}

    Besides the degree (or Conley-Zehnder)  grading,
   the contact homology group is graded by elements of $H_1(V)$, because
   the boundary operator preserves the homology class of a periodic orbit.  We will denote the part of  $HC_*(V,\xi)$
   which correspond to a class $a\in H_1(V)$ by $HC_*(V,\xi|a)$.
  One can  similarly construct a contact homology group
  $HC^{\contr}_*(V,\xi)$,
  generated  only by contractible periodic orbits,  which is another invariant
  of
  the contact manifold $(V,\xi)$.

Contact structures which admit nice contact forms do exist, as it
is
illustrated by examples in Section \ref{ex:Floer}  below. However, the
condition of existence of a nice form is too restrictive. The
general case leads to an algebraic formalism developed in
Sections \ref{sec:af-contact}--\ref{sec:composition} below.

\subsubsection{Examples}\label{ex:Floer}

\noindent 1. \textsf { Contact homology of the standard contact sphere $S^{2n-1}$.}

 Take the $1$-form $\alpha=\frac 12\sum (x_idy_i-y_idx_i)$,
 which is a primitive of the standard symplectic structure
 in $\R^{2n}$. Its restriction  to a generic  ellipsoid
 $$S=\{\sum\frac{x_i^2+y_i^2}{a_i^2}=1\}$$ is a  nice contact form for the standard contact structure
 $\xi$ on the sphere $S=S^{2n-1}$. The form $\alpha|_S$ has
 precisely one periodic orbit for each  Conley-Zehnder index
 $n +2i-1$ for $i=1,\dots,$. Hence the  contact homology group
 $HC_*(S,\xi)$ has one generator in each  dimension $2i, i\geq n-1$.
 See also the discussion in Section \ref{sec:Bott} below.
 \medskip

\noindent 2. \textsf{ Contact homology of Brieskorn spheres.}

 Ilya Ustilovsky  computed (\cite{Ustilovsky}) the
 contact homology of certain Brieskorn spheres.

 Let us  consider the Brieskorn manifold
 \begin{equation*}
 \Sigma(p,\underbrace{2,\dots,2}_n)=
 \{z_0^p+\sum\limits_1^n  z_j^2=0\}\cap \{\sum\limits_0^n|z_j|^2=1\}\subset \C^{n+1}.
 \end{equation*}
 $\Sigma(p,\underbrace{2,\dots,2}_n)$ carries a
  canonical contact structure  as a strictly pseudo-convex hypersurface in a complex manifold.

 Suppose that $n=2m+1$ is odd, and $p\equiv 1\,\mod\,8.$
 Under this assumption  $\Sigma(p,\underbrace{2,\dots,2}_n)$
 is diffeomorphic to $S^{2n-1}$ (see \cite{Brieskorn}).
  However, the following theorem of Ustilovsky implies that     the contact structures on
   Brieskorn spheres
 $\Sigma(p,\underbrace{2,\dots,2}_n)$ and $\Sigma(p',\underbrace{2,\dots,2}_n)$ are not isomorphic, unless $p=p'$.
This result should be confronted with a computation of Morita
(\cite{Morita} ), which implies
 that  the formal homotopy class (see Section \ref{sec:prelim} above) of the contact structure on $\Sigma(p,\underbrace{2,\dots,2}_n)$
  is   standard, provided $p\equiv 1\,\mod\, 2(2m!)$. Hence, Ustilovsky's theorem provides infinitely many non-isomorphic contact
  structures on $S^{4m+1}$ in the standard formal homotopy class.

  \begin{theorem} {\rm (I. Ustilovsky, \cite{Ustilovsky})}
  The contact homology $$HC_*\left(\Sigma(p,\underbrace{2,\dots,2}_n)\right)$$ is defined, and
  the dimension  $$c_k=\dim\,HC_k\left( \Sigma(p,\underbrace{2,\dots,2}_n)\right) $$
  is given by the formula
  \begin{equation*}
  c_k  =
  \begin{cases}
  0,& k \;\;\hbox {is odd or}\;\;k<2n-4,\\
  2,& k=2\left[ \frac{2N}{p}\right] +2(N+1)(n-2),\;\hbox{for  }\;\;
  N\geq 1,\,2N+1\notin p\Z,\\
  1,& \hbox{in all other cases.}\\
  \end{cases}
  \end{equation*}
  \end{theorem}

\medskip

\noindent 3. \textsf{ Contact homology of boundaries
of subcritical Stein manifolds.}

 A co-oriented contact manifold $(V,\xi)$ is called Stein-fillable if it can be realized as a strictly pseudoconvex
 boundary of a complex manifold $W$,
 whose interior is Stein, and if the co-oriented contact structure $\xi$ coincides with the canonical contact structure
 of a strictly pseudo-convex hypersurface.   We say that $(V,\xi)$ admits a
 subcritical Stein filling if the corresponding Stein manifold $\Int W$ admits an exhausting plurisubharmonic function
 without critical points of dimension $\dim_\C(W)$. If $\dim V>3$ then one can equivalently require
 that $W$ deformation  retracts to   a CW-complex of dimension $<\dim_\C W$ (see \cite{Eliash-Stein}).

 Mei-Lin Yau studied in her  PhD thesis \cite{MLYau} contact homology of contact manifolds
 admitting a subcritical Stein filling.
 Here is her result.
 \begin{theorem}\label{prop:MLYau} {\rm (Mei-Lin Yau,  \cite{MLYau})}
 Let $(V,\xi)$ be a  contact manifold of dimension $2n-1$ which admits a subcritical Stein filling $W$.
 Suppose that $c_1(V)=0$ and $H_1(V)=0$. Let $c_1,\dots, c_k$ be generators of $H_*(W)$.
  Then  the contact homology $HC_*(V)$ is defined  and generated
  by elements
  $ q_{i,j}$ of  degree $\deg q_{i,j}= 2(n+i-2)-\dim c_j$, where  $j=1,\dots, k$, and
  $i\geq 1$.
  \end{theorem}

\medskip

\noindent 4. \textsf{ Contact homology of $T^3$ and its coverings.}
\medskip

 Set $\a_n=\cos 2\pi nz\,dx+\sin 2\pi nz\,dy$. This contact form descend
 to the $3$-torus $T^3=\R^3/\Z^3$ and defines there a contact structure $\xi_n$.
 The structure $\xi_1$ is  just the canonical contact structure on $T^3$ as the space
 of co-oriented contact elements of $T^2$. The form  $\a_n$ for $n>1$  is
 equal to the pull-back $\pi^*_n (\a_1)$, where $\pi_n:T^3\to T^3$
 is the covering
 $(x,y,z)\mapsto(x,y,nz)$. Notice that all structures $\xi_n$
 are homotopic as plane field to the foliation $dz=0$.

 \begin{theorem}\label{thm:T3}
 The contact homology group $HC_*(T^3,\xi_n|w)$, where $w$ is the homology class $(p,q,0)\in H_1(T^3)$,
 is isomorphic to $\Z^{2n}$.
 \end{theorem}
 In particular we get as a corollary a theorem of E. Giroux:
 \begin{corollary} {\rm (E.  Giroux, \cite{Giroux})}
 The contact structures $\xi_n$, $n=1,\dots,$ are pairwise non-isomorphic.
 \end{corollary}

  The contact manifold $(T^3, \xi_1)$
 is foliated by pre-Lagrangian
 tori $L_{p,q}$, indexed by simple homology classes $ (p,q)\in H_1(T^2)$. Each torus
 $L_{p,q}$ is foliated by  the $S^1$-family of lifts of
  closed geodesics which represent the class $(p,q)$.
  Thus for any given $ (p,q)\in H_1(T^2)$ (even when  $(p,q)$ have common divisors)
   the set of closed orbits
  in  $\Pc_{\a_1}$ which represent the class $(p,q,0)\in H_1(T^3)$
   is  a circle $S_{p,q}$, and for any $n\geq 1$ the set of closed orbits
  in  $\Pc_{\a_n}$ which represent the class $(p,q,0)\in H_1(T^3)$
  consists of $n$ copies $S^1_{p,q},\dots,  S^n_{p,q}$ of such circles.
   The forms $\a_n$ have no contractible periodic orbits, but of course, they
    are degenerate. To compute the contact homology groups, one can   either work directly with these
    degenerate forms, as it is explained in Section \ref{sec:Bott} below, and
    show that  $HC_*(T^3,\xi_n|w)=H_*(\bigcup\limits_1^n S^i_{p,q})=\Z^{2n}$, or  first  perturb the form $\alpha_1$,
  and respectively all its covering forms $\a_n= \pi_n^*(\a_1)$, in order to substitute each circle $S^i_{p,q}$
  by two non-degenerate periodic orbits, and then show that the orbits from  each of these pairs  are connected by precisely
  two holomorphic cylinders, which cancel each other in the formula for the boundary operator $\partial$.

\subsubsection {Relative contact homology and contact non-squeezing theorems}
\label{sec:squeezing}
Let us observe that  the  complex $(\fF,\partial)$  is filtrated
by the values of the action  functional $S$, $\fF=\mathop{\bigcup}\limits_{a\in\R}\fF^a$, where
the complex $\fF^a$
is generated  by variables $q_\gamma$ with $S(\g)\leq a$. The differential
$\partial$ respects this filtration, and hence descends to $\fF^b/\fF^a, a<b$.
 Hence, one can define  the homology
$H_*^{(a,b]}(\fF,\partial)=H_*(\fF^b/\fF^a,\partial)$ in the window $(a,b]\subset\R$. Of
course, $H_*^{(a,b]}$ depends on a choice of a particular nice form
$\alpha$. If $\alpha>\beta$ then we  have a map $\Phi_*:H_*^{(a,b]}(\fF,\partial;\alpha)\to
H_*^{(a,b]}(\fF,\partial;\beta)$. We write $H^a$ instead of $H^{(-\infty,a]}$.

Consider now a contact manifold $(V,\xi)$ which is
 either closed, or satisfies
the following
  {\it pseudo-convexity}  condition at infinity.
A contact manifold $(V,\xi=\{\a=0\})$  with a fixed contact form $a$
is called pseudo-convex at infinity
if there exists a compatible almost complex structure
 $J$ on the symplectization $V\times\R$ for which $V$ can be exhausted
 by compact domains $V_i$ with smooth  pseudo-convex
 boundary. A sufficient condition for pseudo-convexity is existence of
  an exhaustion $V=\bigcup V_i$, such that for each $i=1,\dots, $
 trajectories of the Reeb field
 $R_\a|_{V_i}$ do not have interior tangency points with $\partial  V_i$.  For instance, for
 the standard contact form  $\a=dz-\sum y_idx_i$
 on $\R^{2n+1}$ the latter condition is satisfied for an exhaustion of $\R^{2n+1}$
 by round balls, and  ence the standard contact form on $\R^{2n+1}$ is pseudo-convex at
  infinity.

Our goal is to define  a relative contact homology group $HC_*(V,U,\xi)$
 for a relatively compact open subset
$U\subset V$, so that  this group would be invariant    under a contact isotopy of $U$ in $V$.

Let us fix a contact form $\alpha$ on $V$  which satisfies the above pseudo-convexity
 condition.
Let us denote by $\Fc_{U,\a}$ the set of  $C^\infty$-functions $f:V\to[0,\infty)$
which  are $\leq 1$ on   $U$,   and for which  the contact form
$f\alpha$ is nice and pseudo-convex at infinity.\footnote{ Of course, the set
 $\Fc_{U,\a}$ may be
empty, because the niceness condition is very restrictive.
In this case one needs to employ
a more general construction from Section \ref{sec:3algebras}.}
 Take a strictly increasing sequence
of functions
$f_i\in \Fc_{U,\a}$, such that
\begin{description}
\item{a)}
$\max\limits_K f_i \mathop{\to}\limits_{i\to\infty}\infty$ for   each compact
 set $K\subset (V\setminus
\overline{U})$;
\item{b)} $f_i|_U \mathop{\to}\limits_{i\to\infty} 1$ uniformly on
compact sets.
\end{description}

\begin{proposition}\label{prop:rel-homol}
The limit
\begin{equation}
HC_*(V,U,\xi)=\lim\limits_{a\to+\infty}
\lim\limits_{\leftarrow}HC^a_*(V,f_i\alpha)
\end{equation}
is independent of $\alpha$,
 and
thus it is  an invariant of  the contact pair $(V,U)$.
A contact isotopy  $f_t:V\to V$ induces a  family of
isomorphisms $(f_t)_*:H_*(V,U)\to H_*(V, f_t(U))$.
An inclusion $i:U_1\mapsto U_2$ induces a homomorphism
\begin{equation*}
i_*: HC_*(V,U_1,\xi)\to HC_*(V,U_2,\xi).
\end{equation*}
\end{proposition}

\bigskip

One of the most celebrated results in Symplectic topology is
Gromov's non-squeezing theorem which  states that one cannot
symplectically embed a  $2n$-ball of radius $1$ into
$\D^2_r\times\R^{2n-2}$ for $r<1$. Here $D^2_r$ denotes a
 $2$-disk of radius $r$ and $\D^2_r\times\R^{2n-2}$ is endowed
 with the product of standard symplectic structures.
  Because of the conformal character of contact geometry one
  cannot expect as strong non-squeezing results
  for
  contact manifolds: one can embed any domain in the standard
  $\R^{2n-1}$ in an arbitrary small ball. However, it turns out
  that it is not always possible  to realize a contact squeezing via a contact
  isotopy inside a manifold with  a non-trivial first Betti number.

As an example, consider the $1$-jet bundle $$V=J^1(\R^n,S^1)=T^*(\R^n)\times S^1$$
 of
$S^1$-valued functions with its standard contact structure $\xi$,
given by the contact form $\alpha=dz-\sum y_idx_i$, $(x,y)\in
\R^{2n}=T^*(\R^{n}),\, z\in S^1 =\R/\Z$.
The contact form $\alpha$ satisfies the condition
of  pseudo-convexity   at infinity
and it is nice: the Reeb field equals $\frac{\partial}{\partial z}$,
and thus it has no contractible periodic orbits.
 Let us consider the class $\cN$ of domains $\Omega\subset V$
 which are images of the split domains $U\times S^1\subset V$
 under a contact isotopy of $V$, where $U$ is connected.
For any $\Omega\in\cN$ the relative contact  homology
group $HC_*(V,\Omega)$ is  well defined
 because   for any function $f:\R^{2n}\to\R$
 the form $f(x,y)\alpha$ is nice.

Let us denote by $\cE_r(\Omega)$, $\Omega\in\cN$,
the space of contact embeddings $ D_r\times S^1\to
 \Omega\times   S^1$, contact
 isotopic in $V$  to the inclusion
 $$D_r\times S^1\hookrightarrow \R^{2n}\times S^1=V\,.$$
Notice, that for any two embeddings $f,g\in\cE_r(\Omega)$
there exists  a positive $\rho\leq r$, such that the
restrictions
$f|_{D_\rho\times S^1}$ and $g|_{D_\rho\times S^1}$
are isotopic       via a {\it contact} isotopy.

Given a  contact embedding $f\in\cE_r(\Omega)$
 there is defined a homomorphism
\begin{equation*}
f_*: HC_*(V,D_r\times S^1,\xi)\to HC_*(V,\Omega,\xi)\,.
\end{equation*}
  Let us choose a
    symplectic trivialization
  of the contact bundle $\xi$ induced by the projection
  $V\to\R^{2n}$. We will assume that indices of periodic orbits,
  and hence the grading of contact homology groups,
  are associated with this trivialization.

For each homology class $k\in \Z=H_1(D_r\times S^1)$ let us
consider the maximal $l=l(f,k)$ such that
$$\Ker\big( f_*|_{HC_l(V,D_r\times S^1,\xi|k)}\}\neq 0 ,,$$
    and define an invariant $w_{\cont}(V,\Omega)$, called
    the relative {\it
contact width} by the formula
\begin{equation}
w_{\cont}(V,\Omega)=\sup \limits_{k,r>0,f\in
\cE_r(U)}\frac{2 k}{l(f,k)}\;.
\end{equation}

  S.-S. Kim has computed $w_{\cont}(V,\Omega)$ for
  certain domains $\Omega$. In particular,
   she proved
\begin{proposition}\label{prop:Kim}
\begin{eqnarray*}
w_{\cont}(V,D^{2n}_r\times S^1)= \pi r^2;  \\
w_{\cont}(V,D^2_r\times D^{2n-2}_R\times S^1) =\pi r^2,
\end{eqnarray*}
if $R\geq r$.
\end{proposition}
    The contact width is clearly a monotone invariant, i.e.
$$ w_{\cont}(V,U_1\times S^1)\leq w_\cont(V,U_2\times S^1)$$
if $U_1\subset U_2$. Hence Proposition \ref{prop:Kim} implies

  \begin{corollary}\label{cor:non-squeezing}
  Suppose that $r<\min(r',R).$ Then there is no  contact isotopy
   $f_t:D^{2n}_{r'}\times S^1\to V$ such that $f_0$
   is the inclusion, and
  $$f_1( D^{2n}_{r'}\times S^1)\subset D^2_{r}\times D_R^{2n-2}\times S^1 .$$
  \end{corollary}

     \begin{problem}  Suppose there exists a contact isotopy
      $f_t:V=\R^{2n}\times S^1\to V$ with
     $f_0=\Id$ and $f_1(U_1\times S^1)\subset U_2\times S^1$. Does
     there exist
     a Hamiltonian isotopy $g_t:\R^{2n}\to\R^{2n}$ such that
      $g_0=\Id$ and $g_1(U_1)\subset U_2$?
      \end{problem}
      Notice that the converse is obviously true.

\section{Algebraic formalism}\label{sec:algebraic}

\subsection{Informal introduction}\label{sec:informal}

The Floer-theoretic formalism described in Section \ref{sec:Floer}
 is applicable only to a very limited class
of contact manifolds. As it follows from Theorem \ref{thm:comp1}
the boundary of the moduli space   of holomorphic cylinders in the
symplectization may   consist of stable curves, different from broken
cylinders; for instance,  it may contain    height 2 stable curves
which consist  of a pair of pants on the upper level, and a copy
of $\C$ plus a trivial cylinder at the bottom one. Hence the
minimal class of holomorphic curves
 in symplectizations which
has the property that the  stable curves  of height $>1$ on the boundary
of the corresponding moduli space  are made of curves from the same class,
must contain all  rational curves with one positive   and an arbitrary
 number of  negative punctures.
 The counting of curves  with one positive and arbitrary number of negative punctures
 can still be interpreted as a differential, but this time defined
 not on the {\it vector space} generated by    periodic trajectories but on the   graded {\it algebra}
  which they
 generate. Thus
 this  leads to a straightforward generalization of the Floer type formalism considered
   in    Section \ref{sec:Floer}
 when instead of  the additive  Floer complex   $\fF$  generated by the
 variables $q_\g$,
 one considers a graded commutative algebra  $\fA$
 generated by these variables, and when  instead of the  formula (\ref{eq:Floer1})
  the differential  $\partial q_\g$ is defined
 as a polynomial of  a higher degree, rather than a linear
 expression as in the Floer homology case. Namely, we define
 \begin{equation}\label{eq:d}
\partial q_\g= \sum \frac{n_{\Gamma,I,d}}{k!\prod^k_1
i_j!\kappa_{\g_j}^{i_j}}q_{\g_1}^{i_1}\dots q_{\g_k}^{i_k}z^d,
   \end{equation}
   where the sum is taken over all ordered
   \footnote{The coefficient $\frac{1}{k!}$ is the price we pay
   for taking {\it ordered} sets of periodic orbits.} sets of
   different periodic orbits $\Gamma=\{\g_1,\dots,\g_k\} $,   multi-indices
$d=(d_1,\dots, d_N)$ and $I=(i_1,\dots,i_k), i_j\geq 0, $
 and where the coefficient $n_{\Gamma,I,d}$  counts
the algebraic number of elements of the moduli space
$$\Mc^d_0(\g; \underbrace{\g_1,\dots,\g_1}_{i_1},\dots,\underbrace{\g_k,\dots,\g_k}_{i_k})/\R,$$ if
this space is $0$-dimensional, and equals $0$ otherwise.
The differential $\partial$ extends to the algebra $\fF$ according
to the graded Leibnitz rule.
Roughly speaking, $\partial q_\g$ is a polynomial,
whose monomomials $q_{\g_1}\dots q_{\g_l}$ are in 1-1
correspondence with rigid, up to translation,   rational
holomorphic   curves with one positive  cylindrical end
over $\g$ and $l$ negative cylindrical ends over
trajectories $\g_1,\dots,\g_l$.

 It turns out that the  quasi-isomorphism class
 of the differential algebra $(\fA,\partial)$
 is  independent of all extra choices (see Section \ref{sec:3algebras} below).
 In particular, the {\it contact homology algebra} $ H_*(\fA,\partial)$
  is an invariant of the contact manifold
 $(M,\xi)$.

 Having included into the picture the moduli spaces of rational
 curves with one positive and several negative punctures, one may
 wonder, what is the role of rational curves with an arbitrary
 number of positive and negative punctures. One can try to
 interpret the
 counting of rational holomorphic curves with
fixed number of  positive and an arbitrary number of negative
punctures  as a sequence of bracket type  operations on the
algebra $\fA$. These operations  satisfy an infinite system of
indentities, which remind the formalism of homotopy Lie algebras.

 However,  there is a more adequate algebraic formalism
 for this picture. Let us associate   with  each periodic
  orbit $\g$ two graded variables $p_\g$ and $q_\g$ of the
  same parity (but   of
  different integer grading, as we will see in Section
  \ref{sec:correlators} below),
  and consider an algebra $\fP$ of formal power series
  $\sum f_\Gamma(q)p^\Gamma$, where $f_\Gamma(q)$
    are polynomials of $q=\{q_\g\}$ with coefficients
  in (a completion of) the group algebra  of $H_2(V)$.
  It is useful to think about the algebra $\fP$ as the graded
  Poisson
   algebra of functions
  on the     infinite-dimensional symplectic super-space $\bf{V}$ with the even
  symplectic form $\sum\limits_{\g\in\Pc} {\kappa_\g}^{-1} dp_\g\wedge dq_\g$, or rather on
  its formal analog along the $0$-section $\{p=\{p_\g\}=0\}$.
  With each $0$-dimensional moduli space
  $\Mc^d_0(\Gamma^-,\Gamma^+)/\R$, $$\Gamma^\pm=\{\underbrace{\g^\pm_1,\dots\g^\pm_1}_{i_1^\pm},\dots,\underbrace{\g^\pm_{s^\pm}\dots
  \g^\pm_{s^\pm}}_{i^\pm_{s^\pm}}\},$$
  we associate a monomial
  $$\frac{n_{\Gamma^-,\Gamma^+,d}}{s^-!s^+!(i_1^-)!\dots(i_{s^-}^-)!(i_1^+)!\dots(i_{s^+}^+)!}
    q_{\g^-}^{I^-}p_{\g^+}^{I^+}z^d,$$
  where $q_{\g^-}^{I^-}=(q_{\g^-_{1}})^{i^-_{1}}\dots(q_{\g^-_{s^-}})^{i^-_{s^-}}$, $p_{\g^+}^{I^+}=(p_{\g^+_1})^{i^+_1}\cdots
  (p_{\g^+_{s^+}})^{i^+_{s^+}}$,  and $n_{\Gamma^-,\Gamma^+,d}$ is
   the algebraic number of elements of the moduli
  space     $\Mc^d_0(\Gamma^-,\Gamma^+)/\R$.

The sum of all these monomials over all $1$-dimensional moduli
spaces    $\Mc^d_0(\Gamma^-,\Gamma^+) $ for   all ordered sets
    $\Gamma^-,\Gamma^+$ of periodic orbits
is  an {\it odd} element $\bh\in\fP$. All the operations on the
algebra $\fA$ which we mentioned above appear as the expansion
terms of $\bh$ with respect to $p$-variables. It turns out that
the   infinite system of  identities  for
operations on $\fA$, which we mentioned above, and which is
 defined by counting holomorphic curves with  a
certain fixed number of positive punctures, can be encoded into a single
equation $\{\bh,\bh\}=0$. Then the differentiation with respect to
the Hamiltonian vector field, defined by the Hamiltonian $\bh$:
\begin{equation*}
d^{\bh}(g)=\{\bh,g\},\;g\in\fP,
\end{equation*}
defines a differential $d=d^{\bh}:\fP\to\fP$, which satisfies the
equation $d^2=0$. Thus one can define the homology $H_*(\fP,d^{\bh})$
which inherits the structure of a graded Poisson algebra.

 The identities, like  $d^2=0$ and $\partial^2=0$,
 encode in algebraic terms information about the structure of the  top-dimensional strata on the boundary of
 compactified moduli spaces of holomorphic curves,
as it is described in Proposition \ref{prop:boundary} above. For
instance, the codimension $1$ strata on the boundary of the moduli
space $\Mc^d_0(\Gamma^-,\Gamma^+)/\R$ consists of height two stable
rational curves $(f_1,f_2)$.  Each floor of
 this curve  may be disconnected,
but  precisely one of its components differs from the straight cylinder.
  Each   connected component of $f_1$  can be glued with a component of $f_1$ only
along one of their ends. One can easily see that the combinatorics  of  such gluing   precisely corresponds to
the Poisson bracket formalism and that the algebraic sum of the monomials associated to all
stable curves of height two equals $\{\bh,\bh\}$. On the other hand,
 the algebraic number of such height 2 curves equals
 $0$ because they form   the  boundary
of the a compactified  $1$-dimensional moduli space of holomorphic curves.
Hence we
  get the identity $\{\bh,\bh\}=0$. The identity is
not tautological due to   the fact that $\bh$ is odd.   In view of the super-Jacobi identity
it is equivalent to the identity $(d^{\bh})^2=0$.

\begin{figure}
\centerline{\psfig{figure=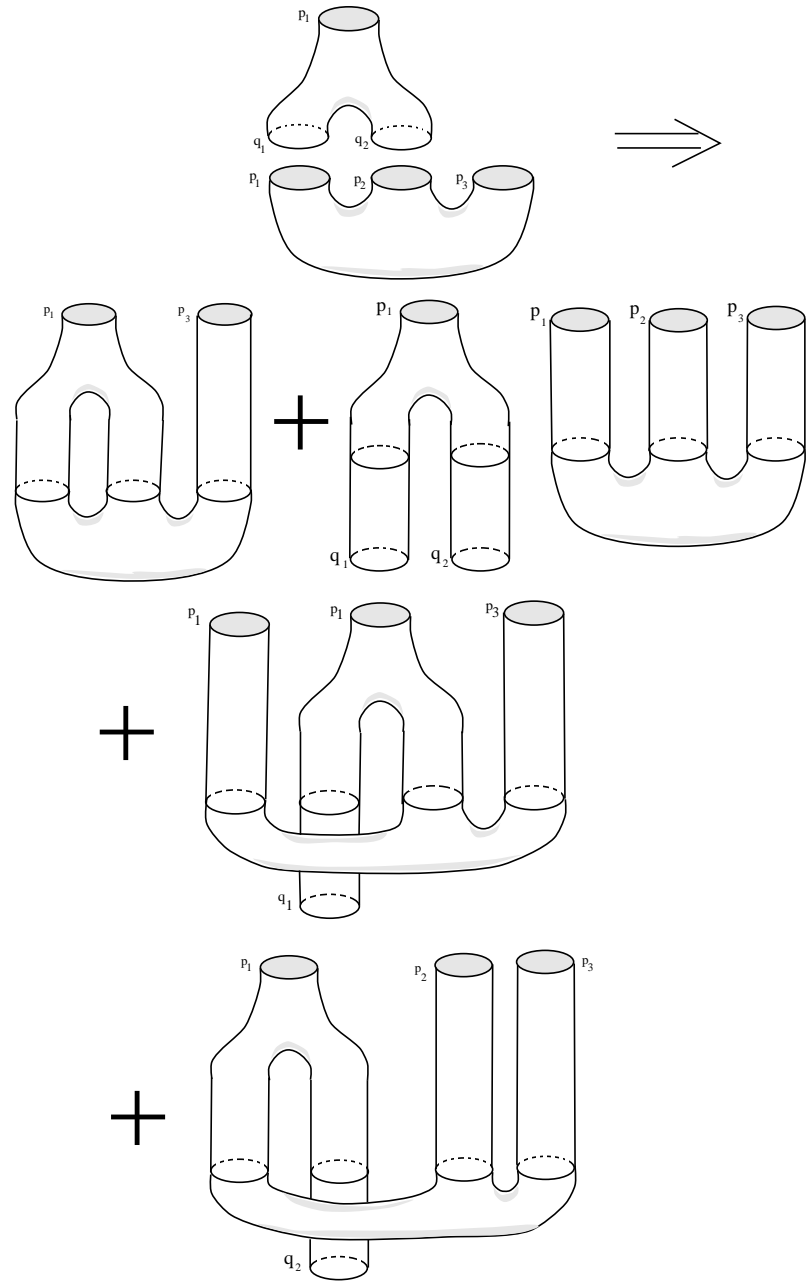,height=170mm}} \caption{\small
         There are  four different way to glue the lower and upper surfaces
    on the picture along their matching ends, i.e. the ends denoted by $p$'s
    and $q$'s with the same index.  These 4 ways correspond to  4
    terms in the composition formula for differential operators:
                             $(\hbar^{-1}p_1 p_2 p_3) \circ(\hbar^{-1}q_1 q_2 p_1)=p_1 p_3+
                                \hbar^{-2} q_1 q_2 p_1^2 p_2 p_3 +
                \hbar^{-1}q_1 p_1^2 p_3+\hbar^{-1}q_2 p_1 p_2 p_3 \;.
               $  We  are ignoring here the sign issues and assuming all  the boundary components
                 to be simple orbits.
                                 }
        \label{fig:combinatorics}
\end{figure}

One can go further and include into the picture  moduli spaces of punctured holomorphic
curves
of higher genus.  Introducing a new variable, denoted $\hbar$, to keep track
of the genus, one can associate with each $0$-dimensional moduli
space
 $\Mc^d_g(\Gamma^-,\Gamma^+)/\R$
   a monomial  $$\frac{n_{\Gamma^-,\Gamma^+,d,g}}{s^-!s^+!(i_1^-)!\dots(i_{s^-}^-)!(i_1^+)!\dots(i_{s^+}^+)!}
    q_{\g^-}^{I^-}p_{\g^+}^{I^+}\hbar^{g-1}z^d,$$
  and form a generating function $\bH=\hbar^{-1}\sum\limits_{g=0}^\infty\bH_g\hbar^{g}$
  counting all rigid holomorphic curves of arbitrary
  genus, whose term $\bH_0$ coincides with $\bh$.    Again,
  the codimension $1$ strata of the boundary of the moduli spaces
  $\Mc_g(\Gamma^-,\Gamma^+)$ consists of height 2 stable curves, but unlike the case of  rational
  curves, two connected components
  on different levels can be glued together along an arbitrary number of ends.
   The combinatorics of such gluing
  can be described by the formalism of algebra of higher order differential operators.
  Fig. \ref{fig:combinatorics}
  illustrates how the composition
of
  differential operators can be interpreted via gluing of Riemann
  surfaces with punctures. A letter $p_i$ in the picture represents a differential
  operator
             $\hbar\frac{\partial}{\partial q_i}$, and  a surface of genus $g$ with upper punctures
    $p_{i_1},\dots,p_{i_k}$ and lower punctures $q_{j_1},\dots,q_{j_l}$ represents
    a differential operator
        $$\hbar^{g-1}q_{j_1}\dots q_{j_l}p_{i_1}\dots p_{i_k}=\hbar^{g-1}q_{j_1}\dots q_{j_l}
        \big(\hbar\frac{\partial}{\partial q_{i_1}}\big)
        \dots \big(\hbar\frac{\partial}{\partial q_{i_k}}\big).$$
                 Thus
         we are led to consider $\bH$ as an element of the Weyl super-algebra
  $\fW$. This algebra should be
   viewed
  as a quantization of the Poisson algebra $\fP$, so that the description of the boundary of the
  moduli spaces is given by the equation $[\bH,\bH]=0$, where   $[\;,\;]$ denotes the commutator in $\fW$.
  As in  the rational case, this  identity is equivalent
  to the identity $D^2_{\bH}=0$  for the differential
  $D^\bH(f)=[H,f]$.
  Hence we can define the
  homology algebra $H_*(\fW,D^{\bH})$, which   also turns out to be an invariant of the contact manifold
  $(V,\xi)$. Similarly to the standard Gromov-Witten theory for closed
  symplectic manifolds one can develop an even more general
  formalism  by encoding in $\bH$ information about higher-dimensional moduli spaces
  of holomorphic curves. This leads to a deformation of the differential
  algebra $(\fW,D^{\bH})$ along the space of closed
  forms on $V$. The corresponding family of homology algebras   is
  then  parameterized by $H^*(V)$.


  After going that far it is natural to  make the above algebraic structure
  for contact manifolds a part
  of a formalism in the spirit of topological field theory, which we call
   {\it Symplectic Field Theory},
  and which also includes the theory of Gromov-Witten invariants of closed manifolds.
  To do that one considers moduli spaces of holomorphic curves with cylindrical ends
  in directed symplectic cobordisms $W=\ora{V^-V^+}$. The generating function
  counting rational holomorphic curves
  in $W$ can be naturally written as a function  $\bff(q^-,p^+ )$ of $p^+$-variables
  associated with the positive end $V^+$, and
  $q^-$-variables associated with the negative end $V^-$ of the cobordism $W$.
  It turns out that the Lagrangian submanifold
  in $(-\mathbf{V}^-)\times\mathbf{V}^+$ generated by the function $\bff$
  defines a {\it Lagrangian correspondence}
$  \mathbf{L}_W\subset(-\mathbf{V}^-)\times(\mathbf{V}^+) $
  which transforms the Hamiltonian
  functions $\bh^+$ and  $\bh^-$ to each other, i.e.
  \begin{equation}\label{eq:HJ-informal}
  \left(\bh^-(p^-,q^-)-\bh^+(p^+,q^+)\right)|_{\mathbf{L}_W}=0,
  \end{equation}
  where
  \begin{equation*}
  L_W=
  \begin{cases}
  q^+_{\g^+}=&\kappa_{\g^+}\frac{\partial \bff}{\partial p^+_{\g^+}}(q^-,p^+));\\
p^-_{\g^-}=&\kappa_{\g^-}\frac{\partial \bff}{\partial
q^-_{\g^-}}(q^-,p^+)).\\
\end{cases}
\end{equation*}
   We recall that $\kappa_{\g^\pm}$ denotes the multiplicity of the orbit $\g^\pm$.

  The composition of symplectic cobordisms produces
  the composition of Lagrangian correspondences,
  so that  if one  consider a ``Heegard splitting" of  a closed symplectic
  manifold $W$ along a contact hypersurface $V$, then  the computation of Gromov-Witten invariants
   of $W$ can be viewed
  as a Lagrangian intersection problem in the   symplectic super-space
  $\mathbf{V}$ associated to the contact manifold $V$.

  After what was said it should not come as a surprise that in
   the quantized picture Lagrangian correspondences
  are being replaced by Fourier integral operators, and the
   composition of correspondences by the convolution
  of the corresponding operators.

 We describe below the SFT-formalism with more details. We treat
contact manifolds in Section \ref{sec:af-contact} and  symplectic
cobordisms in Section \ref{sec:af-cobordism}. Section
\ref{sec:chain}  is devoted to  the   SFT-version of the  chain  homotopy statement
 in Floer homology theory.   In Section  \ref{sec:composition}  we introduce the composition
 formula for the SFT-invariants of symplectic cobordisms. In Section \ref{sec:contact-invariants}
 we discuss    how the   introduced algebraic structures of contact manifolds
 depend on extra choices.  Section \ref{sec:master}  is devoted
 to a differential equation for the  potential $\bF$ of a directed symplectic
 cobordism    with  a {\it non-empty } boundary.  Together with the gluing formula from Section
 \ref{sec:composition} this  equation  provides an effective tool for  computing
  Gromov-Witten invariants.
The remainder
of the paper has even  more sketchy character than the rest of
the paper.
 Section \ref{sec:Legendrian}   is
devoted to invariants of Legendrian submanifolds via SFT. Section
\ref{sec:remarks} is devoted to various examples and possible
generalizations of  SFT. In particular, in Section
\ref{sec:Bott} we discuss   how  one can adapt the theory to
include an important for applications, though non-generic, case of
contact forms with continuous families of periodic orbits. In
Section \ref{sec:computing} we describe a  new recursive
procedure for computing rational  Gromov-Witten invariants of $\C
P^n$. Finally, in Section \ref{sec:satellites} we  just touch the
wealth of other  invariants which exist in  Symplectic Field
Theory.

 \subsection{Contact manifolds}\label{sec:af-contact}

 \subsubsection{Evaluation maps}\label{sec:contact-eval}
 Let $(V,\xi)$ be a contact manifold with a fixed contact form
 $\alpha$, $(W=V\times\R,d(e^t\alpha))$ the symplectization of
 $(V,\xi)$, and $J$  a compatible almost complex structure.  As in Section \ref{sec:holomorphic}
  we denote  by $f_V$ and $f_\R$
  the $V $- and $\R$-components of  a
 $J$-holomorphic curve $f$ in $W$, and
 by $\Mc_{g,r,s^-,s^+}(W,J)$ the  disjoint union
   $$\bigcup\Mc_{g,r}^A(\Gamma^-, \Gamma^+),$$
   which is taken over all $A\in H_2(V)$, and all
    sets $ \Gamma^-,\Gamma^+\subset\Pc_\alpha$ of cardinalities $s^\pm$.

Let us view the set  $\Pc=\Pc_\a$ of periodic orbits of  the Reeb fields $ R_{\alpha }$ as a discrete topological space.
 It naturally splits into the   disjoint
union
$$\coprod\limits_{k=1}^{\infty}\Pc _k,$$ of identical subspaces, where $\Pc_k$ is the space of
  periodic orbits of
multiplicity $k$.

           Consider now three sets of evaluation maps:
            $$ev^0_i: \Mc_{g,r,s^-,s^+}/\R \to V,\;i=1,\dots,r,$$
            $$ev^+_j:\Mc_{g,r,s^-,s^+ }/\R    \to \Pc ,\;j=1,\dots      ,s^+,$$
            and
              $$ev^-_k:\Mc_{g,r,s^-,s^+ }/\R    \to \Pc ,\;k=1,\dots      ,s^-,$$
              where $ev_i^0$ is the  evaluation map $f_V(y_i)$
               at the $i$-th marked point $y_i$, while
                            $ev^\pm_j$   are the evaluation maps
                            at asymptotic markers $\mu^{\bx^\pm}_j$.
               More precise,
               let
                           $$\overline{f}=(f,j,\bx^-,\bx^+,\by,\mu^{\bx^-},\mu^{\bx^+})\in
                            \Mc_{g,r,s^-,s^+ },$$
              and $f$ be asymptotically cylindrical over a
                            periodic orbit $\gamma^\pm_j\in \Pc $ at $\pm\infty$ at the
              puncture   $x^\pm_j$. Then $ev_j^\pm(\overline{f})$ is a point of $
              \Pc $ representing the orbit $\g^\pm_j$ (comp. Section \ref{sec:Bott} below).

              All the above evaluation maps can be combined into a map
              $$ev:\Mc_{g,r,s^-,s^+ }/\R\to V^r\times (\Pc^-)^{s^-}\times (\Pc^+)^{s^+}\,,$$
which extends to the compactified moduli space
$\overline{\Mc_{g,r,s^-,s^+ }/\R}$.
\subsubsection{Correlators}\label{sec:correlators}
          Now we are ready to define
  correlators.
 Given $r$    differential
forms $\theta_1,\dots,\theta_r$ on $V$ and $s^\pm$ ($0$-dimensional) cohomology classes $\alpha^\pm_1,\dots,
\alpha^\pm_{s^\pm}\in H^*(\Pc )$  we define the    degree $-1$, or contact {\it
correlator}
\begin{eqnarray}\up{-1}\langle\,\theta_1,\dots,\theta_r;\alpha^-_1,\dots,
\alpha^-_{s^-};\alpha^+_1,\dots,
\alpha^+_{s^+} \rangle^A_{g}=\\
 \nonumber \int\limits_{\overline{ {\Mc_{g,r,s^-,s^+ }^A}/\R}}ev^*(\theta_1\otimes\dots
\otimes \theta_r\otimes\alpha^-_1\otimes\dots\otimes\alpha^-_{s^-}
\otimes\alpha^+_1\otimes\dots\otimes\alpha^+_{s^+} ).
\end{eqnarray}

Usually we will assume that the forms $\theta_1,\dots,\theta_r$
are closed, but even in this case the above  correlator  depends
on the actual forms, and not just their cohomology classes,
because the domain of integration may have a boundary.  As we will see
below the
superscript $-1$ in $\up{-1}\langle\dots\rangle $ corresponds
to the grading of the generating function for these
correlators. It also refers to the
enumerative meaning of the correlators: they count
components of $1$-dimensional moduli spaces of holomorphic curves.
We will consider below also correlators $\up0\langle\dots\rangle $
and $\up1\langle\dots\rangle$, counting $0$-dimensional and
$-1$-dimensional (i.e. appearing in 1-dimensional families) moduli
spaces of holomorphic curves.

 If we are given $K$ linearly independent differential forms
  $\Theta_1,\dots,\Theta_K$,
  then it is convenient
  to introduce a ``general form" $t=\sum\limits_1^Kt_i\Theta_i$
  from the space $L=L(\Theta_1,\dots,\Theta_K)$
  generated by the  chosen forms, and view $t_i$ as graded variables
  with  $\deg(t_i) =\deg(\Theta_i)-2$. In particular, all terms in
  the sum $\sum\limits_1^Kt_i\Theta_i$ have even degrees.

  Let us   consider two copies
 $\Pc^+$ and $\Pc^-$ of the $0$-dimensional space
 $\Pc$, one   associated with the positive end of $W$, the other with
  the negative one. Cohomology classes  in $\Pc_+$ we will denote
  by $p$, and in $\Pc_-$ by $q$, and will write
  $$p=\sum_{\g\in\Pc}\frac{1}{\kappa_\g}p_\g[\g],\quad q=
  \sum_{\g\in\Pc}\frac{1}{\kappa_\g}q_\g[\g],$$
   where   $\kappa_\g$ is
 the multiplicity of $\g$, and the cohomology classes $[\g]$
  form the canonical basis of $H^*(\Pc)$, dual to points in $\Pc$.
  Of course, speaking about cohomology classes of a discrete space may sound
  somewhat ridiculous. However, this point of view   is useful, especially
   in preparation for a more general  case when
    some periodic orbits may be degenerate and  thus the spaces  $\Pc^\pm$
      need not to be  anymore discrete, see Section \ref{sec:Bott} below.
 We will also fix a basis $A_1,\dots,A_N$ of
 $H_2(V)$. The coordinate vector  $d=(d_1,\dots,d_N)$ of a class $A$
 is called the degree.   Here $d_j$ are integers, while we   consider
 $t,p,q $ as graded variables,
  where the  degrees of the variables $p,q$ are   defined by the  formulas
  $$\deg (p_\gamma )  =-\CZ(\gamma)+(n-3),\, $$
 $$ \deg (q_\gamma )  =+\CZ(\gamma)+(n-3).$$
 The choice
  of grading,
   somewhat strange at the first glance,  is explained by Proposition \ref{thm:finiteness-of-H} below.

     The correlators
  $$ \langle \underbrace{t,\dots,t}_r;\underbrace{q,\dots,q  }_{s^-}
     ;\underbrace{p,\dots,p}_{s^+}\rangle^d_g$$ with  different
     $r,d,g$ determine all the correlators involving forms from the
     space $L$.

\subsubsection{Three differential algebras}\label{sec:3algebras}

   Similar to the theory of Gromov-Witten invariants of   closed symplectic manifolds we
   will organize all correlators
      into a generating function, called {\it Hamiltonian},
   $${\bH}=\frac{1}{\hbar}\sum\limits_{g=0}^\infty{\bH}_g\hbar^g,$$
   where,
\begin{equation}
 {\bH}_g=\sum\limits_d\sum\limits_{r,s^\pm=0}^\infty\frac{1}{r!s^-!s^+!}
     \up{-1}\langle \underbrace{t,\dots,t}_r;\underbrace{q,\dots,q  }_{s^-}
     ;\underbrace{p,\dots,p}_{s^+}\rangle^d_g z^d,
     \end{equation}
and  $t=\sum_1^Kt_i\Theta_i$. We will assume throughout the paper,
that
all forms $\Theta_1,\dots,\Theta_K$ are closed (see,
however, Remarks \ref{rem:non-closed1} and \ref{rem:non-closed2}, and
Section \ref{sec:master} below).
 The variables $\hbar$ and $z=(z_1,\dots,z_N)$ are also considered as graded
 with $\deg \hbar= 2(n-3)$ and $\deg (z_i)=- 2c_1(A_i)$, where $c_1$ is the first
 Chern class of  the almost complex structure $J$.

\begin{proposition} \label{thm:finiteness-of-H}
\begin{description}
  \item{a)}  For each $g=0,\dots$  the series $\bH_g$
   can be viewed as  formal power series in variables $p_\g$   with coefficients
  which are   polynomials   of variables $q_\g$ and formal power series of $t_i$ \footnote{
  In fact, $\bH_g$ depends polynomially on all variables $t_i$ of degree $\neq 0$.
  The degree $0$ variables, i.e. the variables
  associated with $2$-forms, enter into the constant part of $\bH_g$ (i.e.
  the part describing constant holomorphic curves) polynomially, while
   the non-constant part    of $\bH_g$ depends
   polynomially on $e^{t_i}$. This fact is similar to the standard Gromov-Witten theory
  and  will not discussed  in the present paper.}
  with coefficients
  in the group algebra $\C[H_2(V)]$
    {\rm (which we identify with the algebra
  of Laurent polynomials of $z$ with complex coefficients)};
  \item{b)} All terms of $\bH$  have degree $-1$;
  \item{c)} $\bH\big|_{p=0}=\bH_{\const}$, where
   $$\bH_{\const}=\hbar^{-1}\sum\limits_{g,r=0}^\infty
   \frac{1}{r!}\up{-1}
   \langle\underbrace{t,\dots,t}_r\rangle^0_g\hbar^g$$
   accounts for the contribution of constant holomorphic curves.
   In particular, $\bH\big|_{p=0}$ is independent of $q$ and $z$.
  \end{description}
  \end{proposition}
The polynomial dependence of $\bH_g$  on variables $q_\g$ and $z$   in a  geometric language just means that
the union $\overline{\Mc^d_g(\Gamma^+)}$, $\Gamma^+=\g_1\dots\g_{s^+}$,
of the  compactified moduli spaces of holomorphic curves  of a fixed genus of any degree
 with prescribed
positive ends  $\g^+_1,\dots,  \g_{s^+}$
is compact, and in particular that there are only finitely many possibilities for the degrees
and the negative
ends of these curves.
 This follows from the fact  that  for each curve $C\in \Mc_g(\Gamma^-,\Gamma^+)$ we have
 \begin{equation}\label{eq:pos-ends}
 0\leq \int\limits_C d\alpha= \sum\limits_{\g_i\in\Gamma^+}\int\limits_{\g_i}\alpha-
 \sum\limits_{\g_j\in\Gamma^-}\int\limits_{\g_j}\alpha\leq
  \sum\limits_{\g_i\in\Gamma^+}\int\limits_{\g_i}\alpha,
  \end{equation}
  the fact   that there exists a constant $m>0$ such that $\int\limits_{\g}\a>m$ for
  any periodic   orbit $\g\in\Pc_\a$ and Theorem \ref{thm:comp1}  above.
     Proposition \ref{thm:finiteness-of-H}b)   follows from the formula
  (\ref{eq:dim}) for the dimension of the moduli spaces of holomorphic curves, our degree
  convention   and the fact that a correlator
      $\up{-1}\langle\,\theta_1,\dots,\theta_r;\g^-_1,\dots,
\g^-_{s^-};\g^+_1,\dots,
\g^+_{s^+} \rangle^A_{g}$  may be different from $0$ only if the total dimension of the forms
$\theta_1,\dots\theta_r$ equals the dimension of the moduli  space   $\Mc^A_{g,r}( \g^-_1,\dots,
\g^-_{s^-};\g^+_1,\dots,
\g^+_{s^+})/\R $.
  Proposition \ref{thm:finiteness-of-H}c) just means that every non-constant holomorphic curves should have at least
  one positive end, which follows from   inequality  (\ref{eq:pos-ends}), or alternatively the
  maximum principle for holomorphic curves.
  \medskip

 Let us consider the {\it Weyl super-algebra}
  $\fW=\{\sum\limits_{\Gamma,g} f_{\Gamma,g }(q,t )p^\Gamma\hbar^g \}$,
  where
 $$\Gamma=(\g_1,\dots,\g_a),\gamma_i\in\Pc,\,a=1,\dots,\, p^{\Gamma}=p_{\g_1}\dots p_{\g_a},$$
  and        $f_{\Gamma,g }(q,t)$ are polynomials of variables
  $q_\gamma$ and formal power series of $t_i$\footnote{See the previous footnote.}. Proposition
  \ref{thm:finiteness-of-H}a) states that $\bH\in\hbar^{-1}\fW$.

  The product operation $F\circ G$ in $\fW$ is associative and satisfies the following commutation
  relations: all variables
  are super-commute (i.e. commute or anti-commute according to their grading), except
  $p_\gamma$ and $q_\gamma$ which correspond to the same periodic orbit $\gamma$. For these pairs of variables
   we have the following commutation relation:
 \begin{equation}
 [p_\gamma,q_\gamma]=p_\g\circ q_\g-(-1)^{\deg p_\g\deg q_\g}q_\g\circ p_\g=
 {\kappa_\gamma}\hbar \,
 \end{equation}
  where $\kappa_\g$ is the multiplicity of the orbit $\g$.
  The algebra $\fW$ can be represented as an  algebra of formal differential   operators  with respect
  to
  $q$-variables acting  {\it on the left} on the space of polynomials
    $  f(q,z,\hbar) $, by setting
    $$p_\g= {\kappa_\g}\hbar\ora{\frac{\partial}{\partial
    q_\g}}.$$
    Alternatively by setting
    $$q_\g={\kappa_\g}\hbar\ola{\frac{\partial}{\partial p_\g}}$$ we can
    represent
    $\fW$   as an algebra  of  polynomial differential operators  acting {\it on the right} on the algebra $\{\sum
    \limits_{\Gamma,g}f_{\Gamma,g}(q,z)\hbar^gp^\Gamma\}$ of formal power series
    of $\hbar$ and the $p$-variables.

  Notice that  the commutator $[F,G]$ of  two homogeneous elements $F,G\in\fW$
  equals $F\circ G-(-1)^{\deg F\deg G}G\circ F$, and hence if $F\in\fW$ is an odd element (i.e. all its
  summands are odd) then
     $[F,F] =2F\circ F,$
     and $[F,F]=0$ if $F$ is even. For any two elements
     $F,G\in\fW$ the commutator $[F,G]$ belongs to the ideal
     $\hbar\fW$.
                    According to Proposition \ref{thm:finiteness-of-H}   the Hamiltonian
     $\bH$ can be viewed as an element of $\frac{1}{\hbar}\fW$, and hence
     the above remark shows that for $F\in\fW$
     we have $[\bH,F]\in\fW$.
  \begin{theorem} \label{thm:HH}
 The Hamiltonian  $\bH$ satisfies the identity
  \begin{equation}\label{eq:HH}
 \bH\circ \bH=0.
 \end{equation}
 \end{theorem}

This theorem (as well as Theorems
\ref{thm:SFT-cobordism}, \ref{thm:SFT-chain} and \ref{thm:SFT-composition}
below) follows from the description of the boundary of the
corresponding moduli spaces of holomorphic curves. As it was stated
in Proposition \ref{prop:boundary}
  this boundary
  is tiled by codimension one
 strata represented by stable curves of height $2$, so that the
(virtual) fundamental cycles of the boundary of the compactified  moduli
spaces
  can be
symbolically written as   $\partial [\Mc]=\kappa\sum [\Mc_{-}]\times [\Mc_{+}]$,
where $[\Mc_\pm]$ are chains represented by the corresponding moduli spaces and  where
 the coefficient
 $\kappa$ depends on multiplicities of orbits along which the two levels of the corresponding
 stable curve are glued. Together
with the Stokes formula $\int\limits_{[\Mc]} d\omega = \int\limits_{\partial [\Mc]} \omega
$, and the fact that the integrand is a closed form, we obtain
identity
(\ref{eq:HH}).

\begin{remark}\label{rem:non-closed1}
{\rm The same argument shows that when the forms $\Theta_i$
generating the space $L$  are not necessarily closed
we get the following  equation
\begin{equation}\label{eq:non-closed1}
d\bH+\frac12[\bH,\bH]=0,
 \end{equation}
 which
generalizes (\ref{eq:HH}) and can be interpreted as the
zero-curvature equation for the connection $d +[\bH,\cdot\,]$. We
denote here by $d$ the de Rham differential, i.e.
 \begin{equation}
 \begin{split}
 d\bH&=d\Big(\sum\limits_{d,g}\sum\limits_{r,s^\pm=0}^\infty\frac{1}{r!s^-!s^+!}\\
   & \quad \big\langle \underbrace{\sum\limits_1^Kt_i\Theta_i,\dots,\sum\limits_1^Kt_i\Theta_i}_r;\underbrace{q,\dots,q  }_{s^-}
     ;\underbrace{p,\dots,p}_{s^+}\big\rangle^d_g
     z^d\hbar^{g-1}\Big)\\
      \end{split}
     \end{equation}
 \begin{equation*}
 \begin{split}
     \;\;&=\sum\limits_{d,g}\mathop{\sum\limits_{s^\pm=0,}}\limits_{r=1}^\infty\frac{1}{(r-1)!s^-!s^+!}\\
    &\quad \big\langle\sum\limits_1^K t_i
     d\Theta_i, \underbrace{\sum\limits_1^Kt_i\Theta_i,
     \dots,\sum\limits_1^Kt_i\Theta_i}_{r-1};\underbrace{q,\dots,q  }_{s^-}
     ;\underbrace{p,\dots,p}_{s^+}\big\rangle^d_g
     z^d\hbar^{g-1}\,.\\
     \end{split}
     \end{equation*}
 }
\end{remark}
\bigskip

 The identity $\bH\circ\bH=0$
 is equivalent to $[\bH,\bH]=0$, because $\bH$ is an odd element.
 Let  us  define the differential
 $D=D^{\bH}:\fW\to\fW$
 by the formula
 \begin{equation}
 D^{\bH}(f)=[\bH,f]\quad\hbox{ for }\quad f\in\fW .
 \end{equation}
Then Theorem \ref{thm:HH} translates into the identity $D ^2=0$. The differential
 $D^{\bH}$ satisfies the Leibnitz rule, and thus
 $(\fW,D )$
 is a        differential Weyl (super-)algebra.
In particular, one can define the homology algebra $H_*(\fW,D)$,
which inherits  its multiplication operation from  the Weyl
algebra $\fW$.   The differential $D^{\bH}$ extends in an obvious way to the modules
$\hbar^{-k}\fW$, $k=1,\dots.$

\begin{example}\label{ex:circle}
{\rm Let
 $V=S^1$. We have  in this case
 \begin{equation}\label{eq:circle}
 \bH=\hbar^{-1}\big(\frac{t_1t_0^2}{2}+ t_1\sum p_kq_k-\frac{t_1\hbar }{24}\big),
 \end{equation}
  where
$ t=t_01+t_1d\phi$ is a general harmonic differential form on $S^1$,
so that $\deg t_1=-1$,   $\deg t_0=-2$, $\deg\hbar=-4$ and  $\deg p_k,\deg q_k=-2$,
which corresponds to the convention
that the Maslov index of any path in the
$1$-point group $Sp(0)$ equals 0.
The term $t_1t_0^2/2=\int\limits_{S^1} t^{\wedge 3}/6$ is the contribution of the moduli
space $S^1$ of constant maps $\C P^1\to \R\times S^1$ with $3$ marked points.
 The term $-\frac{t_1\hbar}{24}$ is accounted for
the contribution of constant curves of genus $1$ (see
\cite{Witten2}), and the term $t_1p_kq_k$ represents the
contribution $t_1=\int\limits_{S^1} t$ of trivial curves of multiplicity
$k$ with one marking. All other curves do not contribute to $\bH$
for dimensional reasons and because $t_1^2=0$.

Notice that if we organize all variables $p_k,q_k$   into
 formal Fourier series  (comp. \cite{Getzler})
\begin{equation}
u(x)= \sum\limits_{k=1}^\infty(p_ke^{ixk}+q_ke^{-ixk}),
\end{equation}
 then the term accounting for
  the contribution of rational curves in the formula (\ref{eq:circle})
  takes the form
  \begin{equation}
     \frac{t_1}{4\pi }\int\limits_0^{2\pi} (t_0+   u(x)) ^2  \ dx,
\end{equation}
see further discussion of this $u$-formalism in Section \ref{sec:Bott} below.
 }
\end{example}

 We will  associate now with $(\fW,D)$    two other differential algebras, $(\fP,d)$ and $(\fA,\partial)$, which can be viewed
 as {\it semi-classical} and {\it classical} approximations of the Weyl  differential algebra.

Let us denote  by
\begin{description}
 \item{-} $\fP$ -- a graded Poisson algebra  of formal power series in variables
  $p_\g $ with coefficients which are polynomials of $q_\g, z_j,
  z_j^{-1}$, and formal power series of $t_i$,
  \footnote{See the first footnote in Section \ref{sec:3algebras}.} and by
\item{-} $\fA $ --  a graded commutative algebra generated by variables
 $q=\{q_\g\}_{ \g\in\Pc}$   with coefficients in   the algebra
    $\C[H_2(V)][[t]]$.
    \end{description}
  The Poisson bracket on $\fP$ is defined by the formula
  \begin{equation}\label{eq:Poisson}
  \{h,g\}= \sum\limits_\g  {\kappa_\g}\left(\frac{\partial h}{\partial p_\g}
  \frac{\partial g}{\partial q_\g}-(-1)^{{\deg h}{\deg g}}
   \frac{\partial g}{\partial p_\g}\frac{\partial h}{\partial q_\g}
 \right)    ,
 \end{equation}
 assuming that $h$ and $g$ are {\it monomials}.
 When computing partial derivatives, like $\frac{\partial h}{\partial q_\g}$,
  one should remember that we are working in the super-commutative environment,
 and in particular the operator $\frac{\partial }{\partial q_\g}$
 has the same parity as  the variable $q_\g$.
 \begin{remark}\label{rem:u-notation}
 {\rm
  Notice, that if  similarly to  Example \ref{ex:circle} above
 we organize the variables $p_k=p_{\g_k},q_k=q_{\g_k}$ corresponding to multiples
 of each simple periodic orbit $\g=\g_1$ into  a
 Fourier series
 $$u_\g=\sum\limits_{k=1}^\infty(p_ke^{ixk}+q_ke^{-ixk}),$$
 then the value of the Poisson tensor
 (\ref{eq:Poisson}) on covectors $\delta u,\delta v$ takes      the
  form
 \begin{equation}
 \label{eq:u-notation}
 \frac{1}{2\pi i}\int\limits_0^{2\pi}  (\d u)'\d vdx  \,.
 \end{equation}
    }
 \end{remark}
 \bigskip

    In order to define differentials
     on the algebras $\fA$ and $\fP$ let us  first     make
the following observation.
\begin{lemma} \label{lm:HH-Weyl}
We have $$[\bH,\bH]=\frac{1}{\hbar}\{\bH_0,\bH_0\}+\dots\,,$$
and for any $f=\frac{1}{\hbar}\sum f_g\hbar^g\in \frac{1}{\hbar}\fW$
 we have
$$D^\bH(f) = \frac{1}{\hbar}\{\bH_0,f_0\}+\dots\,.$$ In
particular, $\bH_0$ satisfies the equation
$\{\bH_0,\bH_0\}=0.$
\end{lemma}

To cope with a growing number of indices we will rename $\bH_0$
into $\bh$.
   Lemma \ref{lm:HH-Weyl} allows us to define the differential
   $d=d^{\bh}:\fP\to\fP $ by the formula
   \begin{equation}
   dg=\{\bh,g\} \quad\hbox{ for }\quad g\in\fP .
\end{equation}
Theorem \ref{thm:HH} then implies
   \begin{proposition}\label{prop:dd-Poisson} We have
   $d^2=0$ and $d\{f,g\}=\{df,g\}+ (-1)^{\deg f}\{f,dg\}$ for any homogeneous element
   $f\in\fP$. In other words,   $(\fP,d)$ is a
   graded differential Poisson algebra with  unit.
   \end{proposition}
   Proposition \ref{prop:dd-Poisson}
   enables us to define the homology $H_*(\fP,d)$ which inherits  from
    $\fP$ the structure of a graded Poisson algebra with unit.

     Let us recall that according to \ref{thm:finiteness-of-H}
     $\bh|_{p=0}=\bh_{\const}$,
     where $\bh_{\const}$ accounts for constant rational
     holomorphic curves, and thus it is independent of $q$-variables.
     In fact, $$\bh_{\const}(t)=\frac16\sum\limits_{i,j,k=1}^K c^{ijk}t_it_jt_k,$$
     where $c^{ijk}=\int\limits_V\Theta_i\wedge\Theta_j\wedge\Theta_k$ are
     the structural
     coefficients of the cup-product.
          Hence,  $$\bh=\bh_{\const}+\sum h_{\g }(q,t,z)p_\g   +\dots,$$
     where $\dots$ denote terms at least quadratic in $p_\g$.  Thus
     we have
      \begin{equation}
     \{\bh,\bh\} = 2\sum\limits_{\g ,\g'\in\Pc }\kappa_{\g'}
      h_{\g'  }(q,t)\frac {\partial h_{\g  }}{\partial q_{\g'}}(q,t)p_{\g}+ o(p)=0\,.
      \end{equation}
        Therefore,
        \begin{equation}\label{eq:d2}
         \sum\limits_{\g'  \in\Pc }\kappa_{\g'}
      h_{\g' }(q,t)\frac{\partial h_{\g  }}{\partial
      q_{\g'}}(q,t)=0
      \end{equation}
      for all $t$ and all $\g\in\Pc$.

   Let us define a differential $\partial: \fA\to\fA$ by the
      formula

      \begin{equation}
      \partial f=\{\bh,f\}|_{\{p=0\}}=\sum\limits_{\g\in\Pc}\kappa_\g h_\g\frac{\partial
      f}{\partial q_\g}.
      \end{equation}
      Then the equation (\ref{eq:d2}) is equivalent to

      \begin{proposition}\label{prop:zero-section}
     $\partial^2=0$, and hence $(\fA,\partial)$ is a graded
     commutative
     differential algebra with  unit.
     \end{proposition}

   The homology group $H_*(\fA,\partial)$
     inherits the structure
 of a graded commutative algebra with unit.

As it was already mentioned in Section \ref{sec:informal},
 it is convenient to view the Poisson algebra $\fP$ as an algebra
of functions on an infinite-dimensional symplectic super-space $\bV$ with the even symplectic form
$\bo=\sum k_\g^{-1}dp_\g\wedge dq_\g$.
Then  the differential $d^\bh$ is   the  Hamiltonian  vector field on $\bV$ generated by the  Hamiltonian
function $\bh$. One should remember, however, that the $p$-variables are formal,
so all that we have is the infinite jet of the symplectic space
$\bV$ along the $0$-section.  The equation $\bh|_{p=0}=\bh_{\const}$ translates into
 the fact that the vector field $d^\bh$
is tangent to the $0$-section, and the differential $\partial$ is just the restriction of this vector field
to the $0$-section.
  The higher order terms in the expansion of $\bh$ with respect
  to $p$-variables define a sequence  of (co-)homological operations on the algebra $\fA$.
      \medskip

    Notice  also that the differentials $D,d$ and $\partial  $ do not involve any differentiation
    with respect to $t$. Hence the  differential algebras
      $(\fW,D^{\bH})$, $(\fP,d^{\bh})$ and
       $(\fA,\partial)$ can be viewed
             as {\it families} of differential algebras, parameterized
      by $t\in H^*(V)$, and in particular,
      one can compute the homology at any fixed $t\in H^*(V)$.   We will
      sometimes denote
       the corresponding
       algebras and their homology
       groups with the subscript $t$, i.e.
        $(\fW,D)_t$,  $H_*(\fP,d)_t,$ etc., and call them
       {\it specialization} at the point $t\in H^*(V)$.
  We will also use the notation
  $$H_*^{\SFT}(V,\xi|\,J,\a ),\;\;
    H_*^{\RSFT}(V,\xi|\,J,\a  ), \quad\hbox{and }\quad
  H_*^{\cont}(V,\xi|\,J,\a )$$
     instead of
   $H_*(\fW,\partial), H_*(\fP,\partial)$ and $H_*(\fA,\partial)$,
   and
   will usually omit the   extra data $J,\a$ from the notation:
  as we will see in Section \ref{sec:composition} below all these
  homology algebras are independent  of  $J,\a$ and other extra
  choices, like
  closed forms representing cohomology classes of $V$, a coherent orientation of the moduli spaces,
  etc. The
  abbreviation $\RSFT$ stands here for  {\it Rational Symplectic
  Field Theory}.
\begin{remark}\label{rem:grading}
{\rm It is   important to observe that the algebras $\fW,\fP $ and
$\fA$ have an additional grading by elements of $H_1(V)$ (comp. Section
\ref{sec:Floer} above). This
grading is also inherited by the corresponding homology algebras. However, this grading carries a
non-trivial information only when we consider homology of algebras, specialized at points
$t=\sum\limits_1^K t_i\Theta_i$ with  $t_i=0$ for at least some of the
  coordinates $t_i$ corresponding
$1$-dimensional forms.  Otherwise
all cycles in these algebras are graded by  the $0$-class from $H_1(V)$.}
\end{remark}

\subsection{  Symplectic cobordisms}\label{sec:af-cobordism}

\subsubsection{Evaluation maps and correlators}
 \label{sec:sympl-eval}
Let us now repeat the constructions of the previous section for a
general directed symplectic cobordism $W=\ora{V^-V^+}$
 between  two contact manifolds
 $V^-$ and $V^+$ with fixed contact forms $\alpha^-$ and $\alpha^+$.
 As in Section \ref{sec:contact-eval} we consider the sets $\Pc^\pm$
  of periodic orbits of the Reeb fields $R^\pm=R_{\alpha^\pm}$ as   discrete topological spaces.

We denote by $\Mc^A_{g,r,s^-,s^+}( W,J)$ the  disjoint union
   $$\bigcup\Mc_{g,r}^A(\Gamma^-, \Gamma^+;W,J),$$
   where the union is taken over all sets $(\Gamma^-, \Gamma^+)$ of cardinality $(s^-,s^+) $, and
  consider  three sets of evaluation maps:
            $$ev^0_i: \Mc^A_{g,r,s^-,s^+}( W,J) \to W,\;i=1,\dots,r,$$
            $$ev^\pm_j:\Mc^A_{g,r,s^-,s^+}(W,J)    \to \Pc^\pm,\;j=1,\dots      ,s^\pm,$$
              where $ev_i^0$ is the  evaluation map $f (y_i)$ of
               the map $f$  at the $i$-th marked point $y_i$, while
                            $ev^\pm_j$   are the evaluation maps
                            at asymptotic markers
                            $\mu^{\bx^\pm}_j$, i.e.
  $ev_j^\pm(f)$ is a point of $
              \Pc^\pm$ representing the orbit $\g^\pm_j$, which contains the image of the
              corresponding marker.
               The evaluation maps $ev^0_i$ and $ev^\pm_j$ can be combined into a map
              $$ev:\Mc^A_{g,r,s^-,s^+}(W,J)\to W^r\times (\Pc^+)^{s^+}\times (\Pc^-)^{s^-}\,.$$

            Now we are ready to define   degree $0$, or symplectic
   correlators.  We will have to consider on $W$ differential
  forms, which do not necessarily have compact support, but
 have,   however, {\it cylindrical ends}. We say that a
   differential  form  $\theta$ on $W$ is said to have cylindrical ends
  if it satisfies the following
  condition:
\smallskip

\noindent there exists  $C>0$  such that
$$\theta|_{V^-\times(-\infty,-C)}=(\pi^-)^*(\theta^-)\quad\hbox{ and}\quad
 \theta|_{V^+\times(C,\infty)}=(\pi^+)^*(\theta^+),$$ where
$\theta^\pm$ are forms on $V^\pm$, and $\pi^\pm$ are the
projections of the corresponding ends to $V^{\pm}$.
 We will denote the forms $\theta^\pm$ also by  $\restr^\pm(\theta)$, or $\theta|_{V^\pm}$.
  In what follows we assume that all considered differential  forms on $W$ have   cylindrical ends.
 \smallskip

  Given $r$    differential
forms $\theta_1,\dots,\theta_r$ on $W$ and $s^\pm$   cohomology
classes $$\alpha^\pm_1,\dots, \alpha^\pm_{s^\pm}\in
H^*(\Pc^\pm)=H^*_0(\Pc^\pm)$$  we define the degree $0$
correlator
\begin{eqnarray}
\up{0}\langle\,\theta_1,\dots,\theta_r;\alpha^-_1,\dots,
\alpha^-_{s^-};\alpha^+_1,\dots,
\alpha^+_{s^+} \rangle^A_{g}= \hfill \\
\nonumber  \int\limits_{\overline{\Mc_{g,r,s^-,s^+ }^A}}ev^*(\theta_1\otimes\dots
\otimes \theta_r\otimes\alpha^-_1\otimes\dots\otimes\alpha^-_{s^-}
\otimes\alpha^+_1\otimes\dots\otimes\alpha^+_{s^+} ).
\end{eqnarray}

Similar to
  Section \ref{sec:contact-eval} above, we denote  the cohomology classes in
  $H^*(\Pc^+)= H^0(\Pc^+)$
  (resp. in   $H^*(\Pc^-) = H^0(\Pc^-)$) by
  $p^+$ (resp.   $q^-$), and write $$p^+=\sum\limits_{\g\in\Pc^+}k_g^{-1}p^+_\g[\g]\quad
  \hbox{\big(resp.}\quad q^-=\sum\limits_{\g\in\Pc^-}k_g^{-1}q^-_\g[\g]\,\big).$$
   We will also fix a basis $A_1,\dots,A_N$ of
 $H_2(W)$ and denote
 by $d=(d_1,\dots,d_N)$ the degree of $A$ in this basis.

Let us call a system of linearly independent
 closed forms  $\theta_1,\dots,\theta_m $
  on $W$ with cylindrical ends
 {\it basic}, if
           \begin{equation}\label{eq:basic}
 \begin{split}
 &\hbox{a) the image}\quad
    \restr^\pm(L(\theta_1,\dots,\theta_m ))\\
    &\quad\hbox{ generates }\quad \im (H^*(W)\to H^*(V^\pm)); \\
    & \hbox{b) the homomorphism}
 \quad \Ker\big((\restr^+\oplus \restr^-)|_L\big)\to
 H_{\mathrm{comp}}^*(W)\\
 &\quad\;\hbox {is bijective.}\\
 \end{split}
 \end{equation}
        Here we denote by $L(\theta_1,\dots,\theta_m )$ the subspace
 generated by the forms $\theta_1,\dots,\theta_m $, and by
 $H_{\mathrm{comp}}^*(W)$ the cohomology with compact support.
 Equivalently, one can say that a basic system of forms consists
 of a basis of $H^*(W)$ together with a basis of
 $\Ker\big(H^*_{\comp}(W)\to H^*(W)\big)$.

A general point     $t\in L(\theta_1,\dots,\theta_m )$  we will
   write in the form
    $t=\sum\limits_1^mt_i\theta_i.$
   The grading of the variables
 $t,p^+,q^- $  is defined   as in Section
 \ref{sec:correlators}:
 \begin{equation}
 \begin{split}
 \deg(t_i) =&\;\deg(\theta_i)-2;\,\\
  \deg (p_\gamma^+) =&-\CZ(\gamma^+)+(n-3),\\
   \deg (q_\gamma^-) =&\; \CZ(\gamma^-)+(n-3)\\
   \end{split}
   \end{equation}

\subsubsection{Potentials of symplectic cobordisms}\label{sec:relative}
  Let us now organize the  correlators
   into a generating function, called the { \it potential} of the
   symplectic
   cobordism $(W=\ora{V^-V^+},J,\a^\pm)$
  \begin{equation}\label{eq:GW-potential}     \begin{split}
  &{ { \bF}}=
   {  \bF}_{W,J,\alpha^\pm} =\frac1\hbar
   \sum\limits_{g=0}^\infty {\bF}_g\hbar^g,   \\
   &\quad\hbox{  where }\\
  &{\bF }_g=\sum\limits_d\sum\limits_{r  ,s^\pm=0}\frac{1}{r!s^+!s^-!}\up0\langle
    \underbrace{t,\dots,t}_r;
 \underbrace{q^-,\dots,q^-
 }_{s^-};\underbrace{p^+,\dots,p^+}_{s^+}
 \rangle^d_{g}z^d\,.\\
 \end{split}
 \end{equation}

When $W$ is a closed symplectic manifold, then the potential $F$
is just the Gromov-Witten invariant of the symplectic manifold
$W$. However, if $W$ is not closed, then $F$ itself is not an
invariant. It depends on    particular forms
$\theta_i$, rather than their cohomology classes, on $J$, on a
coherent orientation, and several other choices. We will see,
however, that the {\it homotopy class} of $F$, which we define in
Section \ref{sec:chain} below, is independent of most of these
choices.
\medskip

In order to make sense out of the expression for
$\bF$ let us consider a graded commutative
 algebra $ \fD=\fD({W,\a^\pm })$  which consists of
power series of the form
\begin{equation}
\sum\limits_{\Gamma,d,g}\varphi_{\Gamma,d,g} (q^-,t ) z^d(p^+)^{\Gamma}\hbar^g,
\end{equation}
where $\varphi_{\Gamma,d,g}$ are polynomials
 of $q_\g$, formal power series of variables
  $t_i$,\footnote{See the first footnote in Section \ref{sec:3algebras}.}
 and where $\Gamma$ and
$d$ satisfies the following Novikov type inequality:
\begin{equation}\label{eq:Novikov}
[\omega](d)=\sum
d_i\int\limits_{A_i}\omega>-\left|\Gamma\right|=-\sum\limits_{i=1}^k|\gamma_i|,
\end{equation}
where
$\Gamma=\{\g_1,\dots,\g_k\}$, and
$|\gamma_i|=\int\limits_{\gamma_i}\alpha^+$ is the period of the
periodic orbit $\g_i\in\Pc^+$, or equivalently its action. Recall
that $(p^+)^\Gamma=p^+_{\g_1}\dots p^+_{\g_k}$.

 \begin{proposition}
We have
$$\bF_{W,J,\alpha^\pm}
\in\frac1{\hbar}\fD({W,\alpha^\pm}).$$
\end{proposition}
Let us also consider a bigger algebra $\fD\fD$ which consists
of series
\begin{equation}
\sum\limits_{\Gamma,d}\varphi_{\Gamma,d} (q^-,t,\hbar ) z^d(p^+)^{\Gamma},
\end{equation}
where $\varphi_{\Gamma,d}$ are polynomials of $q^-_\g$, formal power series of $t_i$ and formal {\it Laurent
series}
of $\hbar$, while $\Gamma$ and $d$  still satisfy
 condition
(\ref{eq:Novikov}).
For instance, for any element $f\in \hbar^{-1}\fD$ we have $e^f\in\fD\fD$.

The algebra $\fD\fD=\fD\fD(W,J,\alpha^\pm)$ has a structure
of a {\it left $D$-module} over the Weyl algebra
$\fW^-=\fW(V^-,J,\alpha^-)$, and of a {\it right $D$-module} over the Weyl algebra
$\fW^+=\fW(V^+,J,\alpha^+)$.
Indeed,  we first associate with an element
$$\Delta^-= \sum\limits_{\Gamma=\{\g_1,\dots,\g_m\},\Gamma',d,g }
  {\delta}^-_{\Gamma',\Gamma,d,g  }( t )
   (q^-)^{\Gamma'}(p^-)^{\Gamma}z^d\hbar^g\in\fW^-$$
   a differential operator

\begin{equation}
\sum\limits_{\Gamma=\{\g_1,\dots,\g_m\},\Gamma',d,g }
 {\delta}^-_{\Gamma,\Gamma',d,g }(t )
   (q^-)^{\Gamma'}{\hbar^{m+g}}\prod\limits_{i=1}^m
   \kappa_{\g_i}
   \overrightarrow{\frac{\partial}
{\partial q^-_{\g_i}}}z^d,
\end{equation}
then  change   the coefficient ring via  the inclusion homomorphism
$H_2(V^-)\to
H_2(W)$,  and finally lift functions
$ {\delta}^-_{\Gamma,\Gamma',d,g }(t )$ to the space of forms with cylindrical
ends on $W$ via the restriction map $t\mapsto t|_{V^-}$.
 We will   denote  the resulting operator by
$\overrightarrow{\Delta^-}$. Similarly we associate with
$\Delta^+\in\fW^+$ an operator $ \overleftarrow{\Delta^+}$ by
first quantizing $q^+_\g\Rightarrow  \hbar
\kappa_\g\overleftarrow{\frac{\partial}{\partial p^+_\g}}$  and then making an
appropriate change of the coefficient ring.
                                   It is straightforward to verify that
for $f\in \fD$ we have $\overrightarrow{\Delta^-}f,\;f \overleftarrow{\Delta^+}\in \fD$,
and for $F\in\fD\fD$ we have
$\overrightarrow{\Delta^-}F,\; F\overleftarrow{\Delta^+}\in \fD\fD$.

\medskip

Let us denote the Hamiltonians (see Section \ref{sec:3algebras} above)  in
$\fW^\pm$ by $\bH^\pm$ and define a map
$D_W=D_W:\fD\fD\to\fD\fD$ by the formula

\begin{equation}\label{eq:DW}
D_W(G)= \overrightarrow{
\bH^-}G-(-1)^{\deg G}G\overleftarrow{\bH^+},\quad G\in\fD\fD,
\end{equation}
where we assume $G$ dimensionally homogeneous.
Clearly, Theorem \ref{thm:HH} implies that $D_W^2=0$. However, the differential algebra
$(\fD\fD,D_W)$ is too big
and instead of considering its homology we will define a differential on the algebra $\fD$,
or which is equivalent but more convenient, on the module $\hbar^{-1}\fD$.

For an {\it even} element $F\in\hbar^{-1}\fD$  let us define  a map
 $D^F=T_FD_W:\hbar^{-1}\fD\to\hbar^{-1}\fD$
by the formula
\begin{equation} \label{eq:DF}
D^F(g)=e^{-F}[D_W,g](e^F)=e^{-F}\big(D_W(ge^F)-(-1)^{\deg g}g D_W(e^F)\big),
\quad g\in\hbar^{-1}\fD.
\end{equation}
      The map    $T_FD_W$ is         the linearization of  the map
       $\wt D_W:\hbar^{-1}\fD\to\hbar^{-1}\fD$,
      defined by the formula $$\wt D_W(F)=e^{-F}D_W(e^F),\quad F\in\hbar^{-1}\fD\,.$$
            at the point $ F\in\hbar^{-1}\fD$.  Notice that   if $D_W(e^F)=0$ then
            $D^F(g)=e^{-F}D_W(ge^F)$.
      Let us first state a purely algebraic
      \begin{proposition} \label {prop:SFT-cobordism}
      Suppose that     for $F\in\hbar^{-1}\fD$ we have $D_W(e^F)=0$. Then
   \begin{description}
    \item{1.}   $(D^F)^2=0;$
\item{2.}  The homology algebra         $H_*(\fD,D^{F})$
inherits the structure of a left module over the homology
algebra
$ H_*(\fW^-,D^-)$, and
 the structure of a right module over the homology
algebra $ H_*(\fW^+,D^+)$;
  \item{3.} The homomorphisms
$F^\pm:  \fW^{\pm}  \to \fD $,
defined by  the formulas
\begin{equation} \label{eq:end-homomorphisms}
\begin{split}
&f\mapsto e^{-F}
\overrightarrow{f}e^{F},\,f\in\fW^-,\quad\hbox{and}\\
& f\mapsto  e^{F} \overleftarrow{f}e^{-F},\,f\in\fW^+\,,\\
\end{split}
\end{equation}
commute with the  boundary maps of chain complexes, i.e.
$$F^\pm\circ D^\pm= D^{F}\circ F^\pm,$$ and thus induce
homomorphisms of homology
 $$F^\pm_*:H_*(\fW^{\pm},D^\pm)
 \to H_*(\fD,D^{F}),$$
 as modules over $H_*(\fW^{\pm},D^\pm)$.
 \end{description}
 \end{proposition}

\begin{theorem}\label{thm:SFT-cobordism}
The   potential $ \bF\in\hbar^{-1}\fD$  defined  above by the formula
(\ref{eq:GW-potential})
satisfies the equation
\begin{equation}\label{eq:SFT-cycle}
D_We^{\bF}=0 ,
\end{equation} and hence all conclusions of Proposition \ref{prop:SFT-cobordism}
hold for $\bF$.
 \end{theorem}

  The appearance of $e^{\bF}$ in  equation
(\ref{eq:SFT-cycle}) has the following reason.
   Similar to Theorem
\ref{thm:HH} above, equation (\ref{eq:SFT-cycle}) follows from the
description of codimension $1$ strata on the boundary of the
moduli space of holomorphic curves in the cobordism $W$, see
 Proposition \ref{prop:boundary} above.   Notice
that $e^{\bF}$ is the generating function  counting possibly disconnected
 holomorphic curves in $W$. Thus the identity
$$\overrightarrow{\bH^-} e^{\bF}- e^{\bF}\overleftarrow{\bH^+}=0$$
asserts, in agreement with Proposition \ref{prop:boundary},
that the  codimension $1$ strata on the  boundary of the moduli space
$\wt\Mc(W)$
of not necessarily connected curves
in $W$ correspond
  to  stable
curves $(f_1,f_2)$ of height $2$, where one of the  curves $f_1,f_2$ belongs to $\wt\Mc(W)$,
 while the second one
is contained in the symplectization of $V^\pm$  and has
precisely one component different from the   straight cylinder
over a periodic orbit from $\Pc^\pm$.

\medskip

\begin{remark}\label{rem:non-closed2}
{\rm
(Comp. Remark \ref{rem:non-closed1} above) The   potential
$\bF$, extended to the space of all, not necessarily closed
 differential forms
 satisfies the   equation
\begin{equation}\label{eq:non-closed2}
d(e^\bF)=D_We^{\bF},
\end{equation}
where $d$ is the de Rham differential.
This equation generalizes
equation   (\ref{eq:SFT-cycle}).
}
\end{remark}
 Following the scheme of Section \ref{sec:3algebras} above we will associate
 now with the cobordism $W$ two other left-right modules, one over
 the
 Poisson algebras $\fP^\pm$, and another over the graded differential
 algebras $\fA^\pm$.

 Consider the graded commutative algebra  $\fL=\fL({W,\a^\pm })$
  of
power series of the form
\begin{equation}
\sum\limits_{\Gamma,d}\varphi_{\Gamma,d} (q^-,t ) z^d(p^+)^{\Gamma},
\end{equation}
where $\varphi_{\Gamma,d}$ are polynomials of $q_\g^-$ and formal power seies of $t_i$, while $\Gamma$ and
$d$  satisfies the  above inequality (\ref{eq:Novikov}).
Let us also consider the larger graded commutative algebra
\begin{equation}\label{eq:fL}
\widehat{\fL}=\{\sum\limits_{\Gamma^+,\Gamma^-,d}
\varphi_{\Gamma^+,\Gamma^-,d} (q^-,q^+,t )
z^d(p^+)^{\Gamma}(p^-)^{\Gamma'}\},
\end{equation}
where the Novikov condition (\ref{eq:Novikov}) is satisfied for both pairs
$(d,\Gamma)$ and $(d,\Gamma')$.
   The algebra $\widehat{\fL}$
 has a natural Poisson  bracket so that the homomorphisms
$f\mapsto \widehat{f}$,
 where we  denote by $ \widehat{f}$
the image in $\widehat{\fL}$ of an element $f\in\fP^\pm$ under the
coefficient homomorphism, are Poisson homomorphisms.  We set
$\widehat{\bh}=\widehat{\bh^-}-\widehat{\bh^+}$, and   for any $f\in\fL$ denote
by $L_{f}$ the ``Lagrangian variety''
 $$\{p^-_\g=\kappa_\g\frac{\partial
\widehat{f}}{\partial q^-_\g}, q^+_\g=\kappa_\g\frac{\partial
{\widehat f }}{\partial p^+_\g}\}. $$ Strictly speaking $L_f$ is
an ideal in the Poisson algebra $\widehat{L}$. However, it is
useful to think about $L_f$ as a Lagrangian variety in the
symplectic super-space $(\bV^-)\oplus\bV^+$ with the symplectic
form $$\sum \kappa^{-1}_{\g^-}dp^-_{\g^-}\wedge
dq^-_{\g^-}+\kappa^{-1}_{\g^+}dq^+_{\g^+}\wedge dp^+_{\g^+}\,,$$
and with an appropriate change of the coefficient ring.

For any function $f\in\fL$, which satisfies the Hamilton-Jacobi
equation
\begin{equation}\label{eq:HJ1}
\widehat\bh|_{L_f}=0
\end{equation}
 the Hamiltonian vector
field defined by the Hamiltonian $\widehat{\bh}$ is tangent to
$L_{f}$, and hence
  the differential $d^f:\fL\to\fL$, defined by the formula
  $$d^f(g)=\{\widehat\bh,g\}|_{L_f}$$ has the following meaning: we identify
  $\fL$ with the space of functions on $L_f$ and differentiate
  them along the Hamiltonian vector field determined by $\widehat\bh$.

Here is an analog of Proposition \ref{prop:RSFT-cobordism}
for the algebra
$\fL$.
\begin{proposition}\label{prop:RSFT-cobordism}
 Suppose that    $\widehat\bh|_{L_f}=0$.
 Then
 \begin{description}
 \item{1.} $(d^f)^2=0;$
 \item{2.} The maps $f^\pm:\fP^\pm\to\fL$, defined by the
formula $g\mapsto {\widehat{g}}|_{L_{f}}$, are homomorphisms of
chain complexes, i.e. $$d^{f}\circ f^\pm=f^\pm\circ
d^\pm;$$
\item{3.} If  $g_1,g_2\in\fP^\pm$ Poisson commute with $\bh^\pm$
 or, in other words, if $g_1,g_2\in{\mathrm {Ker}}\,d^\pm$
then
\begin{equation*}
\{f^\pm(g_1),f^\pm(g_2)\}=f^\pm(\{g_1,g_2\}).
\end{equation*}
where the left-side Poisson bracket is taken  in the algebra $\widehat{\fL}$.
\end{description}
\end{proposition}
Let us recall that $\bF=\bF_W\in\fD$ has the form
$\bF=\hbar^{-1}\sum\limits_{g=0}^\infty\bF_g\hbar^g$.
Again, to simplify the notation we will write $\bff $ instead of
$\bF_0$.
The following theorem is the reduction of Theorem
\ref{thm:SFT-cobordism} to the level of rational Gromov-Witten theory.
\begin{theorem}\label{thm:RSFT-cobordism}
The series $\bff(q^-,p^+,t)$ belongs to the algebra $\fL$ and satisfies
the  equation
\begin{equation}\label{eq:Lagr}
\widehat{\bh}|_{L_{\bff}}=0.
\end{equation}
\end{theorem}
In particular, all statements of  the above Proposition
\ref{prop:RSFT-cobordism} hold for $\bff$, and  this allows
us    to   define    the homology
\begin{equation*}
{}_{\bff}H^{\RSFT}_*(W|J,\a^\pm)=H_*(\fL,d^{\bff}).
\end{equation*}
  The chain maps $\bff^{\pm}$ induce homomorphism of Poisson homology
algebras
\begin{equation}
(\bff^\pm)_*:H^{\RSFT}_*(V^\pm|J,\a^\pm)=H_*(\fP^\pm,d^\pm)\to
{}_{\bff}H^{\RSFT}_*(W|J,\a^\pm)\,.
\end{equation}

\bigskip
For the rest of this section we assume that $W$ is a  rational
homology
cobordism, i.e. the restriction maps $$H^*(V^-;\R)\leftarrow
H^*(W;\R)\rightarrow H^*(V^+,\R)$$ are isomorphisms.  Equivalently,
this means that
  the inclusions $V^\pm\to W$ induce isomorphisms of rational
  homology groups.

The potential $\bff\in\fL$ which we defined above can be written in the
  form
\begin{equation}
\bff=\sum\limits_i\sum\limits_{|\Gamma^+|=i}f_{\Gamma}^i(q^-,t)(p^+)^{\Gamma^+}.
\end{equation}
                        Notice that the assumption that $W$
                         is a homology cobordism implies
  that $f^0(q^-,t)$ is independent of $q^-$.
 Let us now define a homomorphism
 $\bPsi : \fA^+\to \fA^-$
 by the formula
 \begin{equation}
 \bPsi (q^+_{\g } )=f^1_{ \g }(q^-,t)\in\fA^-
 \end{equation}
on the generators $q^+_\g,\g\in\Pc^+$,
of the algebra $\fA^+$ and then   extend by
linearity.

\begin{theorem}   The homomorphism
$\bPsi : \fA^+\to \fA^-$ commutes with the boundary operators
$\partial^\pm$, i.e. $\partial^-\circ\bPsi =\bPsi
\circ\partial^+$, and in particular defines a homomorphism of
homology algebras
\begin{equation*}
(\bPsi )_*:H_*(\fA^+,\partial^+)\to H_*(\fA^-,\partial^-).
\end{equation*}
\end{theorem}

Without the assumption that $W$ is a homology cobordism
one gets only a correspondence between the algebras $\fA^+$ and $\fA^-$, similar to the
``semi-classical" case considered above.

\subsection{Chain homotopy}\label{sec:chain}
Let $W=\ora{V^-V^+}$ be a directed symplectic cobordism with fixed
contact forms $\a^\pm$ on $V^\pm$. We will discuss in this section
how the function $\bF=\bF_{W,J,\a^\pm}(p^+,q^-,t)$ and other
associated structures change when one replaces  $J$ with another
compatible almost complex structure $J'$ and replaces $t$ with a
form $t'=t+d\theta$ where $\theta$ has   compact support in $W$.

Let us begin with some algebraic preliminaries. Two series
$F_0,F_1\in\hbar^{-1}\fD$ are called {\it homotopic}, if they can
be included into  a family $F_s\in\hbar^{-1}\fD,\,s\in[0,1]$,
which satisfies the following differential equation
\begin{equation} \label{eq:homotopy}
    \frac{d F_s}{ds }=
    e^{-F_s}
    \left(\ora{[\bH^-,K_s]}e^{F_s}+
     e^{F_s}\overleftarrow{[
     K_s,\bH^+]}\right),\,\,s\in[0,1],
    \end{equation}
     for a family $K_s\in\hbar^{-1}\fD$.
     Here  $[\bH^-,K_s]$ and $[K_s ,\bH^+]$ are commutators in the
algebra $\hbar^{-1}\widehat{\fD}$, defined similar to
$\widehat{\fL}$ in (\ref{eq:fL}) above, i.e.
\begin{equation}
\widehat{\fD}=\{\sum\limits_{\Gamma^+,\Gamma^-,d,g}
f_{\Gamma^+,\Gamma^-,d} (q^-,q^+,t )
z^d\hbar^g(p^+)^{\Gamma^+}(p^-)^{\Gamma^-}\},
\end{equation}
where the Novikov condition (\ref{eq:Novikov}) is satisfied for
both pairs $(d,\Gamma^+)$ and $(d,\Gamma^-)$.
      In other words, we view $K_s$ as an operator on
       $\hbar^{-1}\fD$, acting by the multiplication
    by the series $K_s$,
    and  view $\bH^-$ and $ \bH^+  $
    as left and right differential operators.

Notice that the family $ K_s\in\hbar^{-1}\fD, s\in [0,1],$
defines a flow $\Phi^s=\Phi_K^s: \fD\fD\to\fD \fD$,
  by a differential equation
   \begin{equation} \label{eq:chain-flow}
    \frac{d \Phi^s(G) }{ds }=\cK_s\big(\Phi^s(G)\big),
 \end{equation}
 where we set
\begin{equation*}
    \cK_s(G)=\left(\ora{[\bH^-,K_s]}G+
    G\overleftarrow{[ K_s,\bH^+]}\right),\,\,s\in[0,1].
 \end{equation*}
The linear operators $\Phi^s$ preserve  the ``submanifold" $\fE=e^{
\hbar^{-1}\fD_{\even}}$, where $\fD_\even$ is the even part of
$\fD$, and we have $$\Phi^s(e^{F_0})=e^{F_s},$$ where
the family $F_s$ satisfies the equation (\ref{eq:homotopy}).

The tangent space to $\fE$ at a point $e^F$, $F\in \hbar^{-1}\fD_\even$, consists of
series $fe^F$, $f\in\hbar^{-1}\fD_\even$, and thus it is naturally parameterized by
$\hbar^{-1}\fD_\even$. With respect to this parameterization  the
differential  of the flow $\Phi^s|_\fE$ defines a family of maps
$  T_F^s:\hbar^{-1}\fD_\even\to\hbar^{-1}\fD_\even$, $F\in \hbar^{-1}\fD_\even$, by the formula
 \begin{equation}\label{eq:lin-flow}
    T_F^s(f)=e^{-F_s}\Phi^s(fe^F),\;\;\hbox{ where
    }\;\;F_s=\Phi^s(F);
  \end{equation}
  We extend   $T_F^s$ to the whole
   $\hbar^{-1}\fD$ by the same formula
   (\ref{eq:lin-flow}).
Let us list some properties of the flows $\Phi^s$ and $T^s_F$
\begin{proposition}\label{prop:SFT-chain}
Suppose that    for an  element  $F\in\hbar^{-1}\fD_\even$ we
 have $$D_W(e^F)=\ora{\bH^-}e^F -e^F\ola{\bH^+}=0.$$
    Then
\begin{description}
\item{1.}
       The   flow  $T^s_F :\hbar^{-1}\fD\to\hbar^{-1}\fD$
satisfies the equation
\begin{equation}
 T_F^s\circ D^{F}=D^{F_s}\circ T_F^s.
  \end{equation}
  for all $s\in[0,1].$
  In particular,    $D_W(F_s)=0$  for all $s\in[0,1]$, and $T_F^s$
       defines a family of isomorphisms
$ H_*(\fD,D^{F})\to
 H_*(\fD,D^{F_s})$.
\item{2.} The homology class
$[e^{F_s}]\in H_*(\fD \fD,D_W)$ is independent of  $s$.
                    \item{3.} The diagram
\bigskip\bigskip

\centerline{$ \mathop{\fD\qquad}\limits_{\;
 \; \nwarrow \, F^\pm}
\mathop{\longrightarrow}\limits^{T_F^s}
\mathop{\fD}\limits_{\nearrow \, F_s^\pm}$}

\centerline{$\fW^\pm$}

homotopically commutes, i.e. there exist operators
$A_s^\pm:\fW^\pm\to\fD$,
such that
\begin{equation}\label{eq:SFT-diagram}
(T_F^s)^{-1}\circ F_s^\pm- F^\pm=D^{F}\circ A_s^\pm+ A_s^\pm\circ
D^\pm\,.
\end{equation}
In particular,
this diagram commutes on the level of homology algebras.
  \end{description}
  \end{proposition}
The proof of this proposition is a  straightforward computation by
differentiating  the corresponding equations.
To illustrate the argument, let us   verify
(\ref{eq:SFT-diagram})
in \ref{prop:SFT-chain}.

Take, for instance, $f\in\fW^-$ and
set
 \begin{equation*}
\begin{split}
G_s(f)=&(T_F^s)^{-1}(e^{-F_s}\ora fe^{F_s}))-e^{-F}\ora fe^F\\&=
e^{-F}(\Phi^s)^{-1}\left(\ora fe^{F_s}\right)-e^{-F}\ora fe^F.\\
\end{split}
\end{equation*}

Then we have $G_0(f)=0 $ and

\begin{equation*}
\begin{split}
\frac{dG_s(f)}{ds}
&= -e^{-F}\left((\Phi^s)^{-1}
 (  \cK_s \ora{f}  e^{F_s }  )\right)
+  e^{-F}
(\Phi^s)^{-1}\left(
  \overrightarrow{f} \cK_s\left( e^{F_s }\right)\right)
 \\
& =           e^{-F}(\Phi^s)^{-1}\left(
           \left[\cK_s,\overrightarrow{f}\right]
            e^{F} \right).  \\
\end{split}
\end{equation*}
Now remember that $\cK_s=[\widehat{\bH}  ,K_s]$,
 where
$\widehat{\bH}=\overrightarrow{\bH^-}-\overleftarrow
{\bH^+ }$, and using the Jacobi identity we get
\begin{equation*}
\begin{split}
  \frac{dG_s(f)}{ds}
  &=  e^{-F} (\Phi^s)^{-1}\left(
           \left[[\widehat{\bH}  ,K_s],
           \overrightarrow{f}\right]
            e^{F_s}\right) \\
      & =
 e^{-F}(\Phi^s)^{-1}\left(
           \left[[K,\overrightarrow{f}],
           \widehat{\bH}  \right]
            e^{F_s }  \right)\\
            &+
              e^{-F} (\Phi^s)^{-1}\left(
           \left[[\overrightarrow{f},
           \widehat{\bH}  ],K_s\right]
             e^{F_s } \right). \\
                           \end{split}
         \end{equation*}
         Let us define now a linear operator
         $B_s:\hbar^{-1}\fW_-\to\hbar^{-1}\fD $ by the formula
         \begin{equation}
         B_s(g)=   e^{-F}(\Phi^s)^{-1}\left(
         \overrightarrow{[g,K_s] }
          e^{F_s}\right).
         \end{equation}
          for $g\in \hbar^{-1}\fW_-$.
        Recall that
         $D_W(e^F)=\widehat{\bH}e^{F}=0$, and observe that
         $\widehat{\bH}$ commutes with
         $\cK_s=[\widehat{\bH},K_s]$ because $\widehat{\bH}\circ
          \widehat{\bH}=0$.  We also have
                    $[\widehat{\bH},f]= D^-f$ for
                     $f\in\hbar^{-1}\fW_-$.
          Hence
           $ \frac{dG_s(f)}{ds}$ can be rewritten as
          \begin{equation}\label{eq:derivative}
          \frac{dG_s(f)}{ds}= D^{F}(B_s(f))    +
            B_s(D^-(f))
          \end{equation}

        Finally we integrate $B_s$ into the required linear operator
         $A^-:\hbar^{-1}\fW^-\to\hbar^{-1}\fD$:
         \begin{equation}
        A_s^-(g)=\int\limits_0^s B_s(g)ds,
        \quad\hbox{for}\quad g\in\hbar^{-1}\fW^-.
        \end{equation}

         In view of   equation
         (\ref{eq:derivative})
         we get
  \begin{equation}
(\Phi^s)^{-1}\circ(F^s)^- -(F)^-=D^F\circ A_s^-+
A_s^-\circ D^-\,.
\end{equation}
 \qed

\bigskip
 Let us consider now a generic family $J^\tau$, $\tau\in[0,1]$,     of
compatible almost complex structures on $W$ connecting $J^0=J$
with $J^1=J'$. We  assume that the deformation $J^\tau$ is
fixed outside of a compact subset of $W$.

\medskip
Set
\begin{equation}
\Mc_{g,r,s^+,s^-}(W,\{J^\tau\}) =\bigcup\limits_{\tau\in[0,1]}
\Mc_{g,r,s^+,s^-}(W,J^\tau).
\end{equation}

The evaluation maps   defined for each $\tau$ can be
combined into the evaluation map
\begin{equation}
ev:\Mc_{g,r,s^+,s^-}(W,\{J^\tau\})\to
(W\times I)^r\times\left(\Pc^-\right)^{s^-}\times\left(\Pc^+\right)^{s^+}.
\end{equation}

Consider  closed forms
$\widehat {\theta}_1,\dots \widehat{\theta}_r$    on
  $W\times I$, such that
  $\widehat{\theta}_i=\wt{\theta}_i+d\beta_i$, $i=1,\dots, r$, where
  $\wt{\theta}_i$ is the pull-back of  a form $\theta_i$ on $W$
  with cylindrical ends, and
  $\beta_i$ has   compact support in $W\times \Int I$.

Similarly to correlators  of degree $-1$ and $0$
(see \ref{sec:correlators} and \ref{sec:sympl-eval})
we  can define correlators of degree $1$, or
$1$-parametric correlators by the formula
\begin{eqnarray}
 \up1\langle\widehat {\theta}_1,\dots \widehat{\theta}_r\,;\theta^-_1,\dots,
\theta^-_{s^-};\theta^+_1,\dots,
\theta^+_{s^+} \rangle^A_{ g}  =
  \int\limits_{\overline{\Mc_{g,r,s^+,s^-}(W,\{J^\tau\})} }\\
\nonumber   ev^*(\widehat{\theta}_1\otimes\dots
\otimes \widehat{\theta}_r
\otimes\theta^-_1\otimes\dots\otimes\theta^-_{s^-}
\otimes\theta^+_1\otimes\dots\otimes\theta^+_{s^+} ),
\end{eqnarray}
  for a homology class $A\in H_2(W)$, and cohomology
classes $$\theta^\pm_i\in H^*(\Pc^\pm),\; i=1,\dots,s^\pm.$$

Consider a closed form $T=\widetilde{t}+d\beta$,
 where  the notation $\wt t$ and $\beta$
  have the same meaning as
 above, i.e. $\wt t$ is the pull-back of a form $t$
 on $W$ with cylindrical ends, and $\beta$ has
   compact support in $ W\times\Int I$.
We can organize the correlators
 $$\up1\langle T,\dots,T\,; q^-,\dots,q^-,p^+,\dots, p^+\rangle^A_{
 g}$$
    into a generating function
$\bK=\frac 1{\hbar}\sum\limits_{g=0}^\infty\bK_g\hbar^g\in\frac 1{\hbar}\fD$, defined
by the formula
\begin{equation}
\bK=\sum\limits_d\sum\limits_{g,r
,s^\pm=0}^\infty
\frac{1}{r!s^+!s^-!}\up1\langle \underbrace{T,\dots,T}_r;
 \underbrace{q^-,\dots,q^-
 }_{s^-};\underbrace{p^+,\dots,p^+}_{s^+}\rangle^d_{ g}\hbar^gz^d\,
 \end{equation}

  \medskip
 Let us define an operator $\cK:\fD\fD\to\fD\fD$
 by the formula
 \begin{equation}
 \cK(G)=  \overrightarrow{ [\bH^+,\bK]}G +G
\overleftarrow{[\bK ,\bH^-]},\; G\in\fD\fD .
  \end{equation}

Similar to Theorems \ref{thm:HH} and \ref{thm:SFT-cobordism}
the next theorem  can be viewed
as an algebraic    description of the boundary of the compactified
moduli
space        $\overline{\Mc_{g,r,s^+,s^-}(W,\{J^\tau\})}$.
\begin{theorem}\label{thm:SFT-chain}
For a generic family $J_\tau,\tau\in[0,1]$, of compatible almost complex
structures on $W$ we have
\begin{equation}\label{eq:SFT-chain}
e^{\bF^1}= e^{\cK }(e^{\bF^0}),
\end{equation}
where
$\bF^0=\bF_{W,J_0}(T|_{W\times 0}),\,\bF^1 =\bF_{W,J_1}(T|_{W\times 1}),$ and
$\cK=\cK(T).$
\end{theorem}
 Notice that if we define $\bF^s$,
 $s\in [0,1]$, by the formula
\begin{equation}
e^{\bF^s}= e^{s\cK }(e^{\bF^0}),
\end{equation}
then the flow $\Phi^s(F)=F^s$   satisfies the differential
equation (\ref{eq:chain-flow})
 with $K(s)\equiv \bK$. Hence $\bPhi^0$ and $\bPhi^1$
 are homotopic, and therefore
  Theorem
\ref{thm:SFT-chain} and
Proposition \ref{prop:SFT-chain} imply
\begin{corollary}\label{cor:SFT-chain}
\begin{description}
\item{1.} The homology class $[e^{\bF}]\in H_*(\fD
\fD,D_W)$ is
independent of the choice of a compatible almost complex structure
$J$ on $W$ and of the choice of
$t\,{\mathrm {mod}}\left(d\left(\Omega_{\mathrm
{comp}}(W)\right)\right)$, where $ {\Omega}_{\mathrm
{comp}}(W)$ is the space of forms with compact support.
\item{2.} For a generic compatible deformation $J^\tau$, $\tau\in[0,1]$,
the isomorphism  $T :\fD\to\fD$
   defined by the formula
\begin{equation}
    T (f)= e^{-\bF_1}
e^{\cK}(fe^{\bF^0 } )
\end{equation}
 satisfies the equation
\begin{equation}
 T \circ D^{\bF^0}=D^{\bF^1}\circ T  ,
 \end{equation}
 and thus defines an isomorphism
$ H_*(\fD,D^{\bF^0})\to
 H_*(\fD,D^{\bF^1})$.
\item{3.} The diagram
\bigskip\bigskip

\centerline{$ \mathop{\fD\qquad}\limits_{\;
 \; \nwarrow \, \left(\bF^0\right)^\pm}
\mathop{\longrightarrow}\limits^{T }
\mathop{\fD}\limits_{\nearrow \, \left(\bF^1\right)^\pm}$}

\centerline{$\fW^\pm$}

homotopically commutes, i.e. there exist operators
$A^\pm:\fW^\pm\to\fD$,
such that
\begin{equation}
T^{-1} \circ\left(\bF^1\right)^\pm-
\left(\bF^0\right)^\pm=D^{\bF^0}\circ A^\pm+
A^\pm\circ D^\pm\,.
\end{equation}
  \end{description}
\end{corollary}

\bigskip
Consider now the equivalence relation for
 rational   potentials.

 Two series
 $ f_0,f_1\in\fL$ are called {\it homotopic}
if there exist families
$f_s, k_s\in\fL$, $s\in[0,1]$, so that the family
 $f_s$ connects $f_0$ and $f_1$ and satisfies
the following Hamilton-Jacobi differential
equation
\begin{equation}\label{eq:HJ}
 \frac{\partial f_s( p^+,q^-,t)}{\partial s}=
 \bG(p^+,\frac{\partial f_s ( p^+,q^-,t)}
 {\partial p^+},
 \frac{\partial f_s ( p^+,q^-,t)}
 {\partial q^-} ,q^-),
\end{equation}
 where
  \begin{eqnarray}
 \nonumber \bG(p^+,q^+,p^-,q^-,t)=\{\bh^+-\bh^-,k_s\}=\\
 \sum\limits_{
 \g^+\in\Pc^+,\g^-\in\Pc^-} \kappa_{\g^-}\frac{\partial \bh^-(p^-,q^-,t)}
 {\partial p^-_{\g^-}}
 \frac {\partial k_s( p^+,q^-,t)}{\partial  q^-_{\g^-}}+\\
  \nonumber \kappa_{\g^+}\frac{\partial k_s( p^+,q^-,t)}{
  \partial  p^+_{\g^+}}\frac {\partial \bh^+(p^+,q^+,t)}
 {\partial q^+_{\g^+}} \,.
 \end{eqnarray}
  Here $\kappa_\g$ denotes, as usual,  the multiplicity of
  $\g$.

We can view  the correspondence
  \begin{equation*}
  f ( p^+,q^-,t)\mapsto f_s( p^+,q^-,t),
  \end{equation*}
  where $f_s$ is the solution of
  the above equation (\ref{eq:HJ}) with the initial data
  $f_0=f$,
  as a non-linear operator $ S^s:\fL\to\fL$.
  Let us denote by
  $T^s_f $ the linearization of $ S^s$ at a point
  $f$.
  The next proposition is a rational version of Proposition
  \ref{prop:SFT-chain}. It can be either deduced from
\ref{prop:SFT-chain}, or similarly  proven by    differentiation
with respect to the parameter $s$.
  Denote by $\cS$ the subspace of $\fL$ which consists of solutions of
    the Hamilton-Jacobi equation
(\ref{eq:HJ1}), i.e.
 \begin{equation*}
 \widehat\bh|_{L_f}=0,\quad\hbox{  where}\quad
  L_f=\{p^-_\g=\kappa_\g\frac{\partial
 {f}}{\partial q^-_\g}, q^+_\g=\kappa_\g\frac{\partial
{  f }}{\partial p^+_\g}\}.
\end{equation*}

\begin{proposition}\label{prop:RSFT-chain}
  \begin{description}
  \item{1.} The subspace  $\cS\subset\fL$ is invariant
    under the flow $S^s$.
  \item{2.}   For $f\in \cS$
  the isomorphism  $T^s_f:\fL\to\fL$
 satisfies the equation
\begin{equation}
 T^s_f\circ d^{f }=d^{S^s(f)}\circ T^s,
 \end{equation}
 and thus defines an isomorphism
$H_*(\fL,d^{f})\to
H^*(\fL,d^{S^s(f)} )$.
  \item{3.}  For $f\in\cS$ the diagram
   \bigskip\bigskip

\centerline{$ \mathop{\fL\qquad}\limits_{\;
 \; \nwarrow \, \left(f\right)^\pm}
\mathop{\longrightarrow}\limits^{T^s_f}
\mathop{\fL}\limits_{\nearrow \, \left(S^s(f)\right)^\pm}$}

\centerline{$\fP^\pm\;$}

\bigskip

  homotopically commutes.
  \end{description}
\end{proposition}

 Theorem     \ref{thm:SFT-chain} reduces on the level of
 rational curves
 to the following

\begin{theorem}\label{thm:RSFT-chain}
Let $\bF^0$, $\bF^1$ and $\bK$ be as in Theorem
\ref{thm:SFT-chain}. Set $\bff^0=\bF^0_0,\;\bff^1=\bF^1_0,\;
\bk=\bK_0.$ Then
$\bff_0$ and $\bff_1$ are homotopic, i.e. they can be included into a family
$\bff_s,\,s\in [0,1]$, such that the Hamilton-Jacobi equation
(\ref{eq:HJ}) holds with $f_s=\bff_s$ and  $k_s\equiv\bk$.
\end{theorem}

Hence, Proposition \ref{prop:RSFT-chain} implies
\begin{corollary}\label{cor:RSFT-chain}
 For a generic compatible deformation $J_s$,
$s\in[0,1]$, we have
\begin{description}
\item{1.}
 The operator  $S^1:\fL\to\fL$ defines an automorphism of
 the space  of solutions of (\ref{eq:HJ1});
    \item{2.} The linearization $T^1 =T^1_{\bff^0}$ of $S^1$ at the point
    $\bff_0$ satisfies the
  equation
\begin{equation}
 T^1\circ d^{\bff^0}=d^{\bff^1}\circ T^1,
 \end{equation}
 and thus defines an isomorphism
$H_*(\fL,d^{\bff^0})\to
H^*(\fL,d^{\bff^1})$.
  \item{2.} The diagram
   \bigskip\bigskip

\centerline{$ \mathop{\fL\qquad}\limits_{\;
 \; \nwarrow \, \left(\bff^0\right)^\pm}
\mathop{\longrightarrow}\limits^{T^1}
\mathop{\fL}\limits_{\nearrow \, \left(\bff^1\right)^\pm}$}

\centerline{$\fP^\pm$}

\bigskip

  homotopically commutes.
  \end{description}
\end{corollary}

To formulate the ``classical level'' corollary of Theorem
\ref{thm:RSFT-chain}
we assume, as usual, that $W$ is a homology cobordism.
\begin{theorem}\label{thm:classical-chain}
The homomorphisms  $\Psi^1_{J_0} ,\Psi^1_{J_1}: \fA^+\to\fA^-$ associated
to two compatible almost complex structures $J_0$ and $J_1$ are
homotopic, i.e. there exists a map $\Delta:\fA^+\to\fA^-$ such
that
$$\Psi^1_{J_1} -\Psi^1_{J_1} =\partial^-\circ\Delta+\Delta\circ\partial^+.$$
      \end{theorem}
The map $\Delta$ can be
expressed through the function
  $\bk\in\fL$. However, unlike the case of usual Floer homology
  theory, $\Delta$ and $\bk$ are related via a first order non-linear PDE
  (which can be deduced from the equation (\ref{eq:HJ})),
     and thus  one cannot write a general explicit formula  relating
     $\Delta$ and $\bk$.

     \subsection{Composition of cobordisms}\label{sec:composition}
In this section we   study the behavior of
potentials and associated algebraic structures under
the operation of composition of directed symplectic cobordisms.

    Let    us recall (see Section \ref{sec:splitting}) that   given
    a dividing contact type hypersurface $V$ in a
    directed symplectic cobordism $W=\ora{V^-V^+}$ one can split $W$
    into a composition $W=W_-\circledcirc W_+$ of cobordisms $W_-=\ora{V^-V}$ and
    $W_+=\ora{VV^+}$. From the point of view of an almost
     complex structure        the process of splitting consists of deforming
     an original almost complex structure $J=J^0$ to an
     almost complex structure $J^\infty$, such that the restrictions
     $J_\pm=J^\infty|_{W_\pm}$   are compatible  with the structure
     of (completed) directed symplectic cobordisms $W_\pm$.

Conversely,  directed symplectic cobordisms
 $W_-=\ora{V^-V}$ and
    $W_+=\ora{VV^+}$ with matching data
    on the common boundary    can be glued into
    a cobordism         $W=\ora{V^-V^+}$ in the following sense:
    there exists a family  $J^\tau$ of almost complex structures
    on $W$ which in the limit  splits $W$ into the
       composition of cobordisms $W_-=\ora{V^-V}$ and
    $W_+=\ora{VV^+}$.

    \bigskip
In order to write the formula relating   the potentials
of $W$ and $W_\pm$ we  first need to make more explicit the
relation  between $2$-dimensional homology classes realized  by
holomorphic curves in $W_\pm$ and $W$. We will keep
 assuming that there are no torsion
elements in $H_1$.
Let us recall (see Sections \ref{sec:dynamics} and
\ref{sec:holomorphic} above) that  we have chosen curves ${C}^i_-\subset
W_-,i=1,\dots, m_-,\;\;{C}^j_+\subset W_+,j=1,\dots,m_+,$
and ${C}^k
\subset W,k=1,\dots m,$ which represent bases of first homology of
the respective cobordisms. We also have chosen
for every periodic orbit $\g\in
\Pc_\a$ of the Reeb field $R_\a$ on $V$ a surface
$F^\g_\pm$ realizing homology in $W_\pm$ between $\g$ and a linear
combination of basic curves ${C}^i_\pm$.
For our current purposes we have to make one extra choice: for
each curve $ {C}^i_\pm$ we choose a surface $S^i_\pm$ which
realizes homology in $W$ between ${C}^i_\pm$ and the
corresponding linear combination of the curves ${C}^1,\dots,
{C}^m$.
All the choices enable us to associate with every orbit
$\g\in\Pc_\a$ a homology class $C^\g$ which is realized by the chain
 $$F^\g_+-F^\g_-+\sum\limits_1^{m_+} n^+_jS^j_+-\sum\limits_1^{m_-}
n^-_iS^i_-\,,$$
where $\partial F^\g_\pm=[\gamma]-\sum\limits_1^{m_\pm}
n^{\pm}_j{C}^j_\pm$. We will denote by $d^\g$ the degree of
$C^\g$, i.e. the string of its coordinates in the
chosen basis
$A_1,\dots,A_N\in H_2(W)$.

    Let us define an operation
 \begin{equation}
 \star:\fD\fD_-\otimes \fD\fD_+ \to \fD\fD,
 \end{equation}
  where   $\fD\fD_\pm=\fD\fD_{W_\pm}$ and
$\fD\fD= \fD\fD_{W} $.
  For
   $F=\sum\limits_{\Gamma} f_{\Gamma}(t^-,q^-,\hbar,z^-)
   p^{\Gamma}\in\fD\fD_-$ and
    $
  G=\sum\limits_{\Gamma^+} g_{\Gamma^+}
 (t^+,q ,\hbar,z^+)(p^+)^{\Gamma^+} \in\fD\fD_+$ we set
  \begin{equation}
 \begin{split}
  &F  {\star} G(t,q^-,p^+,\hbar,z)=\\
  &\Big(
 \sum\limits_{\Gamma }\wt f_\Gamma(t,q^-,\hbar,z)
 \hbar^s  \prod\limits_{i=1}^s
 \kappa_{\g_{i }}z^{d^{\g_i}}
 \overrightarrow{\frac{\partial}{\partial
 q_{\g_i}}}
\sum\limits_{\Gamma^+} \wt g_{\Gamma^+}
(t,q ,\hbar,z)(p^+)^{\Gamma^+}\Big)\Big|_{q=0}\,. \\
\end{split}
\end{equation}
Here we denote by $\wt f_\Gamma$ and $\wt
g_{\Gamma^+}$ the images
of $ f_\Gamma$ and $ g_{\Gamma^+}$ under the coefficient
homomorphisms $H_2(W_\pm)\to H_2(W)$.
Let us explain what happens with the
 variables $t$ and $z$ in more details.  Let
  $A_1^\pm,\dots, A^\pm_{N^\pm},$
  and
   $A_1,\dots, A_{N},$
    be   the chosen bases in     $H_2(W_\pm)$ and $ H_2(W)$.
                   Then we have $$i^\pm_*(A^\pm_k)=\sum\limits_{j=1}^{N^\pm}n^\pm_{kj} A_j,$$
   where $k=1,\dots, N$,  $(n^\pm_{kj})$ are integer matrices,
   and $i^\pm:W_\pm\to W$ the inclusion maps.

We have
 \begin{equation*}\begin{split}
  & f_{\Gamma}(t^-,q^-,\hbar,z^-)
   = \sum\limits_{d=(d_1,\dots, d_{N^-})}f_{\Gamma,d}(t^-,q^-,\hbar)(z_1^-)^{d_1}\dots
   (z_{N^-}^-)^{d_{N^-}},\\
  &g_{\Gamma^+}
 (t^+,q ,\hbar,z^+) = \sum\limits_{d=(d_1,\dots, d_{N^+})}g_{\Gamma^+,d}(t^+,q,\hbar)(z_1^+)^{d_1}\dots
   (z_{N^+}^+)^{d_{N^+}}  ,\\
 \end{split}
 \end{equation*} where
    we denote by $z_\pm$ the ``$z$-variables" in $W_\pm$.
    Then
    \begin{equation*}
    \begin{split}
   & \wt f_\Gamma(t,q^-,\hbar,z)=
          \sum\limits_{d=(d_1,\dots,d_{N_-})} \wt f_{\Gamma,d}(t|_{W_-},q^-,\hbar)z_1^{M_1^-}
         \dots z_N^{M^-_{N}}, \\
     & \wt g_{\Gamma^+}(t,q,\hbar,z)=
          \sum\limits_{d=(d_1,\dots,d_{N_+})} \wt g_{\Gamma^+,d}(t|_{W_+},q^-,\hbar)z_1^{M_1^+}
         \dots z_N^{M^+_{N}}, \\
          \end{split}
          \end{equation*}
        where $M^\pm_j=   \sum\limits_{k=1}^{N^\pm}  n^\pm_{kj}d_k$, $j=1,\dots,N$.

   \bigskip
  Let us  observe
  \begin{lemma}\label{lm:star}
  \begin{description}
  \item{1.} The operation $\star$ is associative.
  \item{2.}   For $F\in\hbar^{-1 }\fD_-,G\in\hbar^{-1 }\fD_+$ there exists a unique function $H\in \hbar^{-1 }\fD$,
  such that
 $$
  e^H=e^F\star e^G.
 $$
 \end{description}
  \end{lemma}
  We will denote this $H$ by $F\lozenge G$, so that we have
  $e^{F\lozenge G}=e^F\star e^G$.
We will also consider the maps
$$\lozenge G:\hbar^{-1}\fD_-\to \hbar^{-1}\fD,\;\;\lozenge
G(F)=F\lozenge G,\,F\in\fD_-,$$
and
$$F\lozenge:\hbar^{-1}\fD_+\to \hbar^{-1}\fD,\;\;F\lozenge
 (G)=F\lozenge G,\,G\in\fD_+,$$
 and for even $F,G$ their linearizations:
              \begin{equation*}
 \begin{split}
& T_F(\lozenge G):\hbar^{-1}\fD_-\to \hbar^{-1}\fD,\\
 &T_F(\lozenge G)(f)=\frac{d\big(F+\varepsilon
 f)\lozenge G\big)}{d\varepsilon}\big|_{\varepsilon=0}=e^{-F\lozenge
 G}\big((fe^F)\star e^G\big),\;\; f\in\hbar^{-1} \fD_-,\\
 \end{split}
 \end{equation*}
and
\begin{equation*}
 \begin{split}
& T_G(F\lozenge  ):\hbar^{-1}\fD_+\to \hbar^{-1}\fD,\\
 &T_G(F\lozenge  )(g)=\frac{d\big(F \lozenge (G+\varepsilon
 g)\big)}{d\varepsilon}\big|_{\varepsilon=0}=e^{-F\lozenge
 G}\big( e^F \star (ge^G)\big),\;\; g\in\hbar^{-1} \fD_+,\\
 \end{split}
 \end{equation*}

Let us first formulate an algebraic
\begin{proposition}\label{prop:SFT-composition}
Suppose that $F\in \hbar^{-1}\fD_-$ and $G\in \hbar^{-1}\fD_+$ are
even elements, which
satisfy the equations
$$D_{W_-}(e^F)=\overrightarrow{\bH^-}e^F-e^F\overleftarrow{\bH}=0$$  and $$D_{W_+}(e^G)=
\overrightarrow{\bH}e^F-e^F\overleftarrow{\bH^+}=0,$$
where $\bH^\pm=\bH_{V^\pm},\,\bH=\bH_V$.
Then we have
\begin{description}
\item{1.} $D_W(e^{F\lozenge G})=\ora{\bH^-}\, e^{F\lozenge G}-e^{F\lozenge
G}\,\ola{\bH^+}=0$.
\item{2.} The homomorphisms
$ T_G(F\lozenge  ):\hbar^{-1}\fD_+\to \hbar^{-1}\fD$ and
$T_F(\lozenge G):\hbar^{-1}\fD_-\to \hbar^{-1}\fD$ satisfy
the equations $$T_G(F\lozenge  )\circ
D^{G}=D^{F\lozenge G}\circ T_G(F\lozenge )$$ and $$T_F(
\lozenge G )\circ D^{F}=D^{F\lozenge G}\circ T_F(
\lozenge G),$$ and in particular they define homomorphisms
of the corresponding homology algebras:
$$ \left( T_G(F\lozenge  )\right)_*:H_*( \fD_+,D^G)\to H_*(  \fD,D^
 {F\lozenge G})$$ and
$$ \left( T_F( \lozenge G  )\right)_*:H_*( \fD_-,D^G)\to H_*(  \fD,D^
 {F\lozenge G}).$$
\item{3.} $$ T_F( \lozenge G)\circ F^-   =
(F \lozenge G)^-$$ and $$  T_G( F\lozenge)
\circ G^+= (F \lozenge G)^+.$$
\item{4.}
Suppose we are given three cobordisms $W_1,W_2,W_3$ with matching ends
 so that we can  form  the  composition
$W_{123}=W_1\circledcirc W_2 \circledcirc W_3$, and series
$F_i\in\hbar^{-1}\fD_i=\hbar^{-1}\fD_{W_i},\,i=1,2,3$, such that
$D_{W_i}e^{F_i}=0$, $i=1,2,3$.
Then
\begin{equation}
T_{F_1\lozenge F_2}(\lozenge F_3)\circ T_{F_1}(\lozenge
F_2)=T_{F_1}(\lozenge (F_2\lozenge F_3)).
\end{equation}
\end{description}
\end{proposition}
The proof of this proposition is immediate from the definition of
the corresponding operations. For instance, to prove
\ref{prop:SFT-composition}.1 we write
\begin{equation*}
\begin{split}
 \ora{\bH^-}\,\big(e^{F\lozenge G}\big)-\big(e^{F\lozenge
G}\big)\,\ola{\bH^+}&=\ora{\bH^-}\,e^F\star e^G -e^F\star e^G\ola
{\bH^+}\\
&=\left(\ora{\bH^-}\,e^F\right)\star e^G -e^F\star\left( e^G\ola
{\bH^+}\right)\\
&=\left(e^F\,\ola{\bH}\right)\star e^G
-e^F\star\left( \ora{\bH}e^G\right).\\
\end{split}
\end{equation*}
To finish the argument
let us consider a cylindrical cobordism $W_0=V\times\R,$
take the function
$I=\sum \kappa_\g^{-1}p_\g q_\g$. Taking into account associativity
of $\star$ (see \ref{lm:star}) we have
\begin{equation}
\ora {f}e^G=(fe^I)\star e^G,\quad\hbox{and}\quad
 e^F\ola{f}=e^F\star (fe^I)\,.
 \end{equation}
 Hence, we have
 \begin{equation*}
 \begin{split}
 & \left(e^F\,\ola{\bH}\right)\star e^G
-e^F\star\left( \ora{\bH}e^G\right)\\
&=e^F\star(\bH e^I)\star e^G-  e^F\star(\bH e^I)\star e^G=0.\\
 \end{split}
\end{equation*}
\bigskip

Any cohomology class from $H^*(W)$
can be represented by a form $t$ which splits into the sum of forms
$t_\pm$  with cylindrical ends on $W_\pm$,
    so that we have $t_\pm|_{V}=t_V$.
     Let us define now a series
      $\bF^\infty\in\hbar^{-1}\fD$
      by the formula
    \begin{equation}\label{eq:SFT-gluing}
 {\bF^\infty(q^-,p^+,t)}=  {\bF_-(q^-,p,t_-)}
 \lozenge
   {\bF_+(q,p^+,t_+)} ,
   \end{equation}
    where $p,q$      are variable associated to the space $H^*(\Pc)$
    of   periodic orbits of the Reeb vector field
    $R_\a$ of the contact form $\alpha$ on $V$.

     The following theorem is the main claim of this section, and
     similar to Theorems \ref{thm:HH} and \ref{thm:SFT-cobordism}
     and \ref{thm:SFT-chain} it  is a statement about the
     boundary of an appropriate moduli space of holomorphic curves.
     This time we deal with limits of $J^s$-holomorphic curves
     in $W$ when $s\to\infty$,
     i.e. when the family $J^s$ realizes the splitting
     of $W$ into the composition $W_-\circledcirc W_+$,
      see Theorem \ref{thm:comp2} above.

    \begin{theorem}\label{thm:SFT-composition}
The element $\bF^\infty$ is homotopic to the
   potential
   $\bF=\bF_{W,J ,\a^\pm} $
   for any generic  compatible almost complex structure $J$ on $W$.
   \end{theorem}
  \bigskip
  Let us now describe the above results on the level of rational potentials.
Let  $W_-=\ora{V^-V }$, $W_+=\ora{VV^+}$ and
 $W=W_-\circledcirc W_+=\ora{V^-V^+}$ be as above. Set
$\bh^\pm=\bh_{V^\pm},\bh=\bh_V$,
$\widehat\bh_-=\bh^- -\bh$, $\widehat\bh_+=\bh  -\bh^+$,
$\widehat\bh=\bh^+-\bh^-$,
$\fL_\pm=\fL_{W_{\pm}}$ and $ \fL=\fL_W$.

 The operation
$\lozenge:\hbar^{-1}\fD_-\times\hbar^{-1}\fD_+\to\hbar^{-1}\fD$ defined above reduces on the
rational level to the
   operation $$\sharp:\fL_-\times\fL_+\to\fL,$$ defined as follows.
   For $ f_\pm\in\fL_\pm$ we set
\begin{equation}
f_-\sharp f_+(q^-,p^+)=\big(f_-(q^-,p)+f_+(q,p^+)
-\sum_{\g\in\Pc}\kappa_\g^{-1}z^{-d_{\g}}q_\g p_\g\big)\big|_{L},
\end{equation}
where
\begin{equation*}
L=
\begin{cases}
q_\g=&\kappa_\g z^{d_\g}\frac{\partial f_-}{\partial p_\g};\\
p_\g=&\kappa_\g z^{d_\g}\frac{\partial f_+}{\partial q_\g}.\\
\end{cases}
\end{equation*}
Notice that given series $$F_-=\hbar^{-1}\sum\limits_0^\infty
(F_-)_g\hbar^g\in\hbar^{-1}\fD^-\quad\hbox{ and}\quad
 F_+=\hbar^{-1}\sum\limits_0^\infty (F_+)_g\hbar^g\in\hbar^{-1}\fD^+$$  with
$F_-\lozenge F_+= \hbar^{-1}\sum\limits_0^\infty (F_-\lozenge
F_+)_g\hbar^g\in\hbar^{-1}\fD $ then  $$(F_-\lozenge F_+)_0=(F_-)_0\sharp
(F_+)_0.$$
 We will also consider the operations
\begin{equation*}
\begin{split}
&\sharp f_+:\fL_-\to\fL,\;\; \sharp f_+(f_-)=f_-\sharp
f_+\;\;\hbox{and}\;\;
\sharp f_-:\fL_+\to\fL,\;\; \sharp f_-(f_+)=f_-\sharp
f_+,\\
\end{split}
\end{equation*}
and their linearizations
\begin{equation*}
\begin{split}
&T_{f_-}(\sharp f_+):\fL_-\to\fL,\;\; T_{f_-}(\sharp f_+)(g)=
(g|_{L_+})\sharp f_+
\quad\hbox{and}\\
&T_{f_+}(\sharp f_-):\fL_+\to\fL,\;\; T_{f_+}(\sharp f_-)(g)=f_-\sharp
(g|_{L_-})\,.\\
\end{split}
\end{equation*}
Here we view $g|_{L_+}$ (resp. $g|_{L_-}$) as an element of
$\fL_+$  which depends on variables $q^-$ as parameters (resp. an element of
$\fL_-$, which depends on $p^+$ as parameters).
  We have the following rational version
of Theorem \ref{prop:SFT-composition}.
\begin{proposition}\label{prop:RSFT-composition}
Suppose that even elements $f_\pm\in\fL_\pm=\fL_{W_{\pm}}$ satisfy
equation (\ref{eq:HJ1}), i.e.
\begin{equation*}
\widehat\bh_\pm|_{L_{f_\pm}}= 0 ,
\end{equation*}
where
\begin{equation*}
L_{f_-}=
\begin{cases}
q_\g&=k_\g\frac{\partial f_-}{\partial p_\g};\;\;\g\in\Pc\\
p^-_{\g^-}& =\kappa_{\g^-}\frac{\partial f_-}{\partial q^-_{\g^-}},\;\;\g^-\in\Pc^-,\\
\end{cases}
\end{equation*}
and
\begin{equation*}
L_{f_+}=
\begin{cases}
p_\g&=k_\g\frac{\partial f_+}{\partial q_\g};\;\;\g\in\Pc\\
q^+_{\g^+}&=\kappa_{\g^+}\frac{\partial f_+}{\partial p^+_{\g^+}},\;\;\g^+\in\Pc^+.\\
\end{cases}
\end{equation*}
Then
\begin{description}
\item{1.} The function $f_-\sharp f_+$ satisfies the Hamilton-Jacobi
equation
$$\widehat\bh_{L|_{f_-\sharp f_+}}=0;$$
\item{2.} The homomorphisms $T_{f_-}(\sharp f_+):\fL_-\to\fL$ and
$T_{f_+}( f_-\sharp):\fL_+\to\fL$ satisfy the equations
\begin{equation*}
\begin{split}
T_{f_-}(\sharp f_+)\circ d^{f_-}=&d^{f_-\sharp f_+}\circ T_{f_-}(\sharp
f_+),\\
T_{f_+}( f_-\sharp)\circ d^{f_+}=&d^{f_-\sharp f_+}\circ T_{f_+}(\sharp
f_-)\,, \\
\end{split}
\end{equation*}
and hence define homomorphisms of the corresponding
homology algebras:
\begin{equation*}
\begin{split}
\big(T_{f_-}(\sharp f_+)\big)_*:&H_*(\fL_-,d^{f_-})\to H_*(\fL,d^{f_-\sharp f_+}),\\
\big(T_{f_+}(f_-\sharp )\big)_*:&H_*(\fL_-,d^{f_-})\to H_*(\fL,d^{f_-\sharp f_+});\\
\end{split}
\end{equation*}
\item{3.}
\begin{equation*}
\begin{split}
T_{f_-}(\sharp f_+)\circ (f_-)^-=& (
f_-\sharp f_+)^-\,, \\
T_{f_+}(\sharp f_-)\circ (f_+)^-=& (
f_-\sharp f_+)^+
\,; \\
\end{split}
\end{equation*}
\item{4.}
Suppose we are given three cobordisms $W_1,W_2,W_3$ with matching ends
 so that we can  form  the  composition
$W_{123}=W_1\circledcirc W_2 \circledcirc W_3$, and series
$f_i\in\hbar^{-1}\fL_i= \fL_{W_i}$  which satisfy
Hamilton-Jacobi equations
$\widehat{\bh}_{W_i}|_{L_{f_i}}=0,\,i=1,2,3$.
Then
\begin{equation*}
T_{f_2}(\sharp f_3)\circ T_{f_1}(\sharp
f_2)=T_{f_1}(\sharp (f_2\sharp f_3)).
\end{equation*}
\end{description}
\end{proposition}

Let $t$ be a closed form on $W$ which is split into two forms $t_\pm$ on $W_\pm$ with
cylindrical ends. Set $\bff_\pm=\bff_{W_\pm}$ and
\begin{equation}
\bff^\infty(q^-,p^+,t)=\bff_-(q^-,p,t_-)\sharp \bff_+(q,p^+,t_+).
\end{equation}
Alternatively $\bff^\infty$ can be defined as
the first term  in the expansion
$\bF^\infty=\hbar^{-1}\sum\limits_0^\infty\bF^\infty_g\hbar^g$.
 The following theorem is a rational analog, and a direct
 corollary of Theorem \ref{thm:SFT-composition}.
\begin{theorem}
\label{thm:RSFT-composition}
  The series $\bff^\infty(q^-,p^+,t)$  and  $\bff_{W,J}(q,p,t)$
  are  homotopic  for any generic  compatible almost complex structure
  $J$ on $W$.
\end{theorem}

Coming down to the   ``classical" level, let us assume that
 $W,W_-$ and $W_+$ are homology cobordisms (see \ref{sec:3algebras}
 above). Thus there are   defined  the homomorphisms  $\bPsi  :\fA^+\to\fA^-$,
 $\bPsi_+ :\fA^+\to\fA $ and $\bPsi_- :\fA \to\fA^-$,
    see Section \ref{sec:relative} above. Set
 $\bPsi_\infty=\bPsi_+\circ\bPsi_-$.
 Then we have
 \begin{theorem}\label{thm:classical-composition}
 For any generic compatible almost complex structure $J$ on $W$
 homomorphisms
 $\bPsi_1=\bPsi_J,\bPsi_\infty:\fA^+\to\fA^-$ are chain
 homotopic.
\end{theorem}

\subsection{Invariants of contact
manifolds}\label{sec:contact-invariants}

 Theorem \ref{thm:SFT-composition}  allows us to define
 SFT-invariants of contact manifolds.
 Let $(V,\xi)$ be a contact manifold,  and
 $\alpha^+$ and $\alpha^-$
 two contact forms for $\xi$,
 such that $\alpha^+>\alpha^-$, i.e.
 $\alpha^+=f\alpha^-$, for a  function $f>1$.
   Then for an appropriately chosen function $\zeta:V\times\R\to(0,\infty)$
   the form
 $\omega=d(\zeta\alpha^-)$ on $W=V\times\R$ is symplectic, and
 $(W,\omega)$ is a directed symplectic cobordism between $(V,\alpha^-)$
 and $(V,\alpha^+)$.  Let  $t^\pm$ be two cohomologous forms on $V$, and $t$ be a  closed
 form  on $W$ with cylindrical ends which restricts to $t^\pm$ on $V^\pm$.
 Suppose we are also given almost complex structures $J^\pm$ on $V$,
  compatible with $\alpha^\pm$, which
 are extended to  a compatible almost complex structure $J$ on $W$.
 We will call a  directed symplectic cobordism  $(W,J,t)$,
  chosen in the  above way, a {\it concordance} between
  the data on its boundary. Notice, that concordance becomes   an
  equivalence relation if we identify contact forms proportional
  with a constant factor. A concordance $(W,J,t)$ is called {\it
  trivial} if  $W=V\times\R$, the almost complex structure
   $J$ is translationally  invariant, and
  $t$ is the pull-back of a form $t_+$ under the projection $W\to
  V$.

 Let us denote, as usual,   by $(\fW^\pm,D^\pm)$  the differential Weyl algebras
  associated to the data at the
 ends of the cobordism $W$, by $(\fD,D_W)$ the $D$-module
 $\fD(W,J,\a^\pm)$,  by $\bF\in\fW$
 the    potential of the cobordism $W$, and by
 $\bF^\pm:(\fW^\pm,D^\pm)\to(\fD,D_W)$ the
 corresponding homomorphisms
 of differential algebras defined in (\ref{eq:end-homomorphisms}).

 \begin{theorem}\label{thm:SFT-invariance}
  For any concordance $(W,J,t)$ the  homomorphisms $$\bF^\pm:(\fW^\pm,d^\pm)\to(\fD,D_W)$$
    are quasi-isomorphisms
  of differential algebras. In particular, the homology algebras
  $H_*(\fW^-,D^-)$ and
  $H_*(\fW^+,D^+)$ are isomorphic.
  \end{theorem}

 \noindent\proof We will prove   \ref{thm:SFT-invariance}  in three
  steps.

\noindent Step 1. Let us begin with the trivial concordance $(W ,J
,t  )$. In this case $\fD$ can be identified with $\fW^\pm$ and we
have $\bF(q,p,t)=\hbar^{-1}\sum\kappa_\g^{-1}q_\g p_\g$. Hence
\begin{equation*}
\bF^-(f)=e^{-\bF}\ora{f}e^{\bF}=
f=e^{-\bF}\left(
e^{\bF}\ola{f}\right)=\bF^+(f).
\end{equation*}

         \noindent Step 2. If we add now to $t$ a form $d\theta$, where
$\theta$ has a compact support, and change $J$ in a compact part
of $W$ then according to Theorem \ref{thm:SFT-chain}
 the potential $\bF_{W,J}(t+d\theta)$ remains the same
 up to homotopy, and hence Corollary \ref{cor:SFT-chain} implies
that the homomorphisms  $\bF^\pm_*$   induced on homology remain
 unchanged.

\noindent Step 3. Now assume that $(W,J,t)=(W^1,J^1,t^1)$ is a general
concordance. Consider the reversed concordance $(W^2,J^2,t^2)$, so
that the compositions $$(W^{12}=W^1\circledcirc W^2,J^{12}=J^1\circledcirc
J^2,t^{12}=
t_1\circledcirc t_2)$$ and
  $$(W^{21}=W^2\circledcirc W^1,J^{21}=J^2\circledcirc J^1,t^{21}=
t_2\circledcirc t_1)$$ of concordances $(W^1 ,J^1 ,t^1 )$
and $(W^2 ,J^2 ,t^2 )$   are   as in Step 2.
 Then according to Theorem
\ref{thm:SFT-composition}
$\bF_{W^{12},J^{12}}( t^{12})$ is homotopic to
$$\bF_{W^1,J^1}(t^1)\lozenge\bF_{W^2,J^2}(t^2)=
\bF^1\lozenge\bF^2$$ and
 $\bF^{21}=\bF_{W^{21},J^{21}}( t^{21})$ is homotopic to
$$\bF_{W^2,J^2}(t^2)\lozenge\bF_{W^1,J^1}(t^1)=\bF^2\lozenge\bF^1.$$ Hence
Proposition \ref{prop:SFT-composition}
implies
\begin{equation*}
\Id=\left(\bF_{W^1\circledcirc W^2 }\right)^-_*
=\left(T_{\bF^1}(\lozenge\bF^2) \right)_*\circ (\bF^1)^-_*
  \end{equation*}
  and
  \begin{equation*}
  \begin{split}
& \Id= \left(T_{\bF^1}\left(\lozenge(\bF^2\lozenge\bF^1) \right)
\right)_* \\
&= \left(T_{\bF^1\lozenge\bF^2}(\lozenge \bF^1)\right)_*
\circ \left(T_{\bF^1} (\lozenge  \bF^2)\right)_* .\\
\end{split}
\end{equation*}
Hence
$\left(T_{\bF^1}(\lozenge\bF^2) \right)_*$,  $\;(\bF^1)^-_*$,
and similarly    $(\bF_1)^+_* $
are isomorphisms. \qed

\bigskip

The following rational and classical analogs of Theorem
\ref{thm:SFT-invariance}
can be either deduced directly from Theorem
\ref{thm:SFT-invariance}, or alternatively can be proven similarly
to \ref{thm:SFT-invariance}  using \ref{thm:RSFT-chain} (resp. \ref{thm:classical-chain}) and
\ref{thm:RSFT-composition} (resp. \ref{thm:classical-composition}).

\begin{theorem}\label{thm:RSFT-invariance}
    For any concordance $(W,J,t)$ the  homomorphisms $$\bff^\pm:(\fP^\pm,D^\pm)\to(\fL,D_W)$$
     are quasi-isomorphisms
  of differential algebras. In particular,  Poisson homology
   algebras $H_*(\fP^\pm,d^\pm)$  are isomorphic.
  \end{theorem}

\begin{theorem}\label{cor:classical-invariance}
 For any concordance $(W,J,t)$ the  homomorphism  $$\bPsi:(\fA^+,\partial^+)\to
 (\fA^-,\partial^-)$$
  is a quasi-isomorphism
  of differential algebras.
\end{theorem}
\bigskip

The definition of the differential algebras $(\fW,D),\;(\fP,d)$ and $(\fA,\partial)$
 depends on the choice
of a coherent orientation (see Section \ref{sec:orientation}), and
  of spanning surfaces and framings
of periodic orbits
(see Section \ref{sec:dynamics}).
As it is stated in Theorem \ref{thm:coherent}  a coherent orientation
is determined by a choice of asymptotic operators associated
with each periodic orbit $\g$. Let $\bH'$ be the new  Hamiltonian which
one gets by changing the orientation of the asymptotic operator
associated with a fixed  periodic orbit $\g $. One can then check that the
change of variables  $(p_\g, q_\g)\mapsto(-p_\g,-q_\g)$ is an
isomorphism  between the differential algebras $(\fW,D^{\bH})$ and
$(\fW,D^{\bH'})$.
  Different choices for spanning   surfaces and framings of periodic orbits do not affect
$\mod\;2$ grading but change the integer grading of the
differential algebras.

\begin{remark}\label{rem:other-choices}
{\rm
  An accurate introduction of  virtual cycle techniques
would reveal that even more choices should be made.   However,  an independence
of  all these extra choices can be also established following the scheme of this section.
}
\end{remark}

\subsection{A differential equation for potentials of symplectic cobordisms} \label{sec:master}

In this section we describe  differential equations for
the  potentials $\bF_W$ and $\bff_W$ of a  directed symplectic
cobordism with a non-empty boundary. These equations completely
determine  the potentials, and
 in
combination with gluing Theorems  \ref{thm:SFT-composition}  and
\ref{thm:RSFT-composition} they provide in many cases   an effective
 recursive
procedure for computing     potentials $\bF_W$ and $\bff_W$,
 and even  Gromov-Witten invariants of closed symplectic manifolds  $W$
(see some examples in Section \ref{sec:computing} below).

\bigskip
Let us assume for simplicity that $W$ has only a positive end
$E^+=V\times(0,\infty)$, and choose a basic system
$\Delta_1,\dots,\Delta_k$, $\Theta_1,\dots,\Theta_m$ of closed
forms so that the   following conditions are satisfied:
\begin{description}
\item{a)} $\Delta_1,\dots,\Delta_k$ form a basis of $H^*(W)$, and the restrictions
$\delta_i=\Delta_i|_V$, $i=1,\dots,l$ for $l\leq k$ form a basis
of    $\im\big(H^*(W)\to H^*(V)\big)$;
\item{b)}  $\Theta_1,\dots,\Theta_m$ are compactly
supported and    represent a basis of
$\Ker\big(H^*_{\comp}(W)\to H^*(W)\big)$,
\item{c)} there exist  forms $\theta_1,\dots,\theta_m$ on $V$ and
a compactly supported $1$-form $\rho$ on $(0,+\infty)$, such that
$\Theta_j=e_*\big(\rho\wedge \pi^*(\theta_j)\big),\;j=1,\dots, m$, where

$\pi$ is the projection $E=V\times(0,\infty)\to V$ and $e:E\hookrightarrow W$ is the inclusion.
In  other words, $\Theta_j$ equals $\rho\wedge \pi^*(\theta_j)$
viewed as a form on $W$.
\end{description}
\medskip

\begin{theorem}\label{thm:SFT-master}
Let $\bH=\bH_{V,\a,J}$ be the Hamiltonian associated with the
contcat manifold $V$.
Set $$\bH^j(t_1,\dots,t_l,q,p)=\Big(
\frac{\partial \bH}{\partial s_j}(\sum\limits_{i=1}^l
t_i\delta_i+s_j\theta_j,q,p)\Big)\Big|_{s_j=0}
,\;j=1,\dots m,$$
$$\bF^0(t_1,\dots,t_k, p)=\bF_{W,J}(\sum
 t_i\Delta_i,p),$$
and define
 $\bF(t_1,\dots,t_k,\tau_1,\dots,\tau_m, p)$
 by the formula
 \begin{equation}
 e^{\bF(t_1,\dots,t_k,\tau_1,\dots,\tau_m,
 p)}=e^{\bF^0(t_1,\dots,t_k,p)}
 \overleftarrow\bG(t_1,\dots,t_l,\tau_1,\dots,\tau_m,p),
 \end{equation}  where we denote by
 $\overleftarrow\bG$ the operator obtained from
 \begin{equation}
  \bG(t_1,\dots,t_l,\tau_1,\dots,\tau_m,q,p)
  = e^{\tau_m\bH^m(t_1,\dots,t_l,q,p)}\dots
 e^{\tau_1\bH^1(t_1,\dots,t_l,q,p)}
 \end{equation}
 by quantizing   $q_\g=\kappa_\g
 \hbar\overleftarrow{\frac{\partial }{\partial
 p_\g}}$.
 Then  $\bF(t_1,\dots,t_k,\tau_1,\dots,\tau_m, p)$  is
  homotopic to  $$\bF_{W,J}(\sum\limits_{i=1}^k
 t_i\Delta_i+\sum\limits_{r=1}^m\tau_r\Theta_r,p).$$
  \end{theorem}
  \proof
 Set  $$T^j(s)={\bF_{W,J}(\sum
 t_i\Delta_i+\sum\limits_{r=1}^{j-1}\tau_r\Theta_r+s\tau_j\Theta_j,p) }.$$
  We have $$T^j(1)=T^{j+1}(0)\quad\hbox{
 for}\;\;j=1,\dots,m-1 ,     $$
  $$
  T^m(1)= {\bF_{W,J}(\sum\limits_{i=1}^k
 t_i\Delta_i+\sum\limits_{r=1}^m\tau_r\Theta_r,p)},$$
 and
 $$T^1(0)={\bF_{W,J}(\sum
 t_i\Delta_i,p) }={\bF^0(t_1,\dots,t_k, p)}.$$
 Let $S^j\in\hbar^{-1}\fD$ be defined from the equation
 $$e^{S^j}=e^{T^j(0)}e^{\tau_j\ola{\bH^j(t_1,\dots,t_l,q,p)} }.$$
 It is enough to prove that
 $  T^j(1)$ is homotopic    to $S^j$ for $j=1,\dots,m$.
 We have
 \begin{equation*}
 \begin{split}
 &\frac{\partial e^{T^j(s)}}{\partial s}=
  \frac{\partial T^j(s)}{\partial s} e^{T^j(s)}  \\
  &=  e^{T^j(s)}\sum\limits_d\sum\limits_{g,u  ,v
  =0}^\infty     \frac{1}{u!v ! }
   \up0\langle \tau_j\Theta_j,
    \underbrace{\sum\limits_{i=1}^k
 t_i\Delta_i+\sum\limits_{r=1}^{j-1}\tau_r\Theta_r+s\tau_j\Theta_j}_u;
 \underbrace{p ,\dots,p }_{v }
 \rangle^d_{g}z^d \hbar^{g-1}.\\
 \end{split}
 \end{equation*}
  The compactly supported form $\Theta_j$ is exact in $W$,
  $$\Theta_j=d\wt\Theta_j,$$
where $\wt\Theta_j$ is closed at infinity, has a cylindrical    end,
and $\wt\Theta_j|_V=\theta_j$.
Hence,
 \begin{equation*}
 \frac{\partial T^j(s)}{\partial s}  =
\ \sum\limits_d\sum\limits_{g,u  ,v =0}^\infty
\frac{1}{u!v !}
  \up0\langle \tau_jd\wt \Theta_j,
    \underbrace{\sum\limits_{i=1}^k
 t_i\Delta_i+\sum\limits_{r=1}^{j-1}\tau_r\Theta_r+s\Theta_j}_u;
  \underbrace{p ,\dots,p }_{v }
 \rangle^d_{g}z^d \hbar^{g-1}
 \end{equation*}
  \begin{equation*}
 \begin{split}
 &= d\big( \sum\limits_d\sum\limits_{g,u  ,v=0}^\infty  \frac{1}{u!v
 !}
  \up0\langle \tau_j\wt \Theta_j,
    \underbrace{\sum\limits_{i=1}^k
 t_i\Delta_i+\sum\limits_{r=1}^{j-1}\tau_r\Theta_r+s\tau_j\Theta_j}_u;
  \underbrace{p ,\dots,p }_{v }
 \rangle^d_{g}z^d \hbar^{g-1}\big)\\
  &= \frac{\partial}{\partial u}
  \Big( d\big( \sum
  \limits_d\sum\limits_{g,u  ,v=0}  \frac{1}{u!v!}\\
 &\up0\langle
    \underbrace{\sum\limits_{i=1}^k
 t_i\Delta_i+\sum\limits_{r=1}^{j-1}\tau_r\Theta_r+s\tau_j\Theta_j
 +u\tau_j\wt \Theta_j}_u;
 \underbrace{p ,\dots,p }_{v}
 \rangle^d_{g}z^d \hbar^{g-1}\big)\Big)\Big|_{u=0}\\
& =  \frac{\partial}{\partial u}\Big( d\big(
\bF_{W,J}(\sum\limits_{i=1}^k
 t_i\Delta_i+\sum\limits_{r=1}^{j-1}\tau_r\Theta_r+
 s\tau_j\Theta_j+u\tau_j\wt \Theta_j,p)\big)\Big)\Big|_{u=0} ,
    \end{split}
 \end{equation*}
   where $d$ denotes the de Rham differential.
   Using  equation (\ref{eq:non-closed2}) we get
    \begin{equation*}
 \begin{split}
&  d\big(
\bF_{W,J}(\sum\limits_{i=1}^k
 t_i\Delta_i+\sum\limits_{r=1}^{j-1}\tau_r\Theta_r+
 s\tau_j\Theta_j+u\tau_j\wt \Theta_j,p)\big) =\\
 &e^{-\bF_{W,J}}\big( e^{\bF_{W,J}}\ola\bH(\sum\limits_{i=1}^l t_i\Delta_i +
 u\tau_j\theta_j,p) \big),\\
     \end{split}
 \end{equation*}
 and hence
 \begin{equation*}
 \begin{split}
   &\frac{\partial T^j(s)}{\partial s}=
   \frac{\partial}{\partial u}\Big( e^{-\bF_{W,J}}
   \big( e^{T^j(s)}\ola\bH(\sum\limits_{i=1}^l t_i\Delta_i +
 u\tau_j\theta_j,q,p) \big)\Big)\Big|_{u=0}\\
 &= \tau_j e^{-T^j(s)}\Big(-\bF^j  (t_1,\dots,t_k,\tau_1,\dots,\tau_j,s,p)
  \big(e^{T^j(s) }\ola\bH(\sum\limits_{i=1}^l t_i\Delta_i ,q,p)\big)\\
 &+ \big( e^{T^j(s)} \bF^j  (t_1,\dots,t_k,\tau_1,\dots,\tau_j,s,p)\big)
  \ola\bH(\sum\limits_{i=1}^l t_i\Delta_i ,q,p)\\
 & +e^{T^j(s)}
 \ola\bH^j(t_1,\dots,t_l,q,p)\Big) =
 \tau_j e^{-T^j(s)}\Big(\big( e^{T^j(s)} [\bF^j , \ola\bH]
 +e^{T^j(s)}
 \ola\bH^j\Big),  \\
    \end{split}
 \end{equation*}
    where
    \begin{equation*}
    \begin{split}
    &\bF^j  (t_1,\dots,t_k,\tau_1,\dots,\tau_j,s,p) \\
 &=\frac{ \partial   }{\partial u} \big(
 \bF(\sum
 t_i\Delta_i+\sum\limits_{r=1}^{j-1}\tau_r\Theta_r+
 s\tau_j\Theta_j+u\wt\Theta_j,p)\big)\big|_{u=0}
 \end{split}
 \end{equation*}
 and
 $$\bH^j( t_1,\dots,t_l,q,p)=   \frac
 { \partial   }{\partial u} \big(
 \bH(\sum t_i\Delta_i +
 u\theta_j,q,p)   \big)\big|_{u=0}\,.$$
 Therefore,
 \begin{equation*}
\frac{\partial e^{T^j(s)}}{\partial s}=  e^{T^j(s)}
  \frac{\partial T^j(s)}{\partial s}
  =\tau_je^{T^j(s)}\big([\bF^j , \ola\bH]+\ola\bH^j\big)\,.
 \end{equation*}

Let us define now a family
 $U^j(s)\in \hbar^{-1}\fD,s\in [0,1]$
 by the formula
 \begin{equation}
 e^{U^j(s)}=e^{T^j(s)}e^{(1-s)\tau_j\ola
 {\bH^j}(t_1,\dots,t_l,q,p)}\,.
     \end{equation}
     Then $U^j(s)$ is a homotopy
     between $S^j$ and $T^j(1)$. Indeed,   we have
 $U^j(0)=S^j$ and $U^j(1)=T^j(1)$.
 On the other hand we get an equation
 \begin{equation*}
\frac {\partial e^{U^j(s)}}{\partial s} =\tau_je^{U^j(s)}
 \big(-\ola\bH^j+ [\bF^j
, \ola\bH]+\ola\bH^j\big) =e^{U^j(s)}[\tau_j\bF^j , \ola\bH]\;,
\end{equation*}
which is the definition of homotopy (see Section
\ref{sec:chain} above).\qed

\bigskip

We formulate now a version of Theorem \ref{thm:SFT-master} for
rational potentials.

Set $$\bh^j(t_1,\dots,t_l,q,p)=\frac{\partial \bh}{\partial
s_j}(\sum t_i\delta_i+s_j\theta_j,q,p) ,\;j=1,\dots m,$$ and for
any $\bgg\in \fL$ we denote by $L_{\bgg}$ the Lagrangian variety of
$\bV$, defined by the equation
\begin{equation}
L_\bgg=\{q_\g= \kappa_\g\frac{\partial\bgg}
{\partial p_\g}\}.
\end{equation}
\begin{theorem}\label{thm:RSFT-master}
 Let $\bff(t_1,\dots,t_k,\tau_1,\dots,\tau_m, p)$
be the solution of the   Hamilton-Jacobi equation
\begin{equation} \label{eq:RSFT-master}
\frac{\partial \bff}{\partial \tau_j} (t_1,\dots,t_k,
\tau_1,\dots,\tau_m)=
\bh^j(t_0,\dots,t_l,q,p)|_{L_{\bff}}
\end{equation}
with the initial condition
 $$\bff|_{\tau_j=0}=\bff_{W,J}(\sum
 t_i\Delta_i+\sum\limits_{r\neq j}\tau_r\Theta_r,p).$$
Then $$\bff(t_1,\dots,t_k,\tau_1,\dots,\tau_m, p)$$ is
  homotopic to $$\bff_{W,J}(\sum
 t_i\Delta_i+\sum\limits_{r=1}^m\tau_r\Theta_r,p).$$
 \end{theorem}

\subsection{Invariants    of Legendrian
knots}\label{sec:Legendrian}
Symplectic Field Theory can be extended to include Gromov-Witten
invariants of pairs $(W,L)$, where $L$ is a Lagrangian
submanifold of a symplectic manifold $W$.   The corresponding relative object is a
pair $(W,L)$, where $W=\ora {V^-V^+}$ is a directed symplectic cobordism between
contact manifolds $(V^\pm,\a^\pm)$, and $L$ is a Lagrangian cobordism between Legendrian
submanifolds $\Lambda^\pm\subset V^\pm$.  More precisely, we assume that Lagrangian submanifold
$L$ is cylindrical at infinity over $\Lambda^\pm$, i.e. there exists $C>0$, such that
$\;L\cap V^-\times (-\infty, -C]=\Lambda^-\times (-\infty, -C]\;$ and
$\; L\cap V^+\times (C,\infty]=
\Lambda^+\times (C,\infty]\;$.
                                  In other words, we require $L$
                                   to coincide at infinity with  symplectizations
of Legendrian submanifolds $\Lambda^\pm$.

The moduli space of holomorphic curves to be considered in this case
consists of  equivalence classes of     holomorphic curves   with boundary which can have
punctures of two types,
 interior  and at the boundary.  The boundaries of holomorphic curves
are required to be mapped to the Lagrangian submanifold $L$,   the holomorphic  curves
should be cylindrical over periodic orbits from $\Pc^\pm$ at interior punctures,
 while at boundary punctures
we require them to be asymptotically cylindrical over Reeb chords connecting points  of the Legendrian submanifolds
$\Lambda^\pm\subset V^\pm$.  A  more precise definition is given below.
The algebraic structure arising from the stratification of boundaries of these moduli spaces is more complicated
than in the closed case. First of all, unlike the interior punctures the    punctures
 at the boundary are cyclically ordered, which leads to associative, rather than graded commutative algebras.
 Second, the ``usual" cusp degenerations of curves with boundary
 at boundary points (see \cite{Gromov-holomorphic}) has
 in this case codimension
 $1$, rather than $2$ as in the  closed case, and hence the  combinatorics of such
  degenerations  should also be a part of the
 algebraic formalism.

 We will sketch in this paper only the simplest of three cases of SFT, namely the ``classical case",
 which corresponds to the theory of moduli spaces of holomorphic
 disks with only $1$ positive puncture at the boundary.
 \bigskip

         Let $(V,\xi=\{\a=0\})$ be a contact  manifold with a fixed
         contact form $\a$, $W=V\times \R$ its symplectization
         with a compatible almost complex structure $J$,
         $\Lambda\subset V$ a  compact Legendrian submanifold,
         and $L=\Lambda\times\R\subset W$ the symplectization of $\Lambda$, i.e.
         the corresponding Lagrangian cylinder in $W$.   We  assume that all periodic orbits
         of  the Reeb   vector field $R_\a$ are non-degenerate
           and fix a marker
         on every periodic orbit. We also consider the set $\Cc$   of Reeb chords
          connecting points on $\Lambda$, and impose an extra
           non-degeneracy condition along the   chords from $\Cc$
          by requiring that the linearized flow of $R_\a$  along
                   a chord   $c\in\Cc$  connecting points $a,b\in\Lambda$
                     sends  the tangent  space $T_a(\Lambda)$ to
a space transversal to $T_b(\Lambda)$.
We also require that none of
the chords from $\Cc$ be a part of an orbit from
$\Pc$.
 Under these assumptions, the set $\Cc$ is finite: $\Cc=\{c_1,\dots,c_m\}$.

 We will restrict the consideration to the case when
 \begin{equation}
 \pi_1(V)=0,\; \pi_2(V,\Lambda)=0,\quad\hbox{ and}\quad w_2(\Lambda)=0.
 \end{equation}
 First two assumptions are    technical and can be removed (comp. Section \ref{sec:dynamics} above).
  However, the third one
 is essential for orientability of the involved moduli spaces of holomorphic curves. Moreover, the invariants
  we define depends on a particular choice of a spin-structure on $\Lambda$.
  \footnote{We thank K. Fukaya for pointing this out.}

 As in  Section \ref{sec:dynamics} we choose capping surfaces $F_\g$ for $\g\in\Pc$, and  for each chord
 $c\in \Cc$ we  also choose   a surface $G_c$
which is bounded by  a curve
  $c\cup\delta_c, $ where $\delta_c\subset \Lambda$.
   The choice of  surfaces $F_\g,\,\g\in\Pc,$  allows us to define
      Conley-Zehnder indices of  periodic orbits
   (see Section \ref{sec:dynamics} above). Similarly, surfaces $G_c$ enable us to define
   {\it Maslov indices} $\mu(c) ,\,c\in\Cc$.   Namely,  let us consider a Lagrangian subbundle of $\xi|_{\partial G_c}$,
   which consists of  the Lagrangian sub-bundle
  $T\Lambda|_{ \delta_c}\subset \xi|_{\delta_c}$ over $\delta_c$, together with
   the family  of Lagrangian  planes  $T_u\subset \xi_u,u\in c$, which are images of $T_a(\Lambda)$ under the linearized flow of
   the Reeb field $R_\a$.  Choose a symplectic trivialization of $\xi|_{ \partial G_c}$ which extends
   to $G_c$.   With respect to this trivialization the above sub-bundle can be viewed
   as a path of Lagrangian planes in a
   symplectic vector space.  The Maslov index of such path
    is defined as in \cite{Robbin-Salamon}.

    \medskip
   Consider a unit disk
   $D\subset\C$ with punctures $$\left(\{z^+,z_1^-,\dots,z_\sigma^-\}\cup\{x^-_1,\dots, x^-_s\}\right),$$
   where  $\bz=\{z^+,z_1^-,\dots,z_\sigma^-\}$, $0\leq\sigma\leq m$,
   is  a   counter-clockwise ordered set of punctures on $\partial D$,
   and  $\bx=\{x^-_1,\dots, x^-_s\}$ is an ordered set of interior punctures.
    As usual  we provide  interior punctures with asymptotic markers.

                        Let us denote by
                         $\Mc^A(\{c_{i_1},\dots,c_{i_\sigma}\},\{\g_1,\dots,\g_s\},c_i;W,\Lambda,J)$
                         the moduli space of
$J$-holomorphic maps
\begin{equation*}
(D\setminus \left(\bz\cup\bx\right),
\partial\left(D\setminus \left(\bz\cup\bx\right)\right) \to (W,L),
\end{equation*}
 which
are asymptotically cylindrical at    the negative end over the periodic orbit $\g^-_{k}$
  at the puncture $x^-_{i_k}$, and over the chord $c_{i_k}$ at the puncture $z_k$,
asymptotically cylindrical at    the positive  end over the chord $c_i$
   at the puncture $z^+$, and which send
asymptotic markers of interior punctures to
the  markers on the  corresponding
periodic orbits.

  Two maps are called equivalent if they differ by a conformal map
 $D\to D$ which preserves all punctures, marked points and     asymptotic markers.
    Each moduli space  $\Mc^A(c_i,\{c_{i_1},\dots,c_{i_\sigma}\},\{\g_1,\dots,\g_s\};W,\Lambda,J)$ is invariant under
    translations  $V\times\R\to V\times\R$ along the factor  $\R$, and we
   denote the corresponding quotient moduli space  by
     $$\Mc^A(c_i,\{c_{i_1},\dots,c_{i_\sigma}\},\{\g_1,\dots,\g_s\};W,\Lambda,J)/\R\;.$$

Let $(\fA,\partial)=(\fA(V,\alpha),\partial_J)$ be the graded commutative differential algebra
defined in Section
\ref{sec:3algebras} above, or rather its specialization at the point $0$.
 Consider a graded  associative  algebra $\fK=\fK(V,\Lambda,\alpha)$ generated
by elements $c_i\in \Cc$ with coefficients in the algebra $\fA$.
We define a differential $\partial_{\Lambda}=\partial_{\Lambda,J}$
on $\fK$ first on the generators $c_i$ by the formula
\begin{equation}\label{eq:d-leg}
 \partial_{\Lambda}(c_i)=
 \sum \frac{n_{\Gamma,I,d}}
 {k!\prod_1^k \kappa_{\g_j}^{i_j}i_j!}
  c_{j_1}\dots c_{j_\sigma}
 q_{\g_1}^{i_1}\dots q_{\g_k}^{i_k}z^d,
  \end{equation}
   where the sum is taken over all $d\in H_2(V)$, all ordered sets of
   different periodic orbits $\Gamma=\{\g_1,\dots,\g_k\} $,
    all multi-indices  $J=(j_1,\dots,j_\sigma),\,1\leq j_i\leq m$,
and $I=(i_1,\dots,i_k)$, $i_j\geq 0, $ and where the coefficient
$n_{\Gamma,I,d}$  counts the algebraic number of elements of the
moduli space $$\Mc^d\big(c_i,\{c_{j_1},\dots,c_{j_\sigma}\},
\{\underbrace{\g_1,\dots,\g_1}_{i_1},\dots,
\underbrace{\g_k,\dots,\g_k}_{i_k}\}\big)/\R,$$ if this space is
$0$-dimensional, and equals $0$ otherwise. The differential
extends to the whole  algebra $\fK$ by the graded Leibnitz rule.
However, it does not treat coefficients as constants: we have
$\partial_\Lambda(q_\g)=\partial(q_\gamma)$, where  $\partial$ is
the differential defined on the algebra $\fA$.

Then we have
\begin{proposition}
$$  \partial_{\Lambda}^2=0.$$
\end{proposition}
Given a family   of contact forms $\Lambda_\tau,
\alpha_\tau,
J_\tau\;\tau\in[0,1]$
of Legendrian submanifolds, contact forms, and compatible almost complex structures one can define, similar to the  case of closed
contact manifolds (see Sections \ref{sec:relative}  and Section \ref{sec:chain} above)
a   homomorphism of differential algebras
$$\Psi_S:\fK(V, \Lambda_0,\alpha_0)\to\fK(V, \Lambda_1,\alpha_1),$$
which is independent up to homotopy of the choice of a connecting
homotopy.
Composition of homotopies generates composition of homomorphisms, and hence one conclude
\begin{proposition}
The quasi-isomorphism type of the differential algebra
  $$(\fK(V,\Lambda,\alpha),\partial_{\Lambda,J})$$ depends only on the contact structure $\xi$
  and the Legendrian isotopy class of $\Lambda$. The
  {\em Legendrian contact homology} algebra
  $H_*(\fK,\partial_\Lambda)$ has a structure of a module over
  the contact homology algebra
  $H^\cont_*(V,\xi)=H_*(\fA,d)$, and it is an invariant of the Legendrian knot (or link) $\Lambda$.
  \end{proposition}

    The theory looks especially simple when the contact structure $\xi$ on
    $V$ admits a contact form $\alpha$ such that the Reeb vector
    field $R_\a$ has no closed periodic orbits. If, in addition
    the space of trajectories is a  manifold $M$ (e.g. when
    $V=J^1(N)=T^*(N)\times\R$ with a contact form $dz+pdq$), then
     $W$ is automatically symplectic, and the projection $\pi:W\to V$ sends the Legendrian submanifold
    $\Lambda\subset V$ to an immersed {\it Lagrangian} submanifold
    $L\looparrowright M$ with transverse self-intersection points. These points
    correspond to Reeb chords  $c_i$
    connecting points on $\Lambda$.
    Hence,  the algebra $\fK$ in this case is just a free associative
    algebra, generated over $\C$ (or $\C$) by the
    self-intersection points of $L$. It is possible
     to choose a compatible  almost complex structures
    $J$ on  the symplectization  $W=V\times\R$ and  $J_M$ on $M$
    to make
    the projection $W\to M$ holomorphic (comp. Section \ref{sec:Bott} below). Then  punctured
    holomorphic disks in $W$  from moduli spaces
    $\Mc^A(c_i,\{c_{i_1},\dots,c_{i_\sigma}\} ;W,\Lambda,J)$
    project to $J_M$-holomorphic disks in $M$ with boundary in the
    immersed Lagrangian  manifold $L$. Conversely, one can
    check that each such disk  lifts to a disk from  the
   corresponding moduli space $\Mc^A(c_i,\{c_{i_1},\dots,c_{i_\sigma}\}
   ;W,\Lambda,J)$, uniquely, up to translation along
    the $\R$-factor in $W=V\times\R$. This is especially useful when $\dim M=2$.  In
   this case $L$ is an immersed curve, and the holomorphic disks
   are precisely immersed, or branched disks with their boundaries in $L$.
   Moreover, branched disks are never rigid, because the branching
   point may vary. Hence, the  differential
   $\partial:\fK\to\fK$ can be defined in this case in a pure
   combinatorial way, just   summing the terms  corresponding
   to  all  appropriate immersed disks whose boundary consists of arcs of $L$, and
    which are locally convex near their corner.

    Yu. Chekanov independently realized  (see \cite{Chekanov}) this program  for
    Legendrian links in the standard contact $\R^3$.
    He was also  motivated by the  hypothetical description of the
    compactification of the  moduli spaces of holomorphic discs,
    but has chosen to prove the  invariance of the quasi-isomorphism type
    of the differential algebra $(\fK,\partial)$ in a pure
    combinatorial way.  In fact, he proved   a potentially stronger    form  of equivalence
    of
    differential algebras of isotopic Legendrian links,
    which
    he called     stable tamed isomorphism.  Stable
    tame
    isomorphism implies quasi-isomorphism, but we do not know
    whether it is indeed stronger. Let also note that
     Chekanov considered a   $\Z_2$-version of the theory.
    In some examples it works better
      the  $\Q$-version, which is provided by our formalism.
    J. Etnyre--J. Sabloff (\cite{Etnyre-Sabloff}) and L. Ng (\cite{Ng}) worked out
     the combinatorial
    meaning of signs dictated by the coherent orientation
    theory (see Section \ref{sec:orientation} above), and proved
    invariance of the
    stable tame type of the differential algebra $(\fK,\partial)$
    {\it defined over $ \Z$}. Note that Chekanov's paper \cite{Chekanov}
    also contains   examples  which show that the  stable  tame
    $\Z_2$-isomorphism type do distinguish some Legendrian knots,
    which could not be  previously distinguished.

 Similar to the absolute case of SFT, one can define   further invariants of Legendrian
  submanifolds by including in the formalism higher-dimensional moduli spaces. For instance,  by introducing
    marked points on the boundary
  of the disk
  one gets a non-commutative deformation of Legendrian  contact homology
  along the homology of Legendrian
  manifolds. This is useful, in particular, to define invariants of
    Legendrian links with ordered components.
  However, the full-scale generalization of Symplectic Field Theory  to
  directed symplectic-Lagrangian cobordisms between pairs
  of contact manifolds and their Legendrian submanifolds, which would formalize
  information  about moduli spaces of   holomorphic curves of arbitrary genus and arbitrary number
  of positive and negative punctures, is not straightforward due to existence of different type
  of codimension $1$ components on the boundary of the corresponding moduli spaces.
  We will discuss this theory
  in one of our future papers.

\subsection{Remarks, examples, and further algebraic constructions in
SFT}\label{sec:remarks}
 \subsubsection{Dealing with torsion elements in $H_1$} \label{sec:torsion}
Let us discuss grading issues for a contact manifold
$(V,\xi=\{\a=0\})$ in the case when
  the torsion part  of $H_1(V)$ is
non-trivial. As we will see it is impossible to assign in a coherent way
an integer grading to torsion elements and still keep the property
that the Hamiltonian $\bH$ has total grading $-1$.  We will
deal with this problem by assigning to some elements fractional
degrees, and thus obtain a rational grading, incompatible with the
canonical $\Z_2$-grading. In fact the
term ``grading" is  misleading in this case, and more appropriately one
should talk about
     an Euler vector field with rational coefficients.

Let us split  $H_1(V)$ as $T\oplus F$, where $T$  and $F$ are   the
torsion and free parts, respectively. As in Section \ref{sec:dynamics}
above let us  choose curves ${C}_1,\dots,{C}_k$ representing a basis
of $F$, fix a trivialization of the bundles $\xi_{{C}_i}$,
for any periodic orbit $\g\in\Pc_\a$ with $[\g]\in F$ choose a
surface $F_\g$ which realizes the homology between $[\g]$ and a
linear combinations $\sum n_i [{C}_i] $, and trivialize the
bundle $\xi|_\g$ accordingly.
For any other  periodic orbit $\g$ let $\g_l$ be its smallest multiple which
belong to  $F$. In particular, the bundle $\xi_{\g_l}$ is
already trivialized by a framing $f$. The problem is that in general there is no
framing over $\g$ which would produce $f$ over $\g_l$. Choose then
an arbitrary framing $g$  over $\g$ and denote by $g_l$ the
resulting framing over $\g_l$.  Let
$2m(g_l,f)\in\pi_1(Sp(2n-2,\R))=\Z$ be the Maslov class of the
framing $g_l$ with respect to $f$. The Conley-Zehnder indices of
$\g_l$ with respect to these two gradings are then
 related  by the formula
\begin{equation*}
\CZ(\g|f)=\CZ(\g|g_l)+2m(g_l,f).
\end{equation*}
We then assign to $\g$ the fractional  degree
\begin{equation}
\deg \g=\CZ(\g|g)-\frac{2m(g_l,f)}{l}\,.
\end{equation}
With this modification SFT can be extended to the case of contact
manifolds with no restrictions on $H_1$. However, the price we pay
is that this grading, even if integer, may not be compatible with
the universal $\Z_2$-grading which determines the sign rules.

\subsubsection{Morse-Bott formalism}\label{sec:Bott}

 Our assumption that all periodic orbits from $\Pc_\a$ for the
 considered contact forms $\a$
 are non-degenerate, though generic,  but is very inconvenient for
 computations: in many interesting examples periodic orbits
 come in continuous families. Sometimes  the  Reeb
flow is periodic, and it sounds quite stupid to destroy this beautiful
symmetry.

In fact the above formalism can be adapted to this ``Morse-Bott"
case. We sketch below how this could be done   for the periodic
Reeb flow of an $S^1$-invariant form of a
 pre-quantization space.
We consider below only the ``semi-classical" case which concerns
 moduli spaces of rational holomorphic curves.

Let $(M,\omega)$ be a symplectic manifold  of dimension $2n-2$
with  an integral cohomology class $[\omega]\in H^2(M)$. We will assume for simplicity
that $H_1(M)=0$.
 The
pre-quantization space $V$ is a circle bundle over $M$ with
first Chern class equal to $[\omega]$. The fibration $\pi:V\to M$
admits a   $S^1$-connection form $\alpha$ whose curvature is
$\omega$. It   defines a $S^1$-invariant contact structure  $\xi$
on $V$, transversal to the fibers  of the fibration. The Reeb flow
of $R_\a$ is periodic, so all its trajectories
  are closed and coincide with the fibers of the fibration $\pi$, or their multiples.

  The fiber of the fibration $V$ is a torsion element
in $H_1(V)$, and if $l$ is the greatest divisor of the class $[\omega]$ then
the  $l$-multiple of the fiber is homological to $0$.

Consider the cylindrical cobordism  (the symplectization)
$W=V\times\R$ with an almost
 complex structure $J$ compatible with $\alpha$
and denote by $\Mc_{0,r}(s|W,J,\a)$ the moduli  space of rational
holomorphic curves in $W$ with $s$ punctures and $r$ marked points.
Near punctures the curves are required to be  asymptotically
cylindrical over some fibers of $V$, or their multiples.
However, we do not specify to which particular fiber they are
 being asymptotic, or whether this fiber is
 considered
on the positive, or negative end of $W$. We
 do not equip   curves from $\Mc_{0,r}(s|W,J,\a)$
with asymptotic markers of punctures, because   we  cannot
  fix  in a continuous way  points on each simple periodic orbit,
as we did in the non-degenerate case.

As it was already mentioned in Section \ref{sec:almost} above,
 $W$ can be viewed  as the total space of the complex  line bundle  $L$ associated with
the $S^1$-fibration $V\to M$, with the zero-section removed, and
the almost complex structure $J$ can be chosen
  compatible with the structure of this bundle, so that the projection   $W\to M$
becomes holomorphic with respect to a certain almost complex structure $J_M$ on $M$
compatible with $\omega$.
Then automatically the   bundle induced over  any holomorphic curve in the base
has a structure of a holomorphic line bundle.
With this choice of $J$ each holomorphic curve   $f\in \Mc_{0,r}(s|W,J,\a)$ projects
to a $J_M$-holomorphic sphere $\overline{f}:
\C P^1\to M$,   and  can be viewed as a meromorphic section of the   induced
bundle $(\overline{f})^*L$ over $\C P^1$.  This bundle is ample, and hence  poles of its
sections
correspond to   the  {\it negative} ends of $f$, while zeroes  correspond to  the positive ones.
Notice that although  the moduli spaces
 $\Mc_{0,r}(s|W,J,\a)$ can be identified with the moduli spaces
of closed  holomorphic curves in a $\C P^1$-bundle over $M$
with prescribed tangencies to two divisors,
 their compactifications are different, and in particular the
 compactification of the first moduli space may have codimension one strata
 on its boundary.

The correspondence $f\mapsto\overline{f}$ define a fibration
 $$\pr:\Mc_{0,r}(s|W,J,\a)/\R\to \Mc_{0,r+s}(M,J_M).$$
 The fiber $\pr^{-1}(\overline{f})$ is   the
  union of  (an infinite number of) disjoint circles, which are indexed by
 sequences of integers $(k_1,\dots,k_{s+r})$ with $\sum k_i=d_0
 =\int\limits_{A}\omega$,
 where $A\in H_2(M)$ is the homology class realized by $\overline{f}$,
  and where  there are precisely
 $s$  non-zero coefficients  $k_i$.

Let us  consider two copies $\Pc^\pm$ of  the space $\Pc=\Pc_\a$ of periodic orbits,
  as we need to differentiate
between positive and negative ends of holomorphic curves. We will write $\ddot{\Pc}=\Pc^+\cup\Pc^-$
 and define  the   evaluation maps:
\begin{equation}
ev_0:\Mc_{0,r}(s|W,J,\a)/\R\to V^{\times r}\;\; \hbox{ and}\;\;
 ev^+_-:\Mc_{0,r}(s|W,J,\a/\R)\to\ddot{\Pc}^{\times s}\,.
\end{equation}

Here $ev^+_-$ associates with each puncture the corresponding
point of $\ddot\Pc$.
The space $\Pc^\pm$  can be presented as
$    \coprod\limits_{k=1}^{\infty}\Pc^\pm_k, $
    where each $\Pc^\pm_k$ is a copy   of $M$, associated
    with $k$-multiple orbits.

We will denote forms on $\Pc^+$ by $p$, on $\Pc^-$ by $q$,
  denote by $p_k, q_k$ their restrictions to $\Pc_k^\pm$, and organize them
  into   Fourier series
  $u=\sum\limits_{k=1}^\infty(p_ke^{ikx}+q_ke^{-ikx})$. If we are given a
  basis of $H^*(M)$ represented by forms $\Delta_1,\dots,\Delta_a$
  we will consider only forms from the space generated by this
  basis, and write $p_k=\sum\limits_{i=1}^a p_{k,i}\Delta_i,
q_k=\sum\limits_{i=1}^a q_{k,i}\Delta_i$ and denote by
$u_i$ the  $\Delta_i$-component of $u$, i.e.
 $$u_i=\sum\limits_{k=1}^\infty(p_{k,i}e^{ikx}+q_{k,i}e^{-ikx})\;\;\hbox{
 and}\;\;
 u=\sum\limits_1^a u_i\Delta_i.$$

Given a  closed form $t$ on $V$ and  a class $A\in H_2(V)$
 we define the correlator
 \begin{equation}
 \begin{split}
 & \up{-1}\langle
\underbrace{t,\dots,t}_r; \underbrace{u,\dots,u}_s\rangle^A_{0} =\\
&  \int\limits_{\overline{\Mc^A_{0,r}(s|W,J,\a)/\R} } ev_0^*\big(t\otimes \dots
\otimes t\big)\wedge (ev^+_-)^*\big(u \otimes\dots
\otimes u\big) \bigg{|}_{x=0}.\\
\end{split}
\end{equation}

Let us  choose a basis $A_0,\dots, A_N$ in
  $H_2(M)$ in such a way that $\int\limits_{A_0}\omega=l>0$ and
  $\int\limits_{A_i}\omega= 0$ for $i=1,\dots, N$.
  Then the classes $A_i,i\geq 1 ,$ lift to classes $\wt{A}_i\in H_2(V)$  which
under the assumption $H_1(M)=0$ form
   a basis of $
  H_2(V)$. The   degree
  $d=(d_1,\dots,d_N )$ of a class
  $A\in H_2(V)$ is a vector of its coordinates in this basis.

To associate an absolute   homology
class with a holomorphic curve
we  pick   the $l$-multiple (recall that   $l$ denotes the greatest divisor of $\omega$)
of the fiber   $\g$ over a point $x\in M$
 and
choose a lift of the surface representing the class
 $A_0$ with $\int\limits_{A_0}\omega =l$ as a spanning surface
$F_\g$. Any other $m$-multiple of  $\g$ we will cap with the chain
 $\frac ml[F_\g]$. However, to fix a spanning surface
for a    fiber over any other point $y\in  M$  or
its multiples, one needs to make some extra choices, for instance   fix
 a path connecting $x$ and $y$. The condition $H_1(M)=0$ guarantees
  independence of the homology class of the resulting surface of
 the choice of this connecting path.  Notice that   with this choice,
 the degree of $f\in \Mc_{0,r}(s|W,J,\a)/\R$  equals $(d_1,\dots,d_N)$, if the degree of its projection
 $\pr(f)\in\Mc_{0,r+s}(M,J_m)$ is equal to $(d_0,d_1,\dots,d_N)$.

In this notation the rational Hamiltonian $\bh=\bh_{V,J,\a }$ is defined
 by the formula
\begin{equation}
\bh(t,u)=\sum\limits_{d}\sum\limits_{r,s=0}^\infty\frac{1}{r!s!}
    \up {-1} \langle \underbrace{t,\dots,t}_r;\underbrace{u,\dots,u  }_{s}
     \rangle_{0}^{d }  z^d.
     \end{equation}

   Suppose that a basis of $H^*(M)$, represented by
     closed forms   $\Delta_1,\dots,\Delta_a$, is  chosen in such
     a way that   for $b\leq
a$  the system forms $\wt{\Delta}_j=\pi^*(\Delta_j),\,j=1,\dots,b$, generate
 the image $\pi^*(H^*(M))\subset H^*(V)$, and the forms
 $\wt{\Theta}_1,\dots,\wt{\Theta}_c$,
  complete  it
   to a basis of $H^*(V)$.
    We will denote  (graded) coordinates in the space generated by the forms
     $\wt{\Delta}_j,\,j=1,\dots, b$ and
  $\wt{\Theta}_1,\dots,\wt{\Theta}_c$ by
    $ (t,\tau)=(t_1,\dots,t_b,\tau_1,\dots, \tau_c)$.

As usual, the Hamiltonian $\bh$ is viewed as an element  of  a  graded
 commutative Poisson
 algebra $\fP$,
 which consists of formal power series of coordinates of vectors $p_k$
 and  $T=(t,\tau)=
 (t_1,\dots,t_b,\tau_1,\dots,\tau_c)$
  with coefficients
 which are polynomials of coordinates of vectors $q_k=(q_{k,1},\dots,q_{k,a})$.
  The coefficients of these polynomials
 belong to a certain completion (see condition (\ref{eq:Novikov}) above) of the group algebra
 of $H_2(V)$.
 All  the variables $p_{k,i},q_{k,i}$  have in this case
 the same parity as forms $\Delta_i$, the parity of variables $t_i$ and $\tau_j$ is the same as the  degree of
 the corresponding
 forms $\wt \Delta_i$ and $\wt\Theta_j$.
  If $l=1$, i.e. when $H_1(V)=0$, then  the integer grading of variables which corresponds
  to the   choice of capping surfaces described above  is defined
   as follows:
\begin{equation}\label{eq:dim-Bott}
\begin{split}
\deg t_i&=\deg \wt{\Delta}_i-2;\\
\deg \tau_i&=\deg\wt{\Theta}_i-2;\\
\deg q_{k,i}&=\deg\Delta_i-2+2ck;\\
\deg p_{k,i}&=\deg\Delta_i-2-2ck;\\
\deg z_i&==2c_1(A_i)\\
\end{split}
\end{equation}
where $c=c_1(A_0)$.
As it was explained in Section \ref{sec:torsion}  if $l>1$ one can only define fractional degrees, given by the above
formulas (\ref{eq:dim-Bott}) with $c=\frac{c_1(A_0)}{l}$.

The following proposition is useful for applications (see below the discussion
of Biran-Cieliebak conjecture about subcritical symplectic mainfolds).
It follows from the fact that all the  moduli spaces
$\Mc_{g,r}(s|W,J,\a)$ which we defined above  are even-dimensional.

\begin{proposition}\label{prop:even-dim}
Let $(V,\xi)$ be the contact pre-quantization space for a
symplectic manifold $(M,\omega)$. Then all contact homology algebras
$$H_*^{\SFT}(V,\xi)\big|_{t=0},H_*^{\RSFT}(V,\xi)\big|_{t=0},H_*^{\cont}(V,\xi)\big|_{t=0}$$
specialized  at $0\in H^*(V)$
are
free graded, respectively  Weyl, Poisson, or commutative algebras,
generated by  elements
$$p_{k,i},q_{k,i},\;\;i=1,\dots,a,\;k=1,\dots\,.$$ In particular, the parts of all these homology
algebras graded by the homology class $w\in H_1(V)$ ({\rm see \ref{rem:grading} above})
are   non-trivial.
 \end{proposition}

The Poisson tensor on $\fP$  is determined in the ``$u$-notation" by
the  following generalization of the   formula
 (\ref{eq:u-notation}):
\begin{equation}\label{eq:u-vector}
 \frac{1}{2\pi i}\int\limits_0^{2\pi}\langle(\delta u)',\delta
v\rangle dx,
\end{equation}
 where $\langle\;,\;\rangle$ denotes Poincar\'e pairing
on cohomology $H^*(M)$,  which is given in the basis $\Delta_1,\dots, \Delta_a$ by the matrix
$$\eta_{ij}=\langle\Delta_i,\Delta_j\rangle=\int\limits_M\Delta_i\wedge\Delta_j.$$
The Poisson tensor  can be written in
$(p,q)$-coordinates as
\begin{equation*}
\sum\limits_{k=1}^\infty k\sum\limits_{i,j=1}^{a }\eta_{ij}\frac{\partial}{\partial
p_{k,i}}\wedge\frac{\partial}{\partial
q_{k,j}}\,.
\end{equation*}
It can be shown that the  above Hamiltonian $\bh$
satisfies the identity $\{\bh,\bh\}=0,$  and  that the differential
Poisson algebra $(\fP,d^{\bh})$ is quasi-isomorphic to the
corresponding  differential Poisson algebra defined  in Section
\ref{sec:3algebras} for any non-degenerate contact form for the
same contact structure  $\xi$ on $V$.
\bigskip

 The following formula (\ref{eq:Bourg}), which sometimes allows to compute
 the Hamiltonian $\bh$  of $V$
  in terms of the Gromov-Witten invariant $\bff=\bff_{M,J_M}$ of $M$, emerged   in a discussion  of
  the authors with T. Coates and  F. Bourgeois.

    \begin{proposition}\label{prop:Coates}
    Set
    $$
    \bh^{j}_{W,J}(t,q,p,z)=\frac{\partial\bh}
    {\partial\tau_j}(\sum\limits_1^b
    t_i\wt\Delta_i+\tau_j\wt\Theta_j,q,p,z)
    \big|_{\tau_j=0}\;,
    $$
    and
    $$\widehat\bff^j(t,z)=\frac{\partial\bff}{\partial s}\big(
    \sum\limits_1^a t_i\Delta_i +s\pi_*\wt\Theta_j,z\big)\big|_{s=0}
    \,,$$
    for $j=1,\dots,c$
          (comp. Theorem \ref{thm:RSFT-master}).
   Then   we have
   \begin{equation}
      \label{eq:Bourg}
      \begin{split}
    &\bh^j_{W,J}(t_1,\dots,t_b,q,p,z)\\&
    =\frac{1}{2\pi }\int\limits_0^{2\pi}
    \widehat\bff^j(t_1+u_1(x),\dots,t_b+u_b(x), u_{b+1}(x),\dots,u_a(x),\tilde z)dx
    \end{split}
        \end{equation}
     where $z=(z_1,\dots,z_N)$,
      $\tilde z=(e^{-ilx},z_1,\dots,z_N)$ and
      $l$ is the greatest divisor of $\omega$.
     \end{proposition}
     To prove    (\ref{eq:Bourg}) one just observes that
      the correlator
          \begin{equation*}
 \up{-1}\langle
\wt \Delta_{j_1},\dots,\wt\Delta_{j_r},\wt\Theta_j;
u_{i_1}\Delta_{i_i},\dots,u_{i_v}\Delta_v\rangle^d_{0}
\end{equation*}
equals  the Fourier coefficient with $e^{ilx}$ of
the       correlator
   \begin{equation*}
 \up{0}\langle
{\Delta_{j_1},\dots,\Delta_{j_r},\pi_*\wt\Theta_j}, {
u_{i_1}\Delta_{i_i},\dots,u_{i_v}\Delta_v}\rangle^{\tilde
d}_{0}\,.
\end{equation*}

      Notice that if $\wt\Theta_j$ is an odd form, then
      $$\bh(\sum\limits_1^b
    t_i\wt\Delta_i+\tau_j\wt\Theta_j,q,p,\tilde z)=
    \tau_j\bh^j(t,q,p,z),$$
    because all terms of $\bh$ must contain $\tau_j$.
                                        In particular,
     for $M=\C P^{n-1}$  the manifold $V$ is a rational homology sphere, and thus
     a volume form $\Theta$  on $S^{2n-1}$
     generates the odd part
     of $H^*(V;\R)$.
     Hence, the formula
     (\ref{eq:Bourg}) completely determines $\bh$.
         Namely, let
         $\bff(t,z)$ be the Gromov-Witten invariant of $\C P^{n-1}$,
         and let $\Delta_{2i},i=0,\dots,n-1$, be closed forms generating
     $H^*(\C P^{n-1})$, so that $\Delta_{2n-2}=\pi_*(\wt\Theta)$.
      Set $\wt\Delta_0=\pi^*(\Delta_0)$
      and $$\widehat\bff(t,z)=\widehat\bff^{2n-2}(t,z)=\frac{\partial\bff(t,z)}{\partial t_{2n-2}}.$$

     Then we have
     \begin{equation}\label{eq:Bourg2}
  \bh(t_0\wt\Delta_0+\tau\Theta,q,p) =
   \frac{\tau}{2\pi } \int\limits_0^{2\pi}
{\widehat\bff}( t_0\Delta_0+u, e^{-ix})  dx .
 \end{equation}

Let us consider some applications of the formula (\ref{eq:Bourg2}).


\bigskip

\noindent\textsf {Contact homology of the standard contact $3$-sphere}
\medskip

  Let $\pi:V=S^3\to M=\C P^1$ be the Hopf
 fibration. $V$ is the pre-quantization space for $(S^2,\omega)$ with
 $\int\limits_{S^2}\omega=1$. The $0$-form $\Delta_0=1$ and
   the  symplectic  $2$-form $\Delta_2=\omega$
   generate $H^*(M)$, the $0$-form
       $\wt{\Delta}_0=\pi^*(\Delta_0)$ and the volume form $\wt{\Theta}_3$
        with $\pi_*(\wt{\Theta}_3)=\Delta_2$ on $S^3$ generate $H^*(S^3)$.
         Thus we have  functional variables
                 \begin{equation*}
           u_j(w)=\sum\limits_{k=1}^\infty \left(p_{k,j}e^{ikx}+q_{k,j}e^{-ikx}\right),
 \end{equation*}
associated  to the classes $\Delta_j,\; j=0,2,$ and variables $t_0$ and $\tau$
 associated to $\wt{\Delta}_0$ and $\wt{\Theta}_3$. According to (\ref{eq:dim-Bott}) we have
\begin{equation*}
\begin{split}
&\deg q_{k,0}=-2+4k,\quad \deg q_{k,2}= 4k,
\quad\deg p_{k,0}=-2-4k,\\
&\quad\deg p_{k,2}=- 4k,\quad
 \deg t_0=-2,\quad\deg \tau=1.\\
\end{split}
\end{equation*}
 The  potential $\bff$ for $M=\C P^1$ can be written, as it well known (see
 also Section \ref{sec:computing} below),
as
\begin{equation}\label{eq:CP2}
\bff=\frac{t_0^2t_2}{2} +e^{t_2}z,
\end{equation}
  and hence
  \begin{equation*}
\widehat\bff=\frac{t_0^2}{2} +e^{t_2}z,
\end{equation*}
Thus applying (\ref{eq:Bourg}) we get the following expression for the rational Hamiltonian
$\bh$ for $S^3$:
\begin{equation}\label{eq:3-sphere}
\begin{split}
 \bh&=\frac{\tau}{2\pi i}\int\limits_0^{2\pi} \left(\frac{(t_0+u_0)^2}{2}+
e^{u_2-ix}\right)dx = \tau\Big(\frac{t_0^2}{2}+\sum\limits_{k\geq 1} q_{k,0}p_{k,0}\\
&+\sum\limits_{t,s\geq 0}\mathop{\mathop{\sum
\limits_{i_1,\dots,i_s\geq 0}}\limits_{j_1,\dots j_t\geq 0}}
\limits_{\sum\limits_1^s li_l-\sum\limits_1^tmj_m=1}\frac{q_{1,2}^{i_1}\dots
q_{s,2}^{i_s}p_{1,2}^{j_1}\dots p_{t,2}^{i_t}}{i_1!\dots i_s!j_1!\dots
j_t!}\Big)\,.\\
\end{split}
\end{equation}

Let us use (\ref{eq:3-sphere})
to compute the contact homology algebra $$H^{\cont}_*(S^3,\xi_0)=H_*(\fA(S^3,J,\alpha),\partial).$$
The part of $\bh$  linear in the $p$-variables has the form
\begin{equation*}
\tau\sum\limits_1^\infty\left( p_{k,0}q_{k,0}+ p_{k,2}h_k(q_{1,2},\dots,q_{k-1,2})\right)\;,
\end{equation*}
so that the differential $\partial :\fA\to\fA$
is given by the formulas
\begin{equation*}
\partial q_{k,2}=k\tau q_{k,0},\quad \partial q_{k,0}=k\tau h_k(q_{1,2},\dots,q_{k-1,2})\;.
\end{equation*}
Here are
few first polynomials $h_k$:
 \begin{equation*}
 \begin{split}
h_1=&1,\\
h_2=&q_{1,2},\\
h_3=&q_{2,2}+\frac12q_{1,2}^2,\\
h_4=&q_{3,2}+  q_{2,2}q_{1,2}+\frac16q_{1,2}^3.\\
\end{split}
\end{equation*}
Notice that $\im\partial$ coincides with the ideal $I(\tau)$
generated by $\tau$. Hence,  the homology algebra
$H_*(\fA,\partial)$  specialized over a point $t=(t_0,0)$ is a
free graded commutative algebra $ \fA_0$ generated by variables
$q_{k,0},q_{k,2},\;k=1,\dots,$  and in particular it has one
generator in each even dimension. On the other hand, over any
point $t=(t_0,\tau)$ with  $\tau\neq 0$ the algebra
$H_*(\fA,\partial)$ is isomorphic to a  proper subalgebra $\fA_1$
of $ \fA_0$. It       has its first non-trivial generator $g_1=
q_{1,2}-\frac12 q_{1,0}^2 $ in dimension $4$.
\bigskip

\begin{remark}\label{rem:lens-space}
{\rm The contact homology of the Lens space $V=L(l,1)$ which is
 the pre-quantization space for $(S^2,\omega)$ with
 $\int\limits_{S^2}\omega=l$ can be computed similarly.
 The variables   $p_{k,0}, q_{k,0}, p_{k,2}, q_{k,2}, t_0 $ and $\tau$ have the same
 $\Z_2$-grading, as in the case of $S_3$, i.e. all of them, except  $\tau$ are even.
 However, the grading assigned by the Euler field to $p_{k,0}, q_{k,0}, p_{k,2}, q_{k,2}$
  is fractional in this case
 and given by formulas
\begin{equation*}
\begin{split}
&\deg q_{k,0}=-2+\frac{4k}{l},\quad \deg q_{k,2}= \frac{4k}{l},
\quad\deg p_{k,0}=-2-\frac{4k}{l},\\
&\quad\deg q_{k,0}=- \frac{4k}{l} .\\
\end{split}
\end{equation*}

The formula for $\bh$ takes the form
\begin{equation}\label{eq:lens-space}
\bh=\frac{\tau}{2\pi i}\int\limits_0^{2\pi} \left(\frac{(t_0+u_0)^2}{2}+
e^{u_2-ilx}\right)dx
\end{equation}

We will not carry on here the computation of the contact homology  of
the Lens space $L(l,1)$, and only note, that  as in the case of
$S^3$
  the homology algebra $H_*(\fA,\partial)$  specialized over a point $t=0$ is a free
graded commutative algebra $ \fA_0$ generated by variables
$q_{k,0},q_{k,2},\;k=1,\dots,$
 and over any  point $t \neq 0$  the
algebra $H_*(\fA,\partial)$ is isomorphic to a  proper subalgebra $\fA_1$ of $ \fA_0$.
In particular, {\it over any point the contact homology algebra $H_*(\fA,\partial)$
  has no odd elements}.}
\end{remark}
\noindent\textsf{ Distinguishing contact structures on
pre-quantizations spaces}
 \medskip

 Formula (\ref{eq:Bourg}) can be used
   for distinguishing contact structures on pre-quantization
   spaces
   of certain symplectic manifolds, which have different
   Gromov-Witten invariants.
                         Here is an example.
                         \begin{proposition} \label{prop:Barlow}
                         Let $(M_1,\omega_1)$
                         and $(M_2,\omega_2,J_2)$ be two
                         symplectic $4$-manifolds with  integral
                         cohomology classes of their symplectic forms.
                         Suppose that for compatible
                         almost complex structures $J_1$ on $M_1$ and $J_2$
                         on $M_2$ there are no non-constant  $J_1$-holomorphic spheres
                         in
                        $ M_1 $,
                         while $M_2$ contains an embedded
                         $J_2$-holomorphic $(-1)$-sphere $S$.
                          Then  the pre-quantization spaces $(V_1,\xi_1)$ and $(V_2,\xi_2)$
are not contactomorphic.
\footnote{
Yongbin Ruan proved in \cite{Ruan-virtual}
that under the assumptions of Proposition \ref{prop:Barlow}   the
symplectic manifolds $(M_1,\omega_1)\times \C P^1$ and
 $(M_2,\omega_2)\times \C P^1$ are not symplectomorphic
 (and  not even deformationally equivalent), despite   the fact
 that    for a certain choice of $M_1$ and $M_2$
 (e.g. $M_1$ is the Barlow surface and $M_2=\C P^2\#8\overline{\C P^2}$), and for
 appropriate symplectic forms $\omega_1$ and $\omega_2$
 the underlying manifolds $V_1$ and $V_2$ are diffeomorphic. }
\end{proposition}
\begin{remark}
{\rm Even when the manifolds $M_1$ and $M_2$ are homeomorphic,
the prequantization spaces $V_1$ and $V_2$ are not diffeomorphic
(even not homotopy equivalent!)  for most choices   of symplectic forms
$\omega_1$ and $\omega_2$, and hence the statement of the theorem is trivial in these
cases. However, one can show that
for homeomorphic $M_1$ and $M_2$ the
symplectic forms can always be deformed in the class of symplectic forms compatible with
the chosen almost complex structures $J_1$ and $J_2$, in order
  to make $V_1$ and $V_2$ diffeomorphic.}
\end{remark}
To prove Proposition \ref{prop:Barlow}   we will ahow that
   the ``classical"   contact homology algebras
 $H^{\cont}_*(V_1,\xi_1) $ and $H^{\cont}_*(V_2,\xi_2) $ are not isomorphic.

Let  $S_0=S,S_1,\dots, S_m$ be the exceptional $J_2$-holomorphic
spheres in $M_2$. Then the cohomology classes $ S^*_0,\dots,
S^*_m,[\omega]$ are linearly independent, where  we denote by
$S^*_i$ the cohomology class Poincar\'e-dual to $[S_i]\in
H_2(M_2)$. Hence there exists a class $X\in H^2(M_2)$, such that
\begin{equation}
XS^*_0=1, \;X[\omega]=0,\;\;\hbox{and}\;\;XS^*_i=0\;\;\hbox{for}\;\; 1\geq i\geq m.
\end{equation}
Then the  potential $\bff_{M_2,J_2}(tX)$
coincides with $$\bff_{S,J_2|_S}(tX|_{S})=e^tz.$$
 Let us choose a basis of closed forms $\Delta_i, i=0,\dots,a$,
generating
  generating $H^*(M_2)$, so that one of the forms, say $\Delta_1$
  represents the class $X$. The form $\Delta_1$ then
     lifts to a  form $\wt \Delta_1$ such that
$\pi_*\wt \Delta_1=\Delta_1$.
 According to the formula
(\ref{eq:Bourg2}) we have
\begin{equation}
\bh_V( \tau_1\wt\Delta_1,u)= \frac{\tau_1}{2\pi } \int\limits_0^{2\pi}
e^{u_1-ilx} dx ,
\end{equation}
where $u=\sum\limits_0^au_j\Delta_j,\,u_j=\sum\limits_1^\infty
 (p_{k,j}e^{ikx}+q_{k,j}e^{-ikx}),\;l=\int\limits_{S}\omega.$

 Hence, the contact homology algebra $H^{\cont}_*(M_2,\xi_2)$,
 specialized     at the point $\tau\wt\Delta$ for $\tau\neq 0$
 is isomorphic to the contact homology algebra of the  standard
 contact
 Lens space $L(l,1)=\pi^{-1}(S)\subset V$, specialized    at the volume
form $\tau\wt\Delta|_{L(l,1)}$. It follows then from  Remark \ref{rem:lens-space}
   that $H^{\cont}_*(V_2,\xi_2)$,
 specialized     at any point  has
 no odd elements.   On the other hand, for any $2$-dimensional
 cohomology class $t\in H_2(M_1)$ we have
 $\bff_{M_1,J_1}(t)=0$, and hence for any $3$-form
 $\wt\Delta$ on $V_1$,
 the formula (\ref{eq:Bourg2}) implies that the Hamiltonian
 $\bh_{V_1} (\wt\Delta, u )$ vanishes as well, and  therefore
 the contact homology algebra  $H^{\cont}_*(M_1,\xi_1)$,
 specialized     at the point $\tau\wt\Delta,\tau\neq 0$
is a free graded commutative algebra generated by the odd variable
$\tau$, and even variables $p'_{k,j},q'_{k,j }$, which correspond
to even dimensional generators of $H^*(M_1)$. Hence the contact
manifolds $(M_1,\xi_1)$ and $(M_2,\xi_2)$ are not isomorphic. \qed

  \bigskip

 \noindent\textsf{ Subcritical symplectic manifolds}
  \medskip

   The content of this example is a result of our discussion
   with P. Biran and K. Cieliebak.

  In \cite{Donaldson} S. Donaldson  generalized  the Kodaira embedding theorem by
  proving that for any closed  symplectic manifold $(W,\omega)$
   with an integral cohomology class of the symplectic form there
   exists an integer $l>0$ such that the homology class dual to $l[\omega]$
   can be represented by an embedded symplectic hypersurface $W_0$.
  In fact, S. Donaldson proved a stronger result, which  together
  with
   an improvement by  Biran-Cieliebak   asserts that for a sufficiently large $l$
   the hypersurface $W_0$ can be chosen in such a way that in the complement
   $W\setminus W_0$ there exists a vector field $X$ with the following properties:
   \begin{description}
   \item{} $L_X\omega=-\omega$, where $L_X$ denotes the Lie derivative
   along $X$; in other words, $X$ is conformally symplectic and   contracting;
    \item{}  $X$ is forward integrable;
    \item{} $X$ is gradient-like for a Morse function $\phi:W\setminus W_0\to \R_+$,
    which coincides with $-\log\dist^2$ near $W$, where $\dist(x)$ is the     distance
    function from a point $x$ to $W_0$ with respect to some Riemannian metric.
    \end{description}

       The vector field $X$ retracts $W\setminus W_0$ to the Morse complex $K$
       of the function $\phi$,
       which is automatically isotropic for the symplectic form $\omega$
       (see  \cite{Eliash-Gromov-convex}), and, in particular, $\dim K\leq n$.
         Biran-Cieliebak call the pair $(W,W_0)$ {\it subcritical}
           if $\dim K<n$.
      They constructed  in \cite{Biran-Cieliebak}
    several interesting  examples of subcritical pairs, and conjectured
    that {\it if $(W,W_0)$ is subcritical, then $l=1$}.
    We sketch below  the proof of this conjecture.

First, let us observe that the contact structure $\xi$, defined by
the contact form $\a=
 X\hook\omega$ on the boundary  $V$ of a  small tubular neighborhood of
$W_0$, is equivalent to the contact structure which is defined on $V$ as the
pre-quantization space   of the symplectic manifold
$(W_0,l\omega)$. On the other hand, the condition, that the pair
$(W,W_0)$ is sub-critical implies that the contact manifold $(V,\xi)$
is itself  subcritical in the sense of Examlpe \ref{ex:Floer}.4 above,
i.e. it is  isomorphic to the strictly pseudo-convex
boundary  of a sub-critical Stein (or Weinstein) manifold with its canonical complex structure.
Let us recall (see \ref{rem:grading} above) that all SFT-objects, in particular Floer contact
homology $HC_*(V,\xi)$ and the contact homology algebra $H^{\cont}_*(V,\xi)$ are
graded by elements of $H_1(V)$.
Using  arguments  as in  the theorem of Mei-Lin Yau (see \ref{prop:MLYau} above) one can show
 that all  non-trivial elements
 in  the
contact homology algebra $H^{\cont}_*(V,\xi)$ of a subcritical contact manifold
$(V,\xi)$ may correspond only to $0\in H_1(V)$.
  On the other hand,
 it follows from Proposition \ref{prop:even-dim} above that $H^{\cont}_*(V,\xi)$ specialized at $0\in H^*(V)$ has
 non-trivial elements which correspond to the homology class of       the fiber in
 $H_1(V)$. Therefore, $l=1$.

\subsubsection{Computing  rational Gromov-Witten
invariants of $\C P^n$}\label{sec:computing}
We will  show  in this section how SFT  can be used for computing      rational
 Gromov-Witten invariants of $\C P^n$. Our method differs from
 traditional ways (see \cite{Kontsevich-CP2,Fulton-Pand,Graber-Pand,Vakil,Ruan-Tian})   for this computation.
  We will be simultaneously computing the rational potential of $\C^n$ and the rational
Gromov-Witten invariant of $\C P^n$ by a recursion
using Theorem \ref{thm:RSFT-master}

Let us choose basic forms in $\C^n$ as in the previous section,
i.e $\Delta =1$, and $\Theta$ is a volume form with  compact
support in $C^n\setminus 0=S^{2n-1}\times(0,\infty)$
  with $\int\limits_{\C^n}\Theta =1$. We denote by $\delta$ the restriction of
  $\Delta$ to $S^{2n-1}$.
We also assume that $\Theta$ splits into a product $\hat\theta\wedge\rho$, where $\hat\theta$
is pull-back of a unit volume form $\theta$ on $S^{2n-1}$, and $\rho$ is a
compactly supported form in $(0,\infty)$.
Set
\begin{equation}
\begin{split}
\bh^1(t_0,q,p))&=\frac{\partial \bh}{\partial
\tau}(t_0\delta+\tau\theta,q,p)|_{\tau=0}\\ &=\frac{1}{2\pi
}\int\limits_0^{2\pi}\widehat\bff_{\C P^{n-1}}(t_0+u_0(x),u_2,\dots,
u_{2n-2},e^{-ix})dx,\\
\end{split}
\end{equation}
 where  $$\widehat\bff_{\C P^{n-1}}(t_0,\dots,t_{2n-2},z)=
 \frac{\partial\bff_{\C P^{n-1}}}{\partial t_{2n-2}}(t_0,\dots,t_{2n-2},z),$$
 $\bff_{\C P^{n-1}}(t_0,\dots,t_{2n-2},z)$ is the rational
Gromov-Witten invariant  of $\C P^{n-1}$, and
$$u_{2j}(x)=\sum\limits_1^\infty
p_{k,2j}\,e^{ikx}+q_{k,2j}\,e^{-ikx}.$$ Then the equation
(\ref{eq:RSFT-master}), which determines $
\bff(t_0,t_{2n},p)=\bff_{\C^n}(t_0\Delta+t_{2n}\Theta,p)$ takes
the form
\begin{equation}\label{eq:Cn}
\frac{\partial\bff}{\partial t_{2n}}(t_0,t_{2n},p) =\frac{1}{2\pi
}\int\limits_0^{2\pi} \widehat\bff_{\C P^{n-1}}(t_0+u_0(x),u_2,\dots,
u_{2n-2},e^{-ix})dx\big|_{L_\bff},
\end{equation}
where $$L_{\bff}=\{q_{k,2j}=k\frac{\partial \bff}{\partial
p_{k,2n-2j-2}}(t_0,t_{2n},p)\}.$$
 Together with the initial data
\begin{equation}
\bff(t_0,0,p)=
\begin{cases}
p_{1,0},&\hbox{if}\;\;n=1;\\
0,&\hbox {otherwise} \\
\end{cases}
\end{equation}
the equation (\ref{eq:Cn}) provides a recursive procedure for
computing coefficients $f_j(t_0,p)$ of the expansion
$$\bff(t_0,t_{2n},p)=\sum\limits_0^\infty f_j(t_0,p)t_{2n}^j.$$

For instance for $n=1$ we have (see Example \ref{ex:circle})
$\bh^1=\frac{t_0^2}{2}+\sum\limits_1^\infty p_kq_k$, where we
write
$p_k,q_k$ instead of $p_{k,0},q_{k,0}$, and hence the equation
(\ref{eq:Cn}) takes the form
\begin{equation}
\frac{\partial\bff}{\partial t_{2 }}(t_0,t_{2 },p)=\frac{t_0^2}{2}
+ \sum\limits_0^\infty k p_k\frac{\partial\bff}{\partial
p_k}(t_0,t_{2 },p)
\end{equation}
with the initial data $\bff(t_0,0,p)=p_1$. This linear first order
PDE is straightforward to solve, and we get $$\bff(t_0,t_{2
},p)=\frac{t_2t_0^2}{2}+e^{t_2}p_1.$$
\bigskip

For $n=2$ the Hamiltonian $\bh$ is given by the formula
(\ref{eq:3-sphere}), and we have
$$\bh(t_0,\tau,p)=\tau\bh^1(t_0,p).$$ Thus the equation for the
potential of $\C^2$ has the form
\begin{equation}\label{eq:C2}
\begin{split}
&\frac{\partial\bff}{\partial t_{4}}(t_0,t_{4},p)
 =\frac{t_0^2}{2}+\sum\limits_{k\geq 1}k
\frac{\partial\bff}{\partial p_{k,2}}p_{k,0}\\
&+\sum\limits_{t,s\geq 0}\mathop{\mathop{\sum
\limits_{i_1,\dots,i_s\geq 0}}\limits_{j_1,\dots j_t\geq 0}}
\limits_{\sum\limits_1^s li_l-\sum\limits_1^tmj_m=1} \frac{ \left(
\frac{\partial\bff} {\partial p_{1,0}}\right)^{i_1}\dots
 \left(s \frac{\partial\bff}{\partial p_{s,0}}\right)^{i_s}p_{1,2}^{j_1}\dots
 p_{t,2}^{i_t}}{i_1!\dots i_s!j_1!\dots
j_t!};\\
&\bff(t_0,0,p)=0\,.\\
\end{split}
\end{equation}
Hence, we get
\begin{equation}
\begin{split}
\bff(t_0,t_4,p)=&t_4(\frac{t_0^2}{2}+p_{1,2})+\frac{t_4^2}{2!}p_{1,0}+
\frac{t_4^3}{3!}(p_{2,2}+\frac12 p_{1,2}^2)\\
&+\frac{t_4^4}{4!}(2p_{2,0}+p_{1,2}p_{1,0})+\dots\;.
\end{split}
\end{equation}

To compute $\bff_{\C^n}$ for $n>2$ we need to know $\bff_{\C
P^{n-1}}$. So to complete the recursion we will explain now how to
express the rational Gromov-Witten invariant $\bff_{ \C P^{n-1}}$
through $\bff_{\C^n}$.

\bigskip
First of all we split, as it is described in Example
\ref{ex:section}
 above, $\C P^n$ along the boundary of a tubular neighborhood
of $\C P^{n-1}\subset \C P^n$ into two completed symplectic
cobordism $W_1=\C^n$ and $W_2=\C P^{n }\setminus {x}$, where we
introduce on $W_2$ a complex structure of the holomorphic line
bundle over $\C P^{n-1}$ determined by the hyperplane section $\C
P^{n-2}\subset\C P^{n-1}$. We will denote by $\bff_1$ and $\bff_2$
the potentials for $W_1$ and $W_2$, respectively.

Let $\Delta_0,\dots,\Delta_{2 n-2}$ be closed forms representing
the standard basis of $H^*(\C P^{n-1})$. We will keep the same
notation for the pull-backs of these forms to $W_2$. Let
$\Delta_{2n}$ be  a closed $2n$-form with a compact support, which
generates $$\Ker\big(H^*_{\comp}(W_2)\to H^*(W_2)\big).$$
 We are interested
in the potential
$\bff_2(t_0,\dots,t_{2n},q)=\bff_2(\sum\limits_{i=1}^{n }
t_{2i}\Delta_{2i},q)$. First of all notice that by dimensional
reasons the moduli spaces of holomorphic curves which project  to
non-constant curves in $\C P^{n-1}$ do not contribute to the
potential $$\bff_2(t_0,\dots,t_{2n-2},0,q),$$ and hence we have
\begin{equation}
\bff_2(\sum\limits_{i=1}^{n-1} t_{2i}\Delta_{2i},q)= z\sum
\limits_{i=0}^{n-1}q_{1,2i} \mathop{\sum\limits_{(s_1,\dots,s_{n-1})}}\limits_{  \sum
s_j(j-1)=n-i-1}\prod\limits_{j=1}^{n-1}\frac{t_{2j}^{s_j}}{s_j!}\;.
\end{equation}
In particular, for $n=2$ we get
\begin{equation*}
\bff_2(t_0\Delta_0+t_2\Delta_2,q)=ze^t_2q_{1,2}
\end{equation*}
 One can recover $\bff_2(t_0,\dots,t_{2n},q)$ for $t_{2n}\neq 0$
using the equation (\ref{eq:RSFT-master}), as we did it above for
$W_1=\C^n$. However, for the purpose of our  computation of
Gromov-Witten invariant $\bff_{\C P^n}$ this is not necessary, as
we can proceed as follows.

Notice that the above chosen  forms
$\Delta_2, \Delta_{2n-2}$ extend to $\ CP^n$. On the other hand, we
will choose a  volume form $\Delta_{2n}$ on $\C P^n$ to be
supported in the affine part, so that the restriction $\Delta_{2n}|_{\C^n}$
coincides with the form $\Theta$ introduced above. Then Theorem
(\ref{thm:RSFT-composition}) implies that
\begin{equation}\label{eq:CPn}
\begin{split}
&\bff_{\C P^n}(\sum\limits_{i=1}^{n } t_{2i}\Delta_{2i})=\\
 &\big(\bff_1(t_0,t_{2n},p)+
 \bff_2(t_0,\dots,t_{2n-2},0,q)-\sum\limits_{i+j=n-1}
\sum\limits_1^\infty \frac1{k} p_{k,2i}q_{k,2j}\big)\big|_L\;\;,\\
\end{split}
\end{equation}
where
\begin{equation}\label{eq:Lag-CPn}
L=
\begin{cases}
p_{1,2i}=& z  \mathop{\sum\limits_{(s_1,\dots,s_{n-1})}}\limits_{
\sum
s_j(j-1)=i}\prod\limits_{j=1}^{n-1}\frac{t_{2j}^{s_j}}{s_j!}\;;\\
p_{k,2i}=&0,\;\;\hbox{if}\;\; k>1\;; \\
q_{k,2i}=&k\frac{\partial\bff_1}{\partial
p_{k,2(n-i-1)}}(t_0,t_{2n}, p)\,.\\
\end{cases}
\end{equation}

Plugging expressions from (\ref{eq:Lag-CPn}) into equation
(\ref{eq:CPn}) we get
\begin{equation}\label{eq:CPn-Cn}
\bff_{\C
P^n}(t_0,\dots,t_{2n})=\bff_{\C^n}(t_0,t_{2n},p)\big|_{L_1}\;,
\end{equation}
where
\begin{equation*}
L_1=
\begin{cases} p_{1,2i}&= z  \mathop{\sum\limits_{(s_1,\dots,s_{n-1})}}\limits_{
\sum
s_j(j-1)=i}\prod\limits_{j=1}^{n-1}\frac{t_{2j}^{s_j}}{s_j!}\;;\\
p_{k,2i}&=0\;\;\hbox{if}\;\;
k>1\,.\\
\end{cases}
\end{equation*}
Indeed, two last terms in the formula (\ref{eq:CPn}) cancel each
other (as it always happens when $\bff_2$ is linear with respect
to $q$-variables).
 For instance, for $n=1$ we get
\begin{equation}
\bff_1(t_0,t_2,p)=\big(\frac{t_0^2t_2
}2+e^{t_2}p_1\big)\big|_{p_1=z}=\frac{t_0^2t_2}2+e^{t_2}z\,.
\end{equation}

For $n=2$ we have
\begin{equation}
L_1=
\begin{cases}
 p_{1,0}&= ze^{t_2}\,;\\
p_{k,i}&= 0,\quad\hbox{ for all other}\quad k,i\,,\\
\end{cases}
\end{equation}
and hence
\begin{equation}\label{eq:CP2-final}
\bff_{\C P^2}(t_0,t_2,t_4)=\bff_{\C^n}(t_0,t_4,
 ze^{t_2},0,\dots)\,.
\end{equation}

\begin{remark} {\rm The method which we used  above
 for computing of the rational
potential of $\C P^n$, when  applied to an arbitrary
symplectic manifold $W$,   allows us to express
  $\bff_W$ through the potential of the affine part $W\setminus
  M$.
  The latter computation seems tractable when  the Weinstein manifold
  $W\setminus M$ is
  subcritical (see Section \ref{sec:splitting} above), i.e. when
  its isotropic skeleton does not have Lagrangian cells.
  On the other hand, when Lagrangian cells are present
  this problem is related to central questions of Symplectic
  topology.}
  \end{remark}

\subsubsection{Satellites}\label{sec:satellites}
Let $(V,\xi=\{\a=0\})$ be a  contact manifold, $(W
= V\times \R, d(e^t\a))$  its symplectization, and   $J$   a
compatible translation-invariant almost complex structure on $W$.
   In this section  we will show that the
 homological Poisson super-algebra $H_*(\fP ,  d^{\bh}  )$ comes equipped with some
additional structures, rather unfamiliar in abstract Poisson
geometry. Namely, the counting  of genus $g$ curves with a fixed
complex structure and with a fixed configuration of $n$ points
gives rise to an odd $n$-linear totally symmetric poly-form
$\bh^{g,n}$ on the Poisson super-space $\bV$ underlying $\fP $. The
poly-form descends well to the homology and thus yields another
invariant of the contact structure which we call the genus-$g$
$n$-point {\em satellite} of the Poisson structure.

Let us denote by
$\overline{\cM_{g,m}(V)/\R}$   the compactified moduli space of stable connected
$J$-holomorphic curves in $W$ which are characterized by the
arithmetical genus $g$ and by the total number $m$ of punctures and
marked points numbered somehow by the indices $1,...,m$ (see Section \ref{sec:compact} above). We
emphasize that the moduli space in question contains equivalence
classes of all such curves, and in particular, may have infinitely
many connected components corresponding to different homotopy
types of curves in $W$ and different numbering of the $m$
markings. Let $\DM_{g,n}$ be the Deligne-Mumford compactification
of the moduli space of genus $g$ Riemann surfaces with $n$ marked
points. For any $g,n$ with $2g-2+n >0$ and $l\geq 0$ there is a
natural {\em contraction map} $\ct :\overline{\cM_{g,n+l}(V)/\R}\to \DM_{g,n}$
defined by forgetting the map to $W$ and the last $l$ markings
followed by the contraction of those components of the curve which
have become unstable.  Given a differential form
$\tau$ on $\DM_{g,n}$ we will denote by $\ct^*\gt$ its  pull-back
 to $\overline{\cM_{g,m}(V)/\R}$.

Let $u =(p,q,t)$ denote a point in $\bV$, that is $p,q$ and $t$ are
(closed) differential forms on $\cP^{-},\cP^{+}$ and $V$
respectively. We will denote $\ev_i^*{u}, \ i=1,...,m$, the
pull-back by the evaluation map $$\ev_i: \overline{\cM_{g,m}(V)/\R}\to
(\cP^{-}\cup \cP^{+}\cup V)$$ at the $i$-th marking. Let us emphasize the point that we
are
treating here  the marked points and punctures on equal footing.

Let $\gd u\in \bV$ be a tangent vector to $\bV$ at a point $u\in \bV$.
We introduce the formal function
 \begin{equation}
\bh^{g,n}_{\gt}:=\frac{1}{n!}
\sum_{l=0}^{\infty}\frac{1}{l!}\int\limits_{\overline{\cM_{g,n+l}(V)/\R }}\
\ct^*{\gt}\w \ev_1^*\gd u\w...\w\ev^*_n\gd
u\w\ev_{n+1}^*u\w...\w\ev_{n+l}^*u .
\end{equation}
 It is a super-symmetric  $n$-linear
form in $\gd u$ with coefficients depending on the application
point $u$.

Let $d^{\bh}(f)$ denote the Lie derivative of a tensor field $f$
 along the odd Hamiltonian
vector field $d^{\bh}$ on $\bV$ with the Hamilton function $\bh$.

\begin{proposition}
Let $\gt $ be a top degree form on $\DM_{g,n}$. Then
$d^{\bh}(\bh_{\gt}^{g,n})=0$. If the top degree form $\gt=d\a$ is exact
then $\bh^{\gt}_{g,n}=d^{\bh}( \bh^{\a}_{g,n})$. In particular, the
tensor field $\bh^{g,n}_{\gt}$ descends to the homology    algebra
$H_*(\fP, \p)$ into a satellite which depends only on the total
volume of $\gt$.
\end{proposition}

This follows from the Stokes formula applied to
$\bh^{d\gt}_{g,n}=0$ and respectively  to $\bh^{d\a}_{g,n}$.
Codimension $1$ boundary strata of the moduli space $\cM_{g,m}(V)/\R$
are formed by  stable curves of height $2$. Most of the strata do not contribute to
the Stokes formula since they are mapped by the contraction map to
complex codimension $1$ strata of the Deligne-Mumford space
$\DM_{g,n}$, where $\gt$ and   $\a$ restrict to $0$ for
dimensional reasons. Exceptions to this rule occur only if one of
the two curves which form the stable curve is to be contracted. It is therefore a
sphere with  glued to the other level of the stable curve along at precisely
one end, and which have at most one marked points or ends with the index $\leq n$,
 and with any
number of ends or marked points with indices $>n$. Contributions
of such curves to the Stokes formula is expressed bi-linearly via
the $1$-st or $2$-nd derivatives of the Hamilton function $\bh$
and the satellite. It is easy to see that the whole expression is
interpreted correctly as the  Lie derivative of the tensor field
$\bh^{\gt}_{g,n}$ along the Hamiltonian vector field $d^{\bh}$. \qed

\medskip

We will assume further on that $\gt $ is normalized to the total
volume $1$ and will often drop it from the notation for the
satellite $\bh^{g,n}$.

\medskip

Let us consider now a directed symplectic  cobordism
$W=\ora{V_-V_+}$ between two contact boundaries $V_{\pm}$. Then we
have the Hamilton function $\widehat\bh
=\bh^{+}-\bh^{-}$ and the satellites
$\widehat\bh^{g,n}=\big(\bh^{g,n}\big)^+-
\big(\bh^{g,n}\big)^-$ defined  as elements of
the algebra $\widehat \fL$, which in the case when the cobordism
is a concordance just equal to the tensor product of the Poisson
algebras $\fP_{\pm}$. Also, we have the potential $\bff
(p_{-},q_{+},t) $ counting rational $J$-holomorphic curves in $W$
which defines a Lagrangian correspondence between $\fP^{\pm}$
invariant under the vector field $ d^{\widehat{\bh}}$ with the Hamilton function
$\widehat\bh$. Finally, using the moduli spaces $\cM_{g,m}(W)$ of
$J$-holomorphic curves in the cobordism, we can introduce the
satellites $\bff^{g,n}_{\gt}$ as symmetric $n$-forms on the space
$(p_{-},q_{+},t)$-space parameterizing the Lagrangian
correspondence. Then the arguments similar to the above proof of
the proposition, but applied this time to $\bff^{g,n}_{d\gt}=0$,
show that {\em the restriction of $\widehat\bh_{\gt}^{g,n}$ to the
Lagrangian correspondence defined by $\bff$ coincides with the Lie
derivative of $\bff^{g,n}_{\gt}$ along the vector field $d^{\widehat\bh}$
restricted to the Lagrangian correspondence}
(comp. Theorem \ref{thm:RSFT-cobordism} above).
In this sense the
Lagrangian correspondences defined by symplectic cobordisms
preserve the satellite structures defined by $\big(\bh^{g,n}\big)^{\pm}$ on
the homology $H_*(\fP_{\pm},d^{\bh^{\pm}})$. In particular, {\em the
satellite structures of a contact manifold $V$ on the homology
$H_*(\fP,d^{\bh})$ depend only on the contact structure.}

\medskip

The following discussion is the first steps in the study of the
geometric structure defined by the satellites.

Let $\bh^{g,n}_{\a_1,...,\a_n}$ denote components of the satellite
tensors on $\fP$. Using the Poisson tensor $\pi^{\mu\nu}$ we can
couple two satellites with respect to some indices:
$$\bh^{g',n'+1}_{...\mu}\pi^{\mu\nu}\bh^{g'',n''+1}_{\nu ...}.$$
Similarly, we can couple two indices in $\bh^{g-1,n+2}$ with two
indices in the $2$-nd differential $\gd^2\bh $ of the Hamilton
function $\bh$.

\begin{proposition} If $g=g'+g''>0$ then the
coupling of $\bh^{g',n'}$ and $\bh^{g'',n''}$ is a Lie derivative
along $\p$ and thus vanishes in the homology $H_*(\fP ,\p )$.
Similarly, the coupling of $\bh^{g-1,n+2}$ with $\gd^2\bh$
vanishes in the homology $H_*(\fP,\p)$. \end{proposition}

The proof is based on some famous but non-trivial property of the
spaces $\DM_{g,n}$ with $g>0$ to have complex codimension one
strata homologically independent. Such strata correspond to
different ways of cutting a $(g,n)$-surface along one circle and
can be identified either with $\DM_{g',n'+1}\times
\DM_{g'',n''+1}$ where $g'+g''=g, n'+n''=n$ or with
$\DM_{g-1,n+2}$. The independence property implies that a volume
form $\gt$ on the stratum, say  $\gt'\otimes\gt''$ in the first
case, can be obtained as the restriction of a closed codegree two
form $\go$ on $\DM_{g,n}$ which have exact (or even zero, for
suitable choices of $\gt$) restrictions to all other
codimension-$1$ strata in $\DM_{g,n}$. Applying the Stokes formula
to $0=\bh^{g,n}_{d\go}$ we find that the coupling of
$\bh^{g',n'+1}_{\gt'}$ and $\bh^{\gt'',n''+1}_{\gt''}$ (or --- in
the second case --- of $\bh^{g-1,n+2}_{\gt}$ and $\gd^2\bh$)
equals $d^{\bh}( \bh^{g,n}_{\go})$.

\medskip

\begin{remark} {\em
To the contrary, coupling $\bh^{0,3}$ with itself via one index is
not, in general, a $d^{\bh}$-derivative, but instead the following
triple sum is:
\begin{equation*}
\begin{split}
& \bh^{0,3}_{\a\gb\mu}\pi^{\mu\nu}\bh^{0,3}_{\nu\g\d}+
   (-1)^{(\deg{\a}+\deg{\gb})\deg{\g}}
   \bh^{0,3}_{\g\a\mu}\pi^{\mu\nu}\bh^{0,3}_{\nu\gb\d}+\\
  & (-1)^{\deg{\a}(\deg{\gb}+\deg{\g})}
   \bh^{0,3}_{\gb\g\mu}\pi^{\mu\nu}\bh^{0,3}_{\nu\a\d} \equiv 0
  \, .\\
   \end{split}
   \end{equation*}
This follows from the property of the $3$ boundary strata in
$\DM_{0,4}$ to represent the same homology class (use the Stokes
formula for $\go =1$). In fact $\bh^{0,3}$ coincides with the
$3$-rd differential $\gd^3\bh /6$ of the Hamilton function, and
the above Jacobi-like identity can be derived by $4$
differentiations of $\{ \bh ,\bh \}=0$ in $\a,\gb,\g,\d$. One can
interpret the integrability property $(d^{\bh})^2=0$ of the odd vector
field $d^{\bh}$ on $\bV$ as a homotopy Lie super-algebra structure on
$\Pi \bV^*$, the dual space with changed parity. The identity in
question corresponds to the Jacobi identity for the remnant Lie
super-algebra structure in homology. }
 \end{remark}

It is sometimes convenient to extend the definition of genus $0$
satellites to unstable values of $n$ by $\bh^{0,n}=\gd^n\bh /n!$
for $n=0,1,2$. Also, one can define the function $\bh^{1,0}$ at
least locally as a potential for $\bh^{1,1}$, using the following

\begin{proposition} The differential $1$-form $\bh^{1,1}$ is closed.
\end{proposition}

Indeed, the partial derivatives $\gd_{\mu}\bh^{1,1}_{\nu}$ and
$\gd_{\nu}\bh^{1,1}_{\mu}$ are identified with the satellites
$(\bh^{1,2}_{\go})_{\mu\nu}$ corresponding to the $2$-form $\go $
on $\DM_{1,2}$ pulled-back from $\DM_{1,1}$ by forgetting the
$1$-st and respectively the $2$-nd marked point. But the two maps
$\DM_{1,2}\to \DM_{1,1}$ coincide.

\medskip

It would be interesting to study other general properties of
satellites which may depend on more sophisticated geometry of
Deligne-Mumford compactifications. For instance, what can be said
about Poisson brackets among the functions $\bh^{g,0}$?

\medskip

We complete the section with the computation of the satellites in the
example $V=S^1$. Let $t=t_01+t_1d\phi $ denote the general
harmonic form on $S^1$, $\gd t=\gt_0 1+\gd t_1\d\phi $,
$u(x)=t_0+\sum p_ke^{ikx}+ q_ke^{-ikx}$, $\gd u = \gd t_0+\sum \gd
p_ke^{ikx}+\gd q_k e^{-ikx}$.

\begin{proposition}
For $2g-2+n \geq 0$ we have
\[ \bh^{g,n+1}\ =\ \frac{\gd t_1}{2\pi\ n!}\int\limits_0^{2\pi} \ (u_{xx})^g \
(\gd u)^n \ dx .\]
\end{proposition}

Let us begin with the remark that the formula does not (and
cannot) contain $t=t_0+t_1\phi$ because $\deg t < 2$, and
therefore pushing forward from $\cM_{g,m+1}(V)/\R \to
\cM_{g,m}(V)/\R$ by forgetting the corresponding marked point
would send $t$ to $0$. Exceptions to this rule could occur only if
$\cM_{g,m}(V)$ were ill-defined, that is only in the case of
constant maps with "unstable" indices, $2g-2+m\leq 0$, which has
no effect on the satellites with "stable" indices. On the other
hand the factor $\d t_1$ is (and must be) present in the formula
since the dimension of the moduli spaces is odd. With this
information in mind, the enumerative question equivalent to
computation of the satellites can be described as follows. On a
Riemann surface $\Sigma$ of genus $g$, we are given a divisor $D$
of $n$ distinct points with (possibly zero) multiplicities
$m_1,...,m_n$. We have to count the divisors $l_1P_1+...+l_gP_g$
which in the sum with $D$ form the divisor of a rational function.
(In particular, the degree $\sum m_i+\sum l_j$ of the total
divisor must vanish.)   The answer to this question is equal to
the degree of the Abel-Jacobi map $\Sigma^g\to J_{\Sigma}$ defined
by integration of holomorphic differentials
$\vec{\go}=(\go_1,...,\go_g)$ on $\Sigma$ as
\[ (P_1,...,P_g)\mapsto l_1\int^{P_1}\vec{\go}+...+l_g \int^{P_g}\vec{\go}
.\] When the multiplicities $(l_1,...,l_g)=(1,...,1)$, the degree
equals $g!$ (it is well-known that $S^g\Sigma^g\to J_{\Sigma}$ is
a bi-rational isomorphism). For arbitrary $(l_1,...l_g)$ the
Abel-Jacobi map has the Jacobi matrix $[l_j \go_i (P_j)]$. Thus
the degree equals $l_1^2...l_g^2 g!$. Taking these answers as the
coefficients in the generating function on the variables $t_0,
p_l, q_{-l}$ corresponding to $l=0,l>0$ and $l<0$ we arrive at the
factor $u_{xx}^g$. The other factor $(\gd u)^n/n!$ is similarly
accountable for all possible choices of multiplicities
$m_1,...,m_n$ in the divisor $D$. The contour integration of the
product couples the choices with $m_1+...+m_n+l_1+...+l_g = 0$.


\begin{thebibliography}{99}

\bibitem{AbbasHofer}
C.~Abbas and H.~Hofer,
  Holomorphic curves and global questions in contact geometry,
  to appear in Birkh\"auser.
\bibitem{Arnold-conjecture} V.I. Arnold,        {Sur une propri\'et\'e topologique des
applications globalement canoniques de la m\'echanique classique},
      {\it C.\ R.\ Acad. Paris}, 261(1965), 3719--3722.
\bibitem{Arnold-Maslov} V. I. Arnold, On a characteristic class entering in quantization
conditions,
{\it Funct. Anal. and Applic.}, 1(1967), 1--14.
\bibitem{Arnold-steps}  V. I. Arnold,       {First steps in symplectic topology},
{\it Russian Math. Surveys}, 41(1986), 1--21.

\bibitem{Bennequin} D. Bennequin,
Entrelacements et \'equations de Pfaff, {\it Ast\'erisque}, 106--107(1983).
\bibitem{Biran-Cieliebak} P. Biran and K. Cieliebak, Symplectic
topology on subcritical manifolds, preprint 2000.
\bibitem{Brieskorn} E. Brieskorn, Beispiele zur
Differentialtopologie von Singularit\"aten, {\it Invent. Math.},
{\bf 2}(1966), 1--14.
\bibitem{Chekanov} Yu. Chekanov, Differential algebra of a Legendrian link.
Preprint 1997.
\bibitem{Conley-Zehnder} C. Conley and E. Zehnder,
        {The {B}irkhoff--{L}ewis fixed point theorem and a conjecture of {V}.{I}.
{A}rnold},
    {\it Invent. Math.}, 73(1983), 33--49.
    \bibitem{Donaldson-invariants} S.K. Donaldson,
    Polynomial invariants for smooth four-manifolds,
{\it Topology}, {\bf 29}(1990),  257--315.

     \bibitem{Donaldson} S.K. Donaldson, Symplectic
submanifolds and almost-complex geometry, {\it J. Diff. Geom.},
{\bf 44}(1996), 666-705.
\bibitem{Eliash-ICM} Y. Eliashberg, Invariants in contact
topology, {\it Proc. of ICM-98, Berlin}, Doc. Math., 1998,
327-338.

\bibitem{Eliash-Stein} Y. Eliashberg, Topological characterization of Stein manifolds
of complex dimension $>2$, {\it Int. J. of Math}, {\bf 1}(1991),
29--46.

 \bibitem{Eliash-Gromov-convex} Y. Eliashberg and M. Gromov, Convex symplectic manifolds, {\it Proc.
of Symp. in Pure Math.}, vol.52(1991), Part 2, 135--162.
 \bibitem{EHS} Y. Eliashberg, H. Hofer and S. Salamon, Lagrangian intersections
in contact geometry, {\it Geom. and Funct. Anal.}, 5(1995), 244--269.
\bibitem{Etnyre-Sabloff} J. Etnyre and J. Sabloff,
Coherent Orientations and Invariants of Legendrian Knots, preprint
2000.
   \bibitem{EH-3ball} Y. Eliashberg, H. Hofer,  A Hamiltonian characterization
   of the three-ball, {\it Differential Integral Equations},
   {\bf 7}(1994),  1303--1324.

\bibitem{Floer} A. Floer,
The unregularised gradient flow of the symplectic action,
      {\it Comm. Pure Appl. Math.},
     {\bf 41}(1988),
  775-813.
\bibitem{Floer-Hofer} A. Floer,  H. Hofer, Coherent
orientations for periodic orbit
problems in symplectic geometry, {\it Math. Z.}, {\bf  212}(1993), 13--38.
\bibitem{Fukaya-Ono}
 K. Fukaya and K. Ono, Arnold conjecture and Gromov-Witten invariants, preprint 1996.
 \bibitem{Fukaya-Ono2} K. Fukaya,  K. Ono,  Arnold conjecture and Gromov-Witten invariant
for general symplectic
manifolds, in {\it The Arnoldfest (Toronto, ON, 1997)}, 173--190, Fields
Inst. Commun., 24, Amer. Math. Soc., Providence, RI,
1999.

  \bibitem{FKO3}   K.Fukaya, K.Ono, Y.-G.Oh and H.Ohta, to appear.
  \bibitem{Fulton-Pand} W. Fulton and R. Pandharipande, Notes on stable
  maps and quantum cohomology, {\it Proc. of Symp. in Pure Math.},
  {\bf 62}(1995), Part 2, 45--96
  \bibitem{Gathmann}A. Gathmann,
  Absolute and relative Gromov-Witten invariants of very ample
  hypersurfaces, preprint 1999.
  \bibitem{Getzler} E. Getzler, Topological recursion relations in genus $2$, in {\it Integrable systems
   and algebraic geometry} (Kobe/Kyoto, 1997), 73--106, World Sci. Publishing  1998.
\bibitem{Giroux} E. Giroux, Une structure de contact,  m\^eme tendue est plus ou
moins tordue,
{\it Ann. Scient. Ec. Norm. Sup.}, {\bf  27}(1994), 697--705.

\bibitem{Givental-toric} A. Givental,
A symplectic fixed point theorem for
 toric manifolds,
 {\it The Floer memorial volume}, 445--481,
 Progr. Math.,
 133, Birkhäuser, Basel, 1995.
 \bibitem{Givental-Maslov} A. Givental,
 Nonlinear generalization of the Maslov index, {\it Adv. Soviet Math.}, {\bf 1}(1990),
  71--103.
\bibitem{Givental} A. Givental.
Homological geometry and mirror symmetry. Proc. Int. Congress of
Math., Z\"urich-1994. Birkh\"auser, 1995, v. 1, 472-480.


\bibitem{Givental2} A. Givental.
Homological geometry I: Projective hypersurfaces. Selecta Math.,
New series 1 (2), 1995, 325-345.

\bibitem{Givental3}A. Givental.
 Equivariant Gromov - Witten invariants. Intern. Math. Res. Notices, 1996,
No.13, 613 - 663.
\bibitem{Givental1} A. Givental and B. Kim.
Quantum cohomology of flag manifolds and Toda lattices. Commun.
Math. Phys. 168 (3), 1995, 609-641.

\bibitem{Graber-Pand} T. Graber and R. Pandharipande, Localization
of virtual classes, {\it Invent. Math.},{\bf  135}(1999),    487--518.
\bibitem{Gray} J.W. Gray, Some global properties of contact structures,
{\it Annals of Math.}, {\bf 69}(1959), 421--450.
\bibitem{Gromov-holomorphic} M. Gromov, Pseudo-holomorphic curves in symplectic manifolds,
{\it Invent. Math.}, {\bf 82}(1985), 307--347.
\bibitem{Gro-PDR} M. Gromov, Partial Differential Relations, Springer-Verlag,
1986.
\bibitem{Hofer-Weinstein} H. Hofer, Pseudo-holomorphic curves and Weinstein conjecture
in dimension three, {\it Invent. Math.}, {\bf 114}(1993), 515--563.
 \bibitem{Hofer-Wysocki-Zehnder}    H. Hofer, K. Wysocki and E. Zehnder, Properties of
pseudo-holomorphic curves in symplectisations. I. Asymptotics.
{\it  Ann. Inst. H. Poincar\'e, Anal. Non Lin\'eaire}, {\bf
13}(1996), 337--379.

\bibitem{HoferWysockiZehnder:94d}
H.~Hofer, K.~Wysocki, and E.~Zehnder,
  The dynamics on a strictly convex energy surface in {${\bf R}^4$},
  {\it Annals of Math.}, {\bf 148}(1998),  197-289.
\bibitem{HoferWysockiZehnder:99}
H. Hofer, K. Wysocki, and E. Zehnder,
  Finite Energy Foliations of Tight Three-Spheres and Hamiltonian
Dynamics, preprint 1999.
\bibitem{Ionel} E.-N. Ionel, Topological recursive relations in
$H^{2g}(M_{g,n})$, preprint 1999.
\bibitem{Ionel-Parker2} E.-N. Ionel and T.H. Parker
 Gromov-Witten invariants of symplectic sums.{\it  Math. Res. Lett.}
 {\bf 5} (1998),  563--576.
 \bibitem{Ionel-Parker} E.-N. Ionel and T.H. Parker,
 Relative Gromov-Witten Invariants, preprint 1999.
\bibitem{Kollar} J. Kollar,
Rational curves on algebraic varieties, Springer-Verlag, 1996.

\bibitem{Kontsevich-quant} M. Kontsevich, Deformation quantization of Poisson manifolds, I,
preprint 1997.
\bibitem{Kontsevich-CP2} M. Kontsevich,  Enumeration of rational
curves via torus action, in {\it The Moduli Space of Curves}, R.
Dijgraaf, C. Faber and G. van der Geer eds., Birkhauser,1995,
335--368.
\bibitem{Konts-Manin} M. Kontsevich, Yu. Manin,
 Gromov-Witten classes, quantum cohomology, and enumerative
 geometry, {\it Commun.Math.Phys.}, {\bf  164} (1994) 525-562.

\bibitem{Li-Tian} J. Li and G. Tian, Virtual moduli cycles and Gromov-Witten invariants of
general symplectic
manifolds, in {\it Topics in symplectic $4$-manifolds (Irvine, CA, 1996)}, 47--83,
First Int. Press Lect. Ser., I, Internat. Press,
Cambridge, MA, 1998.
\bibitem{Liu-Tian} G. Liu and G. Tian, Floer homology and Arnold conjecture,
{\it J. Diff. Geom.}, {\bf 49}(1998),  1--74.
\bibitem{Lutz} R. Lutz,  Structures de contact sur les fibr\'es principaux en cercles de
dimension $3$
{\it  Ann. Inst. Fourier},   XXVII, {\bf 3}(1977), 1--15.
              \bibitem{Martinet}J. Martinet, Formes de contact sur les vari\'et\'es de dimension $3$,
{\it  Lecture Notes in Math.,} {\bf 209}(1971), 142--163.
\bibitem{McDuff} D. McDuff, The virtual moduli cycle, in{\it  Northern California Symplectic
Geometry Seminar}, 73--102,
Amer. Math. Soc. Transl. Ser. 2, 196, Amer. Math. Soc., Providence, RI,
1999.
    \bibitem{Mori-Miayoka} Y. Myayoka and S. Mori, A numerical
    criterion for uniruleness, {\it Annals of Math}, {\bf
    124}(1986), 65--69.
    \bibitem{Morita} S. Morita, A topological classification of
    complex structures on $S^1\times \Sigma^{2n-1}$,
    {\it Topology}, {\bf 14}(1975), 13--22.
\bibitem{Ng} L. Ng, On Invariants of Legendrian Knots, preprint
2000.
\bibitem{Robbin-Salamon} J.  Robbin, D. Salamon,The Maslov index for
paths, {\it Topology}, {\bf  32}(1993), 827--844.

\bibitem{Ruan} Y. Ruan,  Topological sigma model and
Donaldson-type invariants in Gromov theory. {\it Duke Math. J.,},
{\bf 83}(1996),  461--500.
            \bibitem{Ruan-Li} Y. Ruan and A.-M.
            Li,
             Symplectic surgery and Gromov-Witten invariants of Calabi-Yau 3-folds
             I, preprint 1999
    \bibitem{Ruan-virtual} Y. Ruan,  Virtual
    neighborhoods and pseudo-holomorphic curves, {\it Proceedings
    of 6th G\"okova Geometry-Topology Conference, Turkish J. Math.}, {\bf  23} (1999),
    161--231.
\bibitem{Ruan-Tian} Y. Ruan and G. Tian, A mathematical theory of
quantum cohomology, {\it J. Diff. Geom.}, {\bf 42}(1995),259--367.
\bibitem{Siebert} B. Siebert,
Gromov-Witten invariants of general symplectic manifolds, preprint 1997.
   \bibitem{Salamon:Floer} D. Salamon. Lectures on Floer homology,
in ``Symplectic Geometry and Topology'', IAS/Park City Mathematics Series, v. 7, AMS/IAS
1999, 144 -- 229.
 \bibitem{Ustilovsky}  I. Ustilovsky, Infinitely many contact
 structures on $S^{4m+1}$, {\it Int. Math. Res. Notices}, {\bf
 14}(1999), 781--791.

 \bibitem{Vakil} R. Vakil, The enumerative geometry of rational
 and elliptic curves in projective space, preprint 1997.
 \bibitem{Viterbo} C. Viterbo, in preparation.
 \bibitem{Weinstein}  A. Weinstein,
 On the hypotheses of Rabinowitz's periodic orbits theorems
{\it J. Diff. Eq. }, {\bf 33}(1979), 353--358.
 \bibitem{Witten1} E. Witten,   Supersymmetry and Morse theory,
 {\it J. Diff. Geom.}, {\bf 17}(1982), 661--692.
 \bibitem{Witten2} E. Witten, Two-dimensional gravity and
 intersection theory on moduli space, {\it Surveys in Diff.
 Geom.}, {\bf 1}(1991), 243--310.
\bibitem{MLYau} M.-L. Yau, Contact homology of subcritical Stein
manifolds, thesis, Stanford University, 1999.



           \end{thebibliography}
\end{document}